\documentclass[review]{elsarticle}
\usepackage{mathrsfs}
\usepackage{amsfonts}
 \usepackage{amsthm}
\usepackage{paralist}
\usepackage{graphics}
 \usepackage{epstopdf}
\usepackage{float}
\usepackage{graphicx}
\usepackage{subcaption}
\usepackage{amsmath}
\usepackage[colorlinks=true]{hyperref}
\usepackage{color}
\hypersetup{urlcolor=blue, citecolor=red}
\usepackage{lineno,hyperref}
\modulolinenumbers[5]
\usepackage[normalem]{ulem}
\usepackage{physics}

\usepackage{filecontents}
\usepackage{pgfplots}
\usepackage[margin=1.5in]{geometry}
\pgfplotsset{compat=newest}
\usetikzlibrary{decorations.markings}
\usepackage{pgfplotstable}
\usepackage{booktabs}
\usepackage[normalem]{ulem}
\usepackage{overpic}
\usepackage[colorlinks=true]{hyperref}
\hypersetup{urlcolor=blue, citecolor=blue}
\graphicspath{{./figures/}}
\usepackage{soul}
\soulregister\cite7
\soulregister\eqref7
\soulregister\autoref7
\usepackage{natbib}

\journal{Journal of \LaTeX\ Templates}

\newtheorem{theorem}{Theorem}[section]
\newtheorem{lemma}{Lemma}[section]
\newtheorem{definition}{Definition}[section]

\newtheorem{proposition}{Proposition}[section]
\newtheorem{remark}{Remark}[section]

\numberwithin{equation}{section}
\numberwithin{figure}{section}

\allowdisplaybreaks[4]

\def\be{\begin{equation}}
\def\ee{\end{equation}}
\def\bea{\begin{eqnarray}}
\def\eea{\end{eqnarray}}
\def\bt{\begin{theorem}}
\def\et{\end{theorem}}
\def\bl{\begin{lemma}}
\def\el{\end{lemma}}
\def\br{\begin{remark}}
\def\er{\end{remark}}
\def\bp{\begin{proposition}}
\def\ep{\end{proposition}}
\def\bc{\begin{corollary}}
\def\ec{\end{corollary}}
\def\bd{\begin{definition}}
\def\ed{\end{definition}}

\def\Pi{\mathbf{\psi}}

\graphicspath{ {./result/} }

\newcommand{\Ra}{\operatorname{Ra}}
\newcommand{\Ta}{\operatorname{Ta}}
\newcommand{\Q}{\operatorname{Q}}
\allowdisplaybreaks

\numberwithin{equation}{section}

\begin{document}

\begin{frontmatter}

\title{Dynamical transition of hydromagnetic convection in a rotating fluid layer}

%


\author[a]{Liang Li}
\address[a]{College of Mathematics, Sichuan University,\\ Chengdu, 610065, P. R. China}
\ead{2019322010018@stu.scu.edu.cn}
\author[b]{Yanlong Fan
}
\address[b]{College of Mathematics, Sichuan University,\\ Chengdu, 610065, P. R. China}
\ead{fanyanlong1@stu.scu.edu.cn}
\author[c]{Quan Wang
\corref{mycorrespondingauthor}
}
\cortext[mycorrespondingauthor]{Corresponding author}
\address[c]{College of Mathematics, Sichuan University,\\ Chengdu, 610065, P. R. China}
\ead{xihujunzi@scu.edu.cn}

\author[d]{Daozhi Han}
\address[d]{Department of Mathematics and Statistics, Missouri University of Science and Technology, Rolla, MO 65409, USA}
\ead{handaoz@mst.edu}



\begin{abstract}
	In this article, we aim to study the stability and dynamic transition of an electrically conducting fluid in the presence of an external uniform horizontal magnetic field and rotation based on a Boussinesq approximation model. By analyzing the spectrum of the linear part of the model and verifying the validity of  the principle of exchange of stability, we take a hybrid approach combining theoretical analysis with numerical computation to study the transition from a simple real eigenvalue, a pair of complex conjugate eigenvalues and a real eigenvalue of multiplicity two, respectively. The center manifold reduction theory is applied to reduce the infinite dimensional system to the corresponding finite dimensional one together with one or several non-dimensional transition numbers that determine the dynamic transition types. Careful numerical computations are performed to determine these transition numbers as well as related temporal and flow patterns etc. Our results indicate that both continuous and jump transitions can occur at certain parameter region. 
\end{abstract}


\end{frontmatter}

\section{Introduction}
The Rayleigh-B\'{e}nard (RB) convection is a classical buoyancy-driven convection phenomenon which describes the motion of a horizontal fluid layer heated from below and cooled on the top \cite{Getling1998,Rayleigh}. Due to the important role in the heat transfer in a thermal system, such as the formation of the large scale ocean circulation \cite{ma2010},
 the general atmospheric circulation \cite{Kieu2019},  movement within the Earth’s mantle \cite{Zhang1999}, and the complex activity in the sun \cite{Hanasoge2015}, problems pertaining to the RB convection have been of great interest for some time in fields ranging from numerical analysis to experimental physics \cite{Berge1984,Manneville2006}.

The study relevant for the RB convection without a magnetic effect and in a fixed reference frame has attracted in the past as also nowadays the attention of so many researchers. Lord Rayleigh was the first one to solve the problem of the onset of thermal convection in a plane horizontal layer of fluid heated from below by developing a linear theory \cite{Rayleigh} in 1916. The method he used became a paradigm for analyzing the linear stabilities, i.e. by employing the separation of variables, the characteristic equation of the linear part at a basic state can be determined by the corresponding boundary conditions, then the basic state is stable if all the eigenvalues' real parts are negative. For more introduction on linear stabilities of problems involved the RB convection, we refer readers to \cite{Jeffreys1926,1940a,Chandrasekhar1961}. Besides, there are many works on nonlinear analysis associated with RB convection at the onset of instability. Ma and Wang showed that the attractor bifurcation of RB problem under physically sound boundary conditions occurs when the Rayleigh number crosses a critical threshold \cite{Ma2003,cms/1109706533} by employing their own attractor bifurcation theorem  \cite{2007,Ma1}. For more results in this direction can be found in \cite{Chen,Kevrekidis1989, Kirchgaessner1975}.

It is widely accepted that the fluid movements in the Earth's core are driven by buoyancy forces and strongly affected by the Lorentz (magnetic field) and Coriolis (rotation) forces. Thus, many researchers have investigated the RB convection in the presence of magnetic field and/or rotation in order to make the conclusion applied to practical situations better. An infinitely extended thin horizontal layer of electrically conduction fluid heated uniformly from below under the rotation and magnetic field was firstly theoretically studied by Chandrasekahar \cite{Chandrasekhar1961}. He also obtained the critical Rayleigh number and wave number for the onset of overstability under different values of Ta. Eltayeb \cite{Eltayeb1975} investigated the overstable convection in various directions of magnetic field by performing detailed asymptotic analysis, and obtained the law for the onset of overstability under some circumstances. Refer to the literature \cite{Roberts2000,Jones2000,Podvigina2009} for the linear theory applied to the related model. Recently, many studies on the related model focus on numerical stimulation: Ghosh and Pal \cite{Ghosh2017} studied the instabilities and chaos of RB convection of electrically conducting fluids; Küker and Rüdiger \cite{Kueker2018} studied the turbulent pressure of magnetoconvection for slow and rapid rotation and demonstrated that the effect of the turbulence differed between the low-conductivity and the high-conductivity; Filippi et al. \cite{Filippi2019} considered the effects of anisotropic diffusion acting on rotating magnetoconvection (RMC) in a plane layer; Ghosh and Pal \cite{Ghosh2018,Banerjee2020} did not only directly simulated the RMC, but also used the numerical simulation to analyze the related bifurcation problems by relying on a low dimensional model. To sum up, the RMC models are mainly investigated by the linear method, or the method of numerical computation. 

There are some works associated with RB convection problem in presence of magnetic field or rotation from the perspective of dynamic transition. Wang and Sengul \cite{2017Pattern} studied the dynamic transition of the incompressible MHD equations in absence of rotation and with a large magnetic Prandtl number. Hisa
et al. \cite{Hsia2007} investigated the dynamic transition of the stratified rotating Boussinesq equations in absence of magnetic field. It is worth noting that, in general, the magnetic Prandtl number of electrically conducting fluids (liquid metals, fluids present in the inner core of earth) is very small ($P_{m}\approx 10^{-6}$). Thus, it is theoretically and practically
necessary to study the dynamic transition of the RMC by applying the model in the zero magnetic Prandtl number limit \cite{Banerjee2020}.

Our goal is to study the dynamical properties of the 
RMC by applying the model established in \cite{Banerjee2020} from the perspective of the theory of phase transition dynamics, which was established by Ma and Wang \cite{Ma1}. According to the theory, the dynamic transitions of all dissipative systems are classified into three types: continuous, jump and mixed. In a word, when the control parameter of a dissipative system passes a critical value, the continuous type indicates that its state gradually changes from a state to another; the jump type means the system abruptly jumps to another state and the mixed one means both continuous type and jump type are possible. The key step in the application of the phase transition dynamic theory is to reduce a high dimensional dissipative system into a system of ODEs with a low dimension by employing the center manifold reduction. The theory has been successful applied to many problems arising in nonlinear science, such as biology\cite{Mao2,Xing2021,Jia2021}, chemistry\cite{Ma3}, physics \cite{quan20202,quan20203,Hs,HHW2019} and so on\cite{quan20201,HC,DH}.

In this paper, our theory analysis is consist of two parts. First, we consider the linear stability of the trivial steady-state of the RMC model, and establish the principle of exchange of stability (PES) condition which is sufficient for the emergence of  dynamic transition. Second, rigorous nonlinear analysis is conducted to analyze the detailed dynamics of the transition from first real simple eigenvalue, first complex  eigenvalues and first real eigenvalue with multiplicity two. Although there are the transitions involved the first eigenvalues with higher multiplicity, we do not consider the non-general case in the present work. We also derive the conditions determining the  types of the dynamic transition, by which and combining careful numerical computations indicate that both continuous and jump transitions can occur for the three cases at certain parameter region.

The rest of the paper is arranged as follows.  
We provide the mathematical formulation of the RMC problem, investigate the linear stability and establish the PES condition in Section \ref{section2}. Section \ref{section3}-\ref{section5} states and proves the main conclusions on the transitions from a first simple real eigenvalue, first complex conjugate eigenvalues and a real eigenvalue of multiplicity two. Section \ref{section6} lists some numerical examples. We end the article with a summary of our main results and open problems in Section \ref{section7}.

\section{Mathematic formulation}\label{section2}
\subsection{The mathematic model}
In this article, we consider the RB convection of rotating electrically conducting fluid flows subject to a uniform external magnetic field. 
The model used for the RMC is the
 dimensionless Boussinesq equations in the zero magnetic Prandtl number limit  \cite{Banerjee2020} : 
\begin{align}\label{M2}
    \begin{cases}
        \frac{\partial \mathbf{u}}{\partial t}=\Delta \mathbf{u}-\nabla{\pi}+\Ra\theta\mathbf{e_{3}}+\sqrt{\Ta }(\mathbf{u}\times\mathbf{e_{3}})+\Q \frac{\partial \mathbf{b}}{\partial x_{2}}-(\mathbf{u}\cdot\nabla)\mathbf{u},
        \\
        \frac{\partial \theta}{\partial t}=\frac{1}{\Pr }\Delta\theta+\frac{u_{3}}{\Pr }-(\mathbf{u}\cdot\nabla)\theta,
        \\
        \Delta\mathbf{b}=-\frac{\partial\mathbf{u}}{\partial x_{2}},
        \\
        \nabla\cdot\mathbf{u}=\nabla\cdot\mathbf{b}=0,
    \end{cases}
\end{align}
where $\mathbf{u}(\mathbf{x },t)=(u_{1}(\mathbf{x },t),u_{2}(\mathbf{x },t),u_{3}(\mathbf{x },t))$ is the velocity field, $\theta(\mathbf{x },t)$ is the temperature field, $\pi(\mathbf{x },t)$ is the modified pressure field, $\mathbf{b}(\mathbf{x },t)=(b_{1}(\mathbf{x },t),b_{2}(\mathbf{x },t),b_{3}(\mathbf{x },t))$ is the induced magnetic field, 
$\mathbf{x }=(x_{1},x_{2},x_{3})$ in the non-dimensional region $\Omega=\mathbb{R}^2\times (0,1),$
and $\mathbf{e_{3}}=(0,0,1)$ is antiparallel to the gravitational acceleration $\mathbf{g}=(0,0,-g)$. The four dimensionless numbers $\Ra$, $\Ta $,
$\Q $ and $\Pr $ are the Rayleigh number, the Taylor number, the Chandrasekhar number and the Prandtl number, respectively.

At the top and bottom boundaries, we consider the following boundary conditions:
\begin{align}\label{boundary1}
    (\theta,u_{3},b_{3})=(0,0,0), ~~~~~~\frac{\partial b_{1}}{\partial x_{3}}=\frac{\partial b_{2}}{\partial x_{3}}=\frac{\partial u_{1}}{\partial x_{3}}=\frac{\partial u_{2}}{\partial x_{3}}=0 ~~~~ \text{at} ~~x_{3}=0,1.
\end{align}
To guarantee that the fluid flows rises or descends along the rotating axis, the following constraint is also adopted
\begin{align}\label{boundary2}
    u_{1}(0,0,x_{3})=u_{2}(0,0,x_{3})=0.
\end{align}
This is a natural condition in many fluid applications. Physically, this condition means that the rotation center lies at the original point $\mathbf{x } = 0$.
Finally, for mathematical convenience, periodic boundary conditions are imposed in the $x_{1}$ and $x_{2}$ directions, i.e.,
\begin{align}\label{boundary3}
    (\mathbf{u},\theta,\mathbf{b})(x_{1},x_{2},x_{3})=(\mathbf{u},\theta,\mathbf{b})(x_{1}+2jL_{1}\pi,x_{2}+2kL_{2}\pi,x_{3})~~~~\forall j,k\in \mathbb{Z}.
\end{align}

In what follows we formulate the evolution equations given in \eqref{M2} by using an abstract functional setting that is standard in the framework of dynamic transitions.
First, we let $H^{2}(\Omega)$, and $L^2(\Omega)$ denote the usual Sobolev and Lebesgue spaces.
Then, denoting $\Psi=(\mathbf{u},\theta)$, we let $\mathbf{\widetilde{H}} $, $\mathbf{H}$ and $ \mathbf{H_{1}} $ be the Hilbert spaces given by
\begin{align}\label{space1}
     & \mathbf{\widetilde{H}} = \{\Psi\in [L^{2}(\Omega)]^4\big{|}(u_{3},\theta)|_{x_{3}=0,1}=0,~\nabla\cdot\mathbf{u}=0,~ \Psi~\text{satisfies}~ \eqref{boundary3} \},
    \\   
     & \mathbf{H}= \{\Psi\in  \mathbf{\widetilde{H}}\big{|}\Psi(-x_{1},-x_{2},x_{3})=(-u_1(\mathbf{x}),-u_2(\mathbf{x}),u_3(\mathbf{x}),\theta)\},
    \\
     & \mathbf{H_{1}} =\{\Psi\in [H^2]^4\cap \mathbf{H}\big{|}(\mathbf{u},\theta)~\text{satisfies} ~\eqref{boundary1}-\eqref{boundary3}\},
\end{align}
which are endowed with their natural inner products. Finally, we introduce the linear operator $L:\mathbf{H_{1}}\rightarrow\mathbf{H} $ and bilinear operator $G:\mathbf{H_{1}}\rightarrow\mathbf{H}$ as follows
\begin{align}
     & L\Psi =\begin{pmatrix}
       \mathbf{P}\Delta\mathbf{u}+
        \mathbf{P}\left(\Ra\theta\mathbf{e_{3}}+\sqrt{\Ta }(\mathbf{u}\times\mathbf{e_{3}})-\Q \frac{\partial \Delta^{-1}\frac{\partial u}{\partial x_{2}}}{\partial x_{2}}\right)
        \\   \frac{1}{\Pr }\Delta\theta+
        \frac{1}{\Pr }u_{3}
    \end{pmatrix},
    \\ 
    \label{bilinear}
     & G(\Psi,\widetilde{\Psi}) =-\begin{pmatrix}
        \mathbf{P}(\mathbf{u}\cdot\nabla)\widetilde{\mathbf{u}}
        \\
        (\mathbf{u}\cdot\nabla)\widetilde{\theta}
    \end{pmatrix},
\end{align}
where $\mathbf{P}:[L^2(\Omega)]^3\rightarrow [L^2(\Omega)]^3$ is the Leray projection and $\Delta^{-1}$ is the inverse of  $\Delta$. 

With the help of the preceding two operators, the system \eqref{M2} equipped with these boundary conditions \eqref{boundary1}-\eqref{boundary3} can be rewritten into the following abstract equation
\begin{align}\label{abstract1}
    \frac{d\Psi}{dt}=L\Psi+\mathbf{G}(\Psi),\quad 
    \Psi\in \mathbf{H_{1}},
\end{align}
in which the nonlinear operator $\mathbf{G}$ acting on $\Psi$ is defined by
\begin{align}\label{nonlinear-part}
\mathbf{G}(\Psi)=G(\Psi,\Psi),
\end{align}
which will be used many times later.

\subsection{The eigenvalue problem}

Consider the following eigenvalue problem:
\begin{align}\label{Eigenvalue}
    L\Psi=\beta\Psi,~~~~\Psi\in \mathbf{H_{1}},
\end{align}
that is 
\begin{align}\label{Eigenvalue1}
    \begin{cases}
        \Delta \mathbf{u}-\nabla{\pi}+\Ra\theta\mathbf{e_{3}}+\sqrt{\Ta }(\mathbf{u}\times\mathbf{e_{3}})+\Q \frac{\partial \mathbf{b}}{\partial x_{2}}=\beta \mathbf{u},
        \\
        \frac{1}{\Pr }\Delta\theta+\frac{u_{3}}{\Pr }=\beta\mathbf{\theta},
        \\
        \Delta\mathbf{b}=-\frac{\partial\mathbf{u}}{\partial x_{2}},
        \\
        \nabla\cdot\mathbf{u}=\nabla\cdot\mathbf{b}=0,
    \end{cases}
\end{align}
subject to the same boundary conditions as \eqref{boundary1}-\eqref{boundary3}.

 With regard to the boundary condition \eqref{boundary3},  we can take the following Fourier series:
\begin{align}\label{Fourier1}
\begin{aligned}
   & \Psi_{J}=\begin{pmatrix}
        u_{J}^1(x_{3})\sin{(j\alpha_{1}x_{1}+k\alpha_{2}x_{2})}
        \\
        u_{J}^2(x_{3})\sin{(j\alpha_{1}x_{1}+k\alpha_{2}x_{2})}
        \\
        u_{J}^3(x_{3})\cos{(j\alpha_{1}x_{1}+k\alpha_{2}x_{2})}
        \\
        \theta_{J}(x_{3})\cos{(j\alpha_{1}x_{1}+k\alpha_{2}x_{2})}
    \end{pmatrix},\\
  &  \begin{pmatrix}
        b_{J}
        \\
        \pi_{J}
    \end{pmatrix}=\begin{pmatrix}
        b_{J}^1(x_{3})\cos{(j\alpha_{1}x_{1}+k\alpha_{2}x_{2})}
        \\
        b_{J}^2(x_{3})\cos{(j\alpha_{1}x_{1}+k\alpha_{2}x_{2})}
        \\
        b_{J}^3(x_{3})\sin{(j\alpha_{1}x_{1}+k\alpha_{2}x_{2})}
        \\
        p_{J}(x_{3})\cos{(j\alpha_{1}x_{1}+k\alpha_{2}x_{2})}
    \end{pmatrix},
 \end{aligned}
\end{align}
where $\alpha_{i}=1/{L_{i}}~(i=1,2)$ and $J$ is the index. Plugging \eqref{Fourier1} into \eqref{Eigenvalue1}, we derive
\begin{align}\label{ODE}
    \begin{cases}
        (D^2-\alpha_{jk}^{2}-\beta)u_{J}^{1}(x_{3})+\sqrt{\Ta }u_{J}^{2}(x_{3})-\Q k\alpha_{2}b_{J}^{1}(x_{3})+j\alpha_{1}p_{J}(x_{3})=0,
        \\
        (D^2-\alpha_{jk}^2-\beta)u_{J}^{2}(x_{3})-\sqrt{\Ta }u_{J}^{1}(x_{3})-\Q k\alpha_{2}b_{J}^{2}(x_{3})+k\alpha_{2}p_{J}(x_{3})=0,
        \\
        (D^2-\alpha_{jk}^2-\beta)u_{J}^{3}(x_{3})+\Ra\theta_{J}(x_{3})+\Q k\alpha_{2}b_{J}^{3}(x_{3})-Dp_{J}(x_{3})=0,
        \\
        (D^2-\alpha_{jk}^2-\Pr \beta)\theta_{J}(x_{3})=-u_{J}^3(x_{3}),
        \\
        (D^2-\alpha_{jk}^2)b_{J}^{1}(x_{3})=-k\alpha_{2}u_{J}^{1}(x_{3}),
        \\
        (D^2-\alpha_{jk}^2)b_{J}^{2}(x_{3})=-k\alpha_{2}u_{J}^{2}(x_{3}),
        \\
        (D^2-\alpha_{jk}^2)b_{J}^{3}(x_{3})=k\alpha_{2}u_{J}^{3}(x_{3}),
        \\
        j\alpha_{1}u_{J}^{1}(x_{3})+k\alpha_{2}u_{J}^2(x_{3})+Du_{J}^{3}(x_{3})=0,
    \end{cases}
\end{align}
subject to the following boundary condition
\begin{align}\label{ODE-bdc}
    \begin{cases}
u_{J}^3(x_{3})=b_{J}^3(x_{3})=\theta_{J}(x_{3})=0,\quad x_{3}=0,1,\\
Du_{J}^1(x_{3})=Du_{J}^2(x_{3})=Db_{J}^1(x_{3})=Db_{J}^2(x_{3})=0,
\quad x_{3}=0,1,
    \end{cases}
\end{align}
where $D=d/{dx_{3}}$ and $\alpha_{jk}^2=(j\alpha_{1})^{2}+(k\alpha_{2})^2$. 

 One can deduce from \eqref{ODE} and \eqref{ODE-bdc} that 
\begin{align}
    D^{2}u_{J}^{3}(x_{3})=\cdots=D^{2n}u_{J}^{3}(x_{3})=\cdots=0 ~(n\in \mathbb{N})~~~\text{at}~~x_{3}=0,1,
\end{align}
which means $u_{J}^{3}(x_{3})=\sin{l\pi}x_{3}$ for each $l\in \mathbb{Z}^{+}$.
Thus, denote $r_{J}^2=\alpha_{jk}^2+l^2\pi^2$ with index $J=(j,k,l)\in \mathbb{Z}^2\times \mathbb{N}$, where $(j,k)\neq(0,0)$ and $ l \neq 0 $, then $\beta$ solves
\begin{align}\label{characteristicequation1}
    a_{3}\beta^{3}+a_{2}\beta^{2}+a_{1}\beta+a_{0}=0,
\end{align}
where
\begin{equation}\label{characteristicequation}
    \begin{aligned}
        a_3 & =   \Pr r_J ^4,\quad                                                                                                                                   
        a_2  =   \left(r_J ^4+2 \Pr  \left(r_J ^4+ \Q k^2 \alpha_2^2 \right)\right)r_J ^2;                                                                           \\
        a_1 & =  \left(r_J ^4+ \Q k^2\alpha_2^2  \right) \left(\Pr  \Q k^2\alpha_2^2  +r_J ^4 (\Pr +2)\right)+\Pr  \Ta  l^2 \pi^2   r_J ^2 - \Ra \alpha_{jk}^2 r_J ^2 ; \\
        a_0 & =  \left(r_J ^5+ \Q k^2\alpha_2^2    r_J \right)^2+r_J ^4   l^2\pi ^2  \Ta - \Ra \alpha_{jk}^2 \left(r_J ^4+Qk^2\alpha_2^2  \right) .
    \end{aligned}
\end{equation}

To derive the explicit expressions of each eigenvalue $\beta$ and eigenvector $\Psi$
of $L$,
we introduce the following index sets:
\begin{align*}
    I_{1}= {} & \{(j,k,l)|(j,k)\in \mathbb{Z}^2, j\geq  0,~(j,k)\neq(0,0),~l=1,2,\cdots\},
    \\
    I_{2}= {} & \{(j,k,0)|(j,k)\in \mathbb{Z}^2, j\geq  0,~(j,k)\neq(0,0)\},
    \\
    I_{3}= {} & \{(0,0,l)|l=1,2,\cdots\}.
\end{align*}

First, for each $J\in I_{3}$,  according to \eqref{Fourier1}, we let
\begin{align}\label{I3}
    \Psi_{J}=\begin{pmatrix}
        u_{J}^1\cos{l\pi x_{3}}
        \\
        u_{J}^2\cos{l\pi x_{3}}
        \\
        u_{J}^3\sin{l\pi x_{3}}
        \\
        \theta_{J}\sin{l\pi x_{3}}
    \end{pmatrix},~~
    \begin{pmatrix}
        b_{J}
        \\
        \pi_{J}
    \end{pmatrix}=\begin{pmatrix}
        b_{J}^1\cos{l\pi x_{3}}
        \\
        b_{J}^2\cos{l\pi x_{3}}
        \\
        b_{J}^3\sin{l\pi x_{3}}
        \\
        p_{J}\cos{l\pi x_{3}}
    \end{pmatrix}
\end{align}
where $u_{J}^{s},b_{J}^{s} (s=1,2,3)$, $p_{J}$ and $\theta_{J}$ are coefficients.  Then, it derives from \eqref{Eigenvalue1} that
\begin{align}\label{I31}
    \begin{cases}
        (l^2\pi^2+\beta)u_{J}^1-\sqrt{\Ta }u_{J}^2=0,
        \\
        (l^2\pi^2+\beta)u_{J}^2+\sqrt{\Ta }u_{J}^1=0,
    \end{cases}
\end{align}
which has nontrivial solution $(u_{J}^1,u_{J}^2)$ if and only if 
$\beta= \beta_{J}^{k}~(k=1,2)$, where
\begin{align}
    \beta_{J}^{1}=-l^2\pi^2+i\sqrt{\Ta },\quad   \beta_{J}^{2}=-l^2\pi^2-i\sqrt{\Ta },\quad J \in  I_3.
\end{align}
If \eqref{I31} has a nontrivial solution $(u_{J}^1,u_{J}^2)$, it follows from $\nabla\cdot \mathbf{u}=0$, $u_{J}^{3}=0$ and the second equation of \eqref{Eigenvalue1}  that  $\theta_{J}=0$. Taking the nontrivial solution $(u_{J}^1,u_{J}^2)=(1,\pm i)$, we have a pair of conjugate eigenvectors $ \Psi_{J}^{1}\pm  i  \Psi_{J}^{2}$, where 
\begin{align}\label{eigenvectorI3}
    \Psi_{J}^{1}=\begin{pmatrix}
        \cos{l\pi x_{3}}
        \\
        0
        \\
        0
        \\
        0
    \end{pmatrix},\quad 
    \Psi_{J}^{2}=\begin{pmatrix}
        0
        \\
        \cos{l\pi x_{3}}
        \\
        0
        \\
        0
    \end{pmatrix}.
\end{align}
If \eqref{I31} has the unique solution $(u_{J}^1,u_{J}^2)=(0,0)$, it follows from $\nabla\cdot \mathbf{u}=0$ that $u_{J}^{3}=0$, and we can obtain the corresponding eigenvalue $\beta_{J}^{3}$ and eigenvector $  \Psi_{J}^{3} $ by solving the second equation of \eqref{Eigenvalue1}, which are given by
\begin{align}
    \beta_{J}^{3} = -\frac{l^2\pi^2}{\Pr },\quad
    \Psi_{J}^{3}   = \begin{pmatrix}
        0
        \\
        0
        \\
        0
        \\
        \sin{l\pi x_{3}}
    \end{pmatrix}.
\end{align}
Hence for the index set $ I_{3}$, we introduce the corresponding eigen-space 
\begin{align}\label{EigenvectorsI31}
 \mathbf{ E_{3}}=\operatorname*{span}_{J\in I_3, s\in \{1,2,3\} }{\Psi_{J}^{s}}.
\end{align}

Next, for each $J\in I_{2}$, \eqref{boundary3} allows us to set
\begin{align}\label{EigenvectorI2}
    \Psi_{J}=\begin{pmatrix}
        u_{J}^{1}\sin{(j\alpha_{1}x_{1}+k\alpha_{2}x_{2})}
        \\
        u_{J}^{2}\sin{(j\alpha_{1}x_{1}+k\alpha_{2}x_{2})}
        \\
        0
        \\
        0
    \end{pmatrix}.
\end{align}
Inserting the preceding expression into \eqref{ODE} and
performing some simple calculations, the corresponding eigenvalue and eigenvector are given by
\begin{align}\label{EigenvalueI2}
     & \beta_{J}=-\frac{\Q k^2\alpha_{2}^2+\alpha_{jk}^{4}}{\alpha_{jk}^2},
         \\
     & \Psi_{J}=\begin{pmatrix}
        k\alpha_{2}\sin{(j\alpha_{1}x_{1}+k\alpha_{2}x_{2})}
        \\
        -j\alpha_{1}\sin{(j\alpha_{1}x_{1}+k\alpha_{2}x_{2})}
        \\
        0
        \\
        0
    \end{pmatrix}.
\end{align}
Likewise, for the index set $ I_{2}$, we introduce the corresponding eigen-space  
\begin{align}\label{EigevalueI21}
    \mathbf{ E_{2}}=\operatorname*{span}_{J\in I_2}\{\Psi_{J}\}.
\end{align}

Finally, we consider eigenvalue $\beta_J$ with $J\in I_1$.
Let three roots of the cubic equation \eqref{characteristicequation1} be ordered by $\Re \beta_{J}^{1}\geq \Re \beta_{J}^{2}\geq \Re \beta_{J}^{3}$. For the real eigenvalue $\beta_{J}^{s} $,  one can get the corresponding eigenvector given by
\begin{align}\label{EigenvectorsI1}
    \Psi_{J}^{s}=\begin{pmatrix}
        u_{J}^{1s}\sin{(j\alpha_{1}x_{1}+k\alpha_{2}x_{2})}\cos{l\pi x_{3}}
        \\
        u_{J}^{2s}\sin{(j\alpha_{1}x_{1}+k\alpha_{2}x_{2})}\cos{l\pi x_{3}}
        \\
        \cos{(j\alpha_{1}x_{1}+k\alpha_{2}x_{2})}\sin{l\pi x_{3}}
        \\
        \theta_{J}^{s}\cos{(j\alpha_{1}x_{1}+k\alpha_{2}x_{2})}\sin{l\pi x_{3}}
    \end{pmatrix},
\end{align}
with
\begin{equation}
    \begin{aligned}
        u^{1s}_{J}     & = -\frac{l\pi}{\alpha_{jk}^2}\left( j \alpha_1+  \frac{\sqrt{\Ta} k \alpha_2 r_{J}^{2} }{ (\beta_{J}^{s} +r_{J}^{2}) r_{J}^{2}+\Q k^2\alpha_2^2  } \right) ,  \\
        u^{2s}_{J}     & = -\frac{l\pi}{\alpha_{jk}^2}\left( k \alpha_2 -  \frac{\sqrt{\Ta} j \alpha_1 r_{J}^{2} }{ (\beta_{J}^{s} +r_{J}^{2}) r_{J}^{2}+\Q k^2\alpha_2^2  } \right) , \quad \theta^{s}_{J} & = \frac{1}{r_{J}^{2}+\Pr \beta_{J}^{s}}.
    \end{aligned}
\end{equation}
If \eqref{characteristicequation1} has a pair of complex conjugate roots, $\beta^{q}_{J}=\overline{\beta^{q+1}_{J}}$($q=1$ or $2$), correspondingly,
$l$ has a pair of conjugate eigenvectors $ \Psi_{J}^{q}\pm  i  \Psi_{J}^{q+1}$, where 
\begin{align}\label{EigenvectorI11}
   & \Psi_{J}^{q}=\Re \begin{pmatrix}
        u_{J}^{1q}\sin{(j\alpha_{1}x_{1}+k\alpha_{2}x_{2})}\cos{l\pi x_{3}}
        \\
        u_{J}^{2q}\sin{(j\alpha_{1}x_{1}+k\alpha_{2}x_{2})}\cos{l\pi x_{3}}
        \\
        \cos{(j\alpha_{1}x_{1}+k\alpha_{2}x_{2})}\sin{l\pi x_{3}}
        \\
        \theta_{J}^{q}\cos{(j\alpha_{1}x_{1}+k\alpha_{2}x_{2})}\sin{l\pi x_{3}}
    \end{pmatrix},\\
    \label{EigencectorI122}
   & \Psi_{J}^{q+1}=\Im \begin{pmatrix}
        u_{J}^{1q}\sin{(j\alpha_{1}x_{1}+k\alpha_{2}x_{2})}\cos{l\pi x_{3}}
        \\
        u_{J}^{2q}\sin{(j\alpha_{1}x_{1}+k\alpha_{2}x_{2})}\cos{l\pi x_{3}}
        \\
        \cos{(j\alpha_{1}x_{1}+k\alpha_{2}x_{2})}\sin{l\pi x_{3}}
        \\
        \theta_{J}^{q}\cos{(j\alpha_{1}x_{1}+k\alpha_{2}x_{2})}\sin{l\pi x_{3}}
    \end{pmatrix},
\end{align}
which satisfy $L\Psi_{J}^{q}=\Re \beta_{J}^{q}\Psi_{J}^{q}-\Im \beta_{J}^{q}\Psi_{J}^{q+1}$, and $L\Psi_{J}^{q+1}=\Im \beta_{J}^{q}\Psi_{J}^{q}+\Re \beta_{J}^{q}\Psi_{J}^{q+1}$.
Accordingly, for the index set $ I_{1}$, 
 we set the corresponding eigen-space 
\begin{align}\label{EigenvectorI111}
   \mathbf{ E_{1}}=\operatorname*{span}_{J\in I_{1},s\in \{ 1,2,3, \}} \{\Psi_{J}^{s}\}.
\end{align}

Note that $L$ is a completely field, guaranteeing that $\mathbf{H_{1}}$
can be decomposed into 
\begin{align}\label{decomposition}
\mathbf{H_{1}=E_{1}\oplus E_{2}\oplus  E_{3}}. 
\end{align}

In order to reduce the system \eqref{abstract1} to a system of ODEs by center manifold reduction, we need to derive eigenvectors of the conjugate operator $L^{*}$. 
Making use of the standard definition of the conjugate operator $L^{*}$, i.e.,
\begin{align}\label{conjugateformual}
    \left\langle L\Psi,\Psi^{*}\right\rangle=\left\langle \Psi, L^{*}\Psi^{*}\right\rangle,\quad \forall\ \Psi=(\mathbf{u},\theta),\quad \Psi^{*}=(\mathbf{u^{*}},\theta^{*})\in \mathbf{H_{1}},
\end{align}
where $\left\langle \cdot, \cdot\right\rangle $ denoting the inner product in $L^2$, it gives that  $L^{*}$ takes the form:
\begin{align}\label{conjugateoperator}
    L^{*}\Psi^{*}=\begin{pmatrix}
        \mathbf{P} [\Delta \mathbf{u^{*}}-\sqrt{\Ta }(\mathbf{u^{*}}\times \mathbf{e_{3}})+\frac{1}{\Pr }\theta^{*}\mathbf{e_{3}}-\Q  \frac{\partial \Delta^{-1}\frac{\partial u^{*}}{\partial x_{2}}}{\partial x_{2}}]
        \\
        \frac{1}{\Pr }\Delta\theta^{*}+\Ra u_{3}^{*}
    \end{pmatrix}.
\end{align}

Similarly, we obtain the conjugate eigenvalues with index $J\in I_{3}$ as follows, 
\begin{align}\label{conjugateeigenvaluesI3}
    \beta_{J}^{1*}=\overline{\beta_{J}^{2*}}=-l^2\pi^2-i\sqrt{\Ta }, \quad \beta_{J}^{3*}=-\frac{l^2\pi^2}{\Pr },
\end{align}
and the corresponding eigenvectors $\Psi_{J}^{1*}\pm i   \Psi_{J}^{2*}$
and  $\Psi_{J}^{3*}$ are given by
\begin{align}\label{conjugateeigenvectorsI3}
    \Psi_{J}^{1*}=\begin{pmatrix}
        \cos{l\pi x_{3}}
        \\
        0
        \\
        0
        \\
        0
    \end{pmatrix},\quad 
    \Psi_{J}^{2*}=\begin{pmatrix}
        0
        \\
        \cos{l\pi x_{3}}
        \\
        0
        \\
        0
    \end{pmatrix},\quad 
    \Psi_{J}^{3*}=\begin{pmatrix}
        0
        \\
        0
        \\
        0
        \\
        \sin{l\pi x_{3}}
    \end{pmatrix}.
\end{align}
For $J\in I_{2}$, we have:
\begin{align}\label{conjugateeigenvalueI2}
    \beta_{J}^{*}=-\frac{\Q k^2\alpha_{2}^{2}+\alpha_{jk}^{4}}{\alpha_{jk}^{2}},\quad \Psi_{J}^{*}=\begin{pmatrix}
        k\alpha_{2}\sin{(j\alpha_{1}x_{1}+k\alpha_{2}x_{2})}
        \\
        -j\alpha_{1}\sin{(j\alpha_{1}x_{1}+k\alpha_{2}x_{2})}
        \\
        0
        \\
        0
    \end{pmatrix}.
\end{align}
And for $J\in I_{1}$, first assuming $\Re \beta_{J}^{1*}\geq \Re \beta_{J}^{2*}\geq \Re \beta_{J}^{3*}$, and if $\beta_{J}^{s*}$ is real, then the corresponding eigenvector is:
\begin{align}\label{conjugateeigenvectorsI1}
    \Psi_{J}^{s*}=\begin{pmatrix}
        u_{J}^{1s*}\sin{(j\alpha_{1}x_{1}+k\alpha_{2}x_{2})}\cos{l\pi x_{3}}
        \\
        u_{J}^{2s*}\sin{(j\alpha_{1}x_{1}+k\alpha_{2}x_{2})}\cos{l\pi x_{3}}
        \\
        \cos{(j\alpha_{1}x_{1}+k\alpha_{2}x_{2})}\sin{l\pi x_{3}}
        \\
        \theta_{J}^{s*}\cos{(j\alpha_{1}x_{1}+k\alpha_{2}x_{2})}\sin{l\pi x_{3}}
    \end{pmatrix},
\end{align}
where 
\begin{equation}
    \begin{aligned}
        u_{J}^{1s*}     & = - \frac{l \pi  }{\alpha_{jk} ^2} \left( - \frac{\sqrt{\Ta} k\alpha_2 r_J ^2  }{ r_{J}^{2}(r_{J}^{2}+\beta_{J}^{s*})+\Q k^2\alpha_{2}^2 } + j\alpha_1 \right) , \\ 
        u_{J}^{2s*}     & = -\frac{l \pi   }{\alpha_{jk} ^2} \left(\frac{\sqrt{\Ta} j\alpha_1 r_J ^2  }{ r_{J}^{2}(r_{J}^{2}+\beta_{J}^{s*})+ \Q k^2\alpha_{2}^2} + k\alpha_2 \right)  , \quad \theta_{J}^{s*} & =  \frac{\Pr  \Ra }{ \beta_{J}^{s*} \Pr +r_J^{2}} .
    \end{aligned}
\end{equation}
If there is a pair of complex conjugate roots, $\beta^{q*}_{J}=\overline{\beta^{(q+1)*}_{J}}$($q=1$ or $2$), then $L^{*}$ has a pair of conjugate eigenvectors $ \Psi_{J}^{q*}\pm  i  \Psi_{J}^{q*}$, where
\begin{align}\label{conjugateeigenvectorI11}
    &\Psi_{J}^{q*}     =\Re \begin{pmatrix}
        u_{J}^{1q*}\sin{(j\alpha_{1}x_{1}+k\alpha_{2}x_{2})}\cos{l\pi x_{3}}
        \\
        u_{J}^{2q*}\sin{(j\alpha_{1}x_{1}+k\alpha_{2}x_{2})}\cos{l\pi x_{3}}
        \\
        \cos{(j\alpha_{1}x_{1}+k\alpha_{2}x_{2})}\sin{l\pi x_{3}}
        \\
        \theta_{J}^{q*}\cos{(j\alpha_{1}x_{1}+k\alpha_{2}x_{2})}\sin{l\pi x_{3}}
    \end{pmatrix},
    \\
    &\Psi_{J}^{(q+1)*}  =\Im \begin{pmatrix}
        u_{J}^{1q*}\sin{(j\alpha_{1}x_{1}+k\alpha_{2}x_{2})}\cos{l\pi x_{3}}
        \\
        u_{J}^{2q*}\sin{(j\alpha_{1}x_{1}+k\alpha_{2}x_{2})}\cos{l\pi x_{3}}
        \\
        \cos{(j\alpha_{1}x_{1}+k\alpha_{2}x_{2})}\sin{l\pi x_{3}}
        \\
        \theta_{J}^{q*}\cos{(j\alpha_{1}x_{1}+k\alpha_{2}x_{2})}\sin{l\pi x_{3}}
    \end{pmatrix},
\end{align}
and
 $L^{*}\Psi_{J}^{q*}=\Re \beta_{J}^{q*}\Psi_{J}^{q}+\Im \beta_{J}^{q*}\Psi_{J}^{(q+1)*}$, and $L^{*}\Psi_{J}^{(q+1)*}=-\Im \beta_{J}^{q*}\Psi_{J}^{q*}+\Re \beta_{J}^{q*}\Psi_{J}^{(q+1)*}$.

 \subsection{Principle of exchange of stabilities}
The first step to study the dynamical transitions of \eqref{abstract1} is to verify its PES condition. From the preceding discussion, we have known that if $J\in I_{2}$ or $I_{3}$, the real parts of the corresponding eigenvalues of $L$ are always negative. Moreover, the real parts of all eigenvalues of $L$ are also negative for $\Ra=0$. Thus, to verify the PES condition, we only need to focus on the characteristic equation \eqref{characteristicequation1} with $J\in I_{1}$ and $\Ra>0$.

 Note that \eqref{characteristicequation1} has a zero root if and only if $a_{0}=0$, or 
 equivalently, $\Ra=f(J)$ where 
 \begin{align}\label{weiledingyizhibiaoji}
    f(J)=\frac{r_{J}^{2}(r_{J}^{4}+\Q  k^2\alpha_{2}^{2})^{2}+ \Ta  r_{J}^{4}l^2\pi^2}{\alpha_{jk}^{2}(r_{J}^{4}+ \Q  k^2\alpha_{2}^{2})}.
\end{align}
And, \eqref{characteristicequation1} has a pair of pure complex conjugate roots if and only if $a_{1}a_{2}=a_{0}a_{3}$ which is equivalent to $\Ra=g(J)$, where 
\begin{align}\label{weiledingyizhibiaoji2}
    g(J)=\frac{2(\Q k^2\alpha_{2}^2  +r_{J}^{4})}{\alpha_{jk}^2 r_{J}^{2}} \left[ (1+ \Pr  )r_{J}^{4}+\Pr \Q k^2\alpha_{2}^{2}  +\frac{\Pr ^{2} \Ta r_{J}^{2}l^2\pi^2 }{(1+ \Pr  )r_{J}^{4}+\Pr \Q  k^2\alpha_{2}^{2} }  \right].
\end{align}

 To search for the threshold $\Ra_c$ of $\Ra$
at which the real part of first eigenvalue of $L$ becomes critical, we introduce two
critical values $\Ra_{c_{1}}$ and $\Ra_{c_{1}}$, defined by
\begin{align}\label{critical1}
    \Ra_{c_{1}}=\min_{ J \in I_{1}}f(J),\quad
    \Ra_{c_{2}}=\min_{ J \in I_{1}} g(J).
\end{align}
Then, the exact expression of the threshold $Ra_{c}$ for the system \eqref{abstract1} can by intuitively derived by $\Ra_{c_{1}}$
and $\Ra_{c_{2}}$, given by
\begin{align}\label{linjiezhibiaobijiao}
    \Ra_{c}=\min{ \{\Ra_{c_{1}},\Ra_{c_{2}}\}}.
\end{align} 

In what follows, let us verify the PES condition for the system \eqref{abstract1}. To this end, we introduce a critical index set defined by
\begin{align}\label{cirticalindexset}
    X =  \{J\in I_{1}|\min{\{f(J),g(J)\}}= \Ra_{c}\}.
\end{align}
Note that $\text{card}(X)$ is finite, because both $f(J)$ and $g(J)$ go to positive infinity as $\abs{J}$ approaches positive infinity. We thus have the following PES condition:
\begin{lemma}\label{compare} For the system \eqref{abstract1},
we have the following assertions:
\begin{enumerate}[\rm{(}1\rm{)}]
\item When $Ra_{c}=Ra_{c_{1}}<Ra_{c_{2}}$, we have
    \begin{align}\label{PES1}
        \begin{aligned}
         & \beta_{J}^{1}
         \begin{cases}
            <~0,~~\Ra<\Ra_{c_{1}}
            \\
            =~0,~~\Ra=\Ra_{c_{1}}
            \\
            >~0,~~\Ra>\Ra_{c_{1}}
        \end{cases}~~~~J\in X,
        \\
         & \Re \beta_{J}^{s}(\Ra_{c_{1}})<0 ~~\forall(J,s)\notin X\times\{1\}(s=1,2,3);
        \end{aligned}
    \end{align}
    \item When $Ra_{c}=Ra_{c_{2}}<Ra_{c_{1}}$, we have
 \begin{align}\label{PES2}
        \begin{aligned}
         & \Re \beta_{J}^{1}=\Re \beta_{J}^{2}
         \begin{cases}
            <~0,~~\Ra<\Ra_{c_{2}}
            \\
            =~0,~~\Ra=\Ra_{c_{2}}
            \\
            >~0,~~\Ra>\Ra_{c_{2}}
        \end{cases}~~~~J\in X,
        \\
         & \Re \beta_{J}^{s}(\Ra_{c_{2}})<0 ~~\forall(J,s)\notin X\times\{1,2\}(s=1,2,3).
    \end{aligned}
    \end{align}

\end{enumerate}
    \end{lemma}
\begin{proof}
    For $\Ra=0,$ we have
    \begin{align*}
        \beta_{J}=-\frac{r_{J}^{2}}{\Pr }, ~-\frac{\Q k_{2}^{2}\alpha_{2}^{2}}{r_{J}^{2}}-r_{J}^{2}\pm\sqrt{\Ta }r_{J}l\pi i,\quad \forall J\in I_{1}.
    \end{align*}
\begin{enumerate}[(1)]
    \item  For $\Ra_{c_{2}}>\Ra_{c_{1}}$, owing to the continuous dependence of $\beta_{J}(\Ra)$ on $\Ra$, and by the definition of $\Ra_{c_{1}}$ there exists $\delta_0>0$ such that
    \begin{enumerate}[(i)]
        \item if $\Ra<\Ra_{c_{1}}$, then $\Re \beta_{J}<0$ $\forall J$;
        \item if $\Ra=\Ra_{c_{1}}$, then $\beta_{J}^{1}=0$ for $J\in X$ and $\beta_{J}^{s}<0$ for $(J,s)\notin X\times\{1\}$;
        \item if $\Ra_{c_{1}}+\delta_0>\Ra>\Ra_{c_{1}}$, the $\beta_{J}^{1}>0$ for $J\in X$ and $\beta_{J}^{s}<0$ for $(J,s)\notin X\times\{1\}$.
    \end{enumerate}
    Therefore, the first conclusion holds.
    \item For $\Ra_{c_{2}}<\Ra_{c_{1}}$, we let $\beta_{J_{1}}^{1}=\overline{\beta_{J_{1}}^{2}}=\sigma(\Ra)+\rho(\Ra)i$.  The definition of $\Ra_{c_{2}}$ says that
    \begin{enumerate}[(i)]
        \item if $\Ra<\Ra_{c_{2}}$, then $\Re \beta_{J}<0$,  $\forall J$.
        \item if $\Ra=\Ra_{c_{2}}$, then $\sigma(\Ra_{c_{2}})=0$.
      \end{enumerate}
      Therefore, to show the second conclusion, we only need to show $     \sigma'(\Ra_{c_{2}})>0$. Taking the derivative of both sides of \eqref{characteristicequation1} with respect to $\Ra$ at $\Ra_{c_{2}}$, it yields                \begin{align}\label{derivative}
                    \begin{cases}
                        -2 a_2 \rho \rho '-2 a_1 \sigma '+a_0' = 0, 
                        \\
                        -2 a_1 \rho '+2 a_2 \rho  \sigma '+\rho  a_1' = 0.
                    \end{cases}
                \end{align}
It obtains from the preceding equations that
                \begin{equation}\label{derivative-2}
                    \begin{aligned}
                        \sigma'(\Ra_{c_{2}})= {} & \frac{a_1  a_0' - a_0 a_1' }{2 a_2^2 \rho  ^2+2 a_1 ^2} > 0 \iff a_1  a_0' - a_0 a_1' > 0.
                    \end{aligned}
                \end{equation}
A direct verification shows that $a_1  a_0' - a_0 a_1' > 0$ is guaranteed by \( \Ra_{c_2} <\Ra_{c_1} \).
\end{enumerate}    
\end{proof}

\begin{figure}[H]
    \begin{minipage}[t]{0.45\linewidth}
        \centering
        {\includegraphics[width=2.5in]{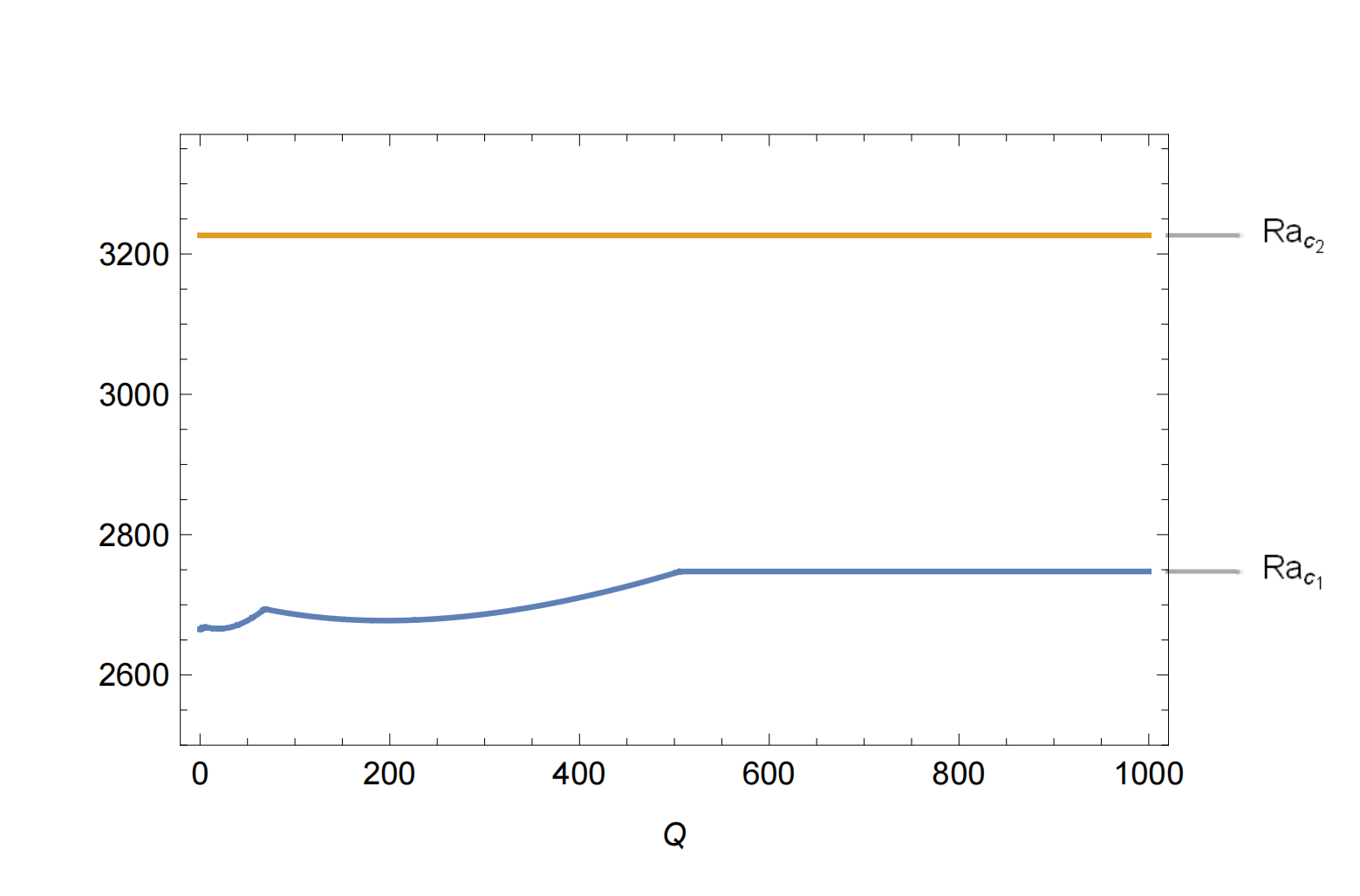}}
    \end{minipage}
    \hfill
    \begin{minipage}[t]{0.45\linewidth}
        \centering
        {\includegraphics[width=2.5in]{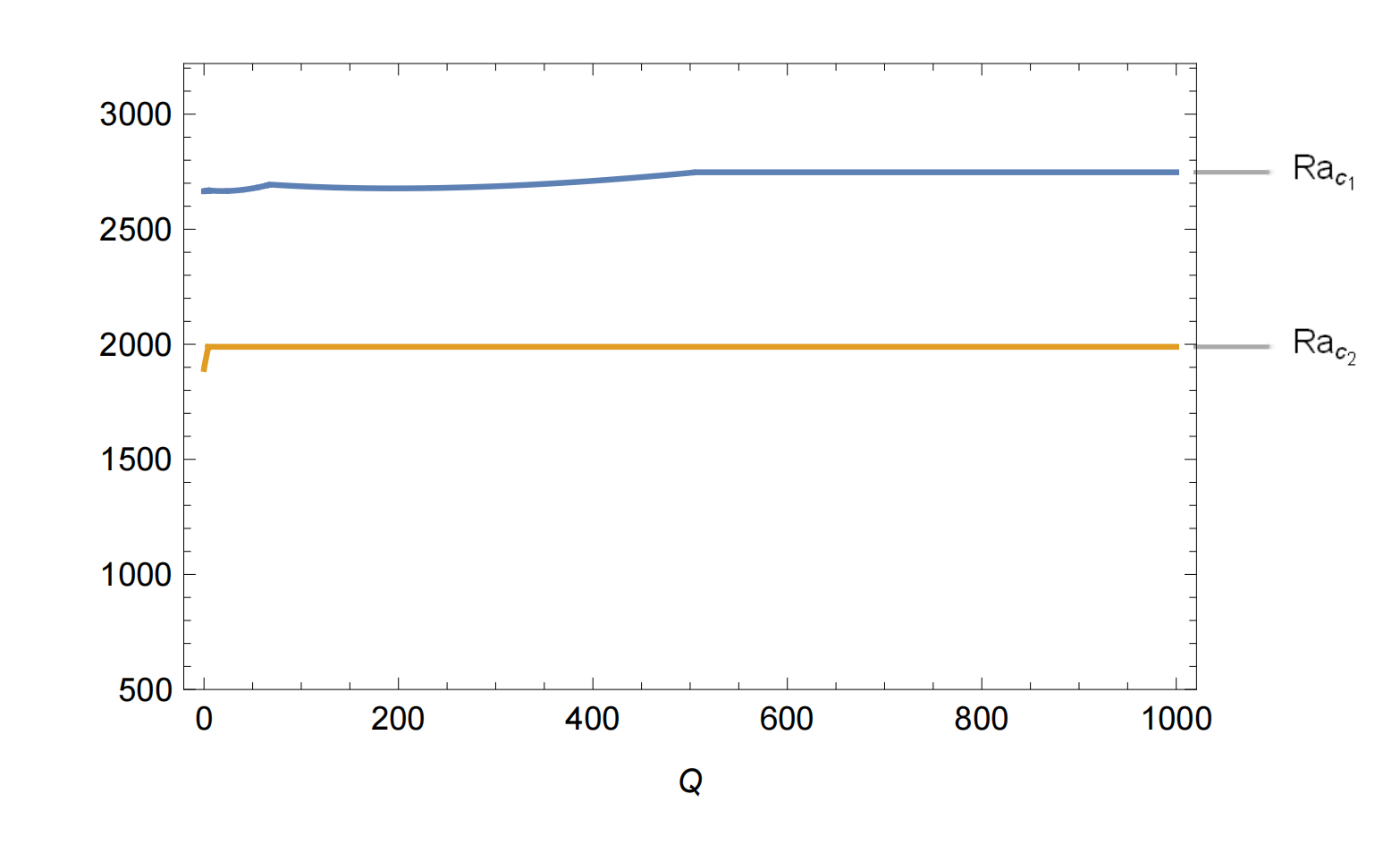}}
    \end{minipage}
    \caption{Values of \( \Ra_{c_1} \) and \( \Ra_{c_2} \) with \( \Q \in [0,1000]\) (left: \( \Ta = 2700 \), \( \Pr = 0 . 5 \), \( L_{1} = 1 \) and \( L_2 = 1. 2\) ) and (right: \( \Ta = 2700 \), \( \Pr = 0 . 2 \), \( L_{1} = 1 \) and \( L_2 = 1. 2\)).}
        \label{Ta2700Pr05L11L212}
\end{figure}

We compute the values of $\Ra_{c_{1}}$ and $\Ra_{c_{2}}$ as a function of $Q$, 
see \autoref{Ta2700Pr05L11L212}.
One can see that both the PES condition \eqref{PES1} and \eqref{PES2} can be realized for the system \eqref{abstract1}.

\begin{figure}[h]
    \begin{minipage}[t]{0.5\linewidth}
        \centering
        \begin{overpic}[width=2.3in]{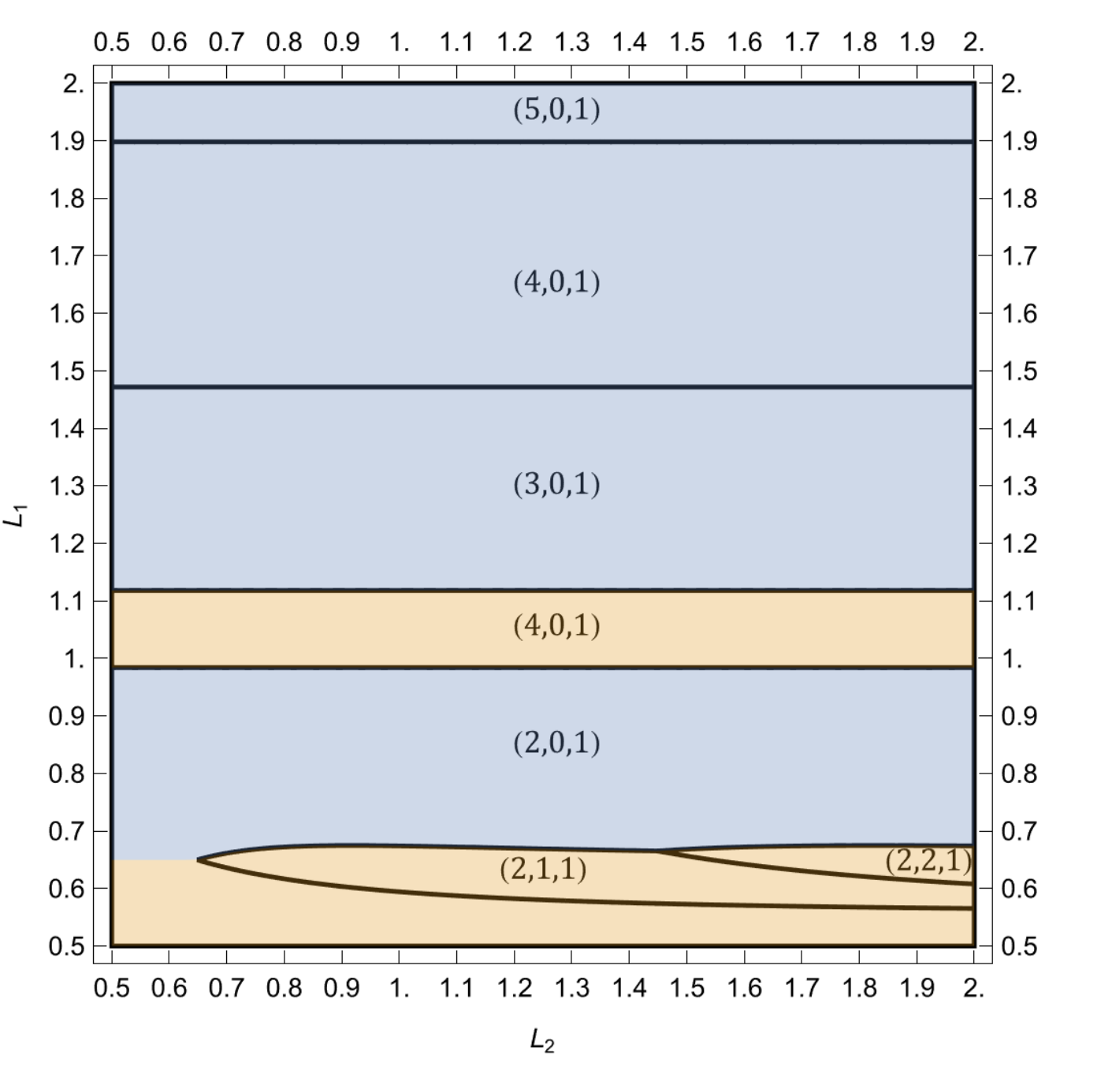}
        \end{overpic}
    \end{minipage}
    \hfill
    \begin{minipage}[t]{0.5\linewidth}
        \centering
        \begin{overpic}[width=2.3in]{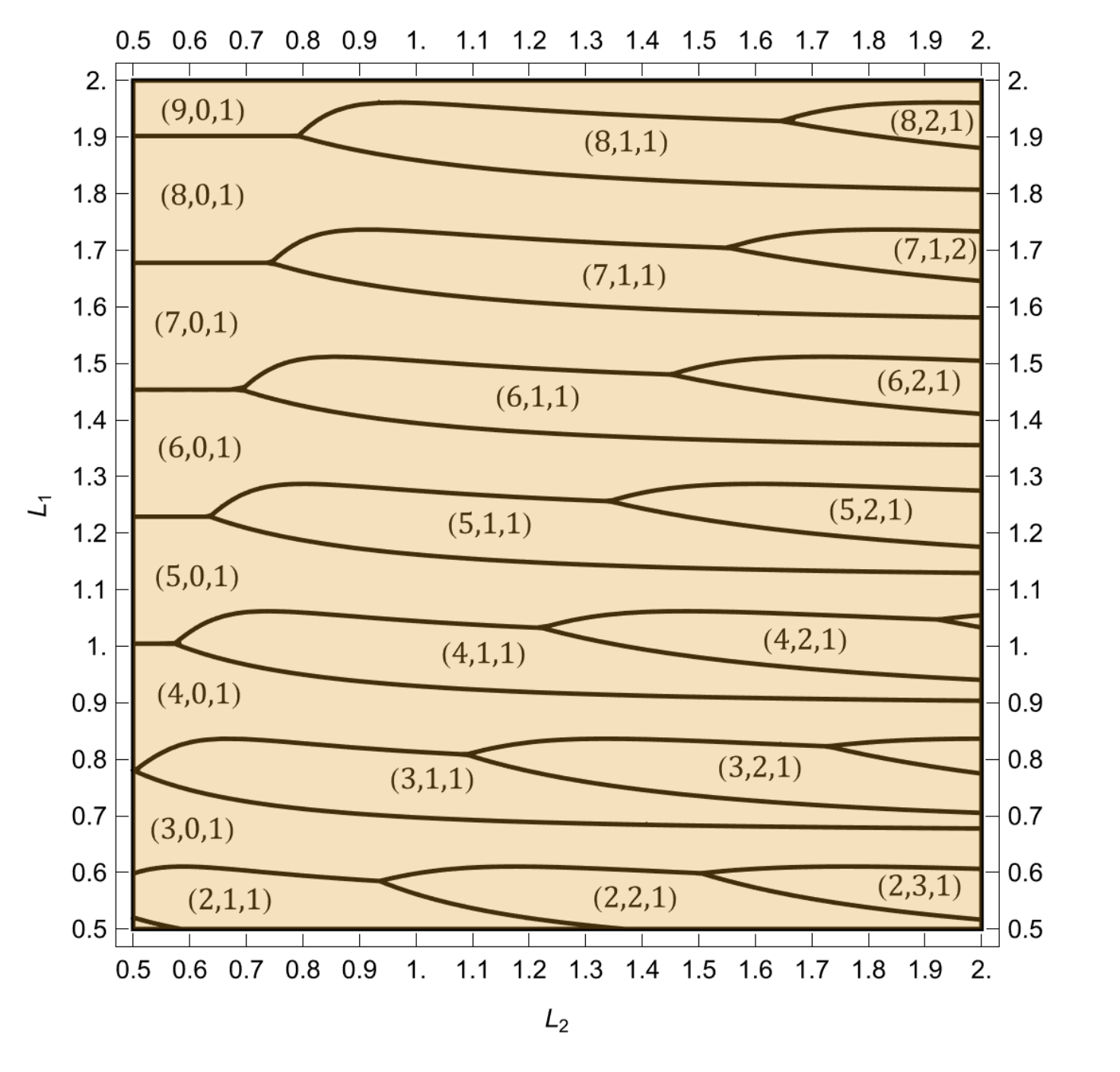}
        \end{overpic}
    \end{minipage}
    \caption{ Critical parameter region and critical index.
     (left: \( \Ta = 1100\), \( \Q = 100 \), \(\Pr = 0.2 \)) and (right: \( \Ta = 2700\), \( \Q = 100 \), \(\Pr = 0.6 \)).
    }
        \label{Ta2700Q100Pr06L1052L2052}
\end{figure}

To give more details on the PES condition with specified parameters, we numerically determine the elements of $X$, shown in  \autoref{Ta2700Q100Pr06L1052L2052}. Light blue domain is the parameter region-$(L_1,L_2 )$ in which 
    \(\Ra_{c_2} < \Ra_{c_1}\), light yellow expresses \(\Ra_{c_1} < \Ra_{c_2}\), where other parameters are specified. Values of $(j,k,l)$ shown on the different parts separated by solid black lines, are the corresponding critical index, and these solid black lines represent the parameter lines where
 there are two critical indexes.  Hence, $\text{card}(X)$ can be $1$ and $2$, even $3$. However, we only focus on the following three generic cases: 
\begin{enumerate}[(1)]
    \item $\text{card}(X)=1$ and $\beta_{J}^{1}=0$ i.e. $a_{0}=0$, where $\beta_{J}^{1}$ is the first eigenvalue which will change its sign when $\Ra$ crosses the critical value $\Ra_{c_1}$.
    \item $\text{card}(X)=1$ and $ \Re\beta_{J}^1= \Re\beta_{J}^{2}=0$ and $ \Im\beta_{J}^{1}=- \Im\beta_{J}^{2}\neq 0,$ i.e. $a_{1}a_{2}=a_{0}a_{3}$, where $\beta_{J}^{1}$ and $\beta_{J}^{2}$ are the first conjugate complex eigenvalues whose  real parts change their signs as $\Ra$ increases and crosses the critical value $\Ra_{c_2}$.
     \item $\text{card}(X)=2$ and there exists a real simple eigenvalue with multiplicity two, whose sign changes as $\Ra$ increases and crosses the critical value $Ra_{c_{1}}.$
\end{enumerate}

If first eigenvalue is real and card$(X)=1$, we always use $J_{0}=(j_{0},k_{0},1)$ to represent the critical index. And, $J_{2}=(j_{2},k_{1},1)$ and $J_{3}=(j_{3},k_{3},1)$ are utilized to express the critical indexes, if first eigenvalue is real and card$(X)=2$. 
If first eigenvalue is complex, we assume that there exists only one $J_{1}=(j_{1},k_{1},l_{1})\in I_{1}$ such that $
    \Ra_{c_{2}}=g(J_{1})$.

	\section{Transition from first real simple eigenvalue}\label{section3}
\begin{figure}[H]
		\begin{minipage}[t]{0.45\linewidth}
			\centering
			{\includegraphics[width=1.4in]{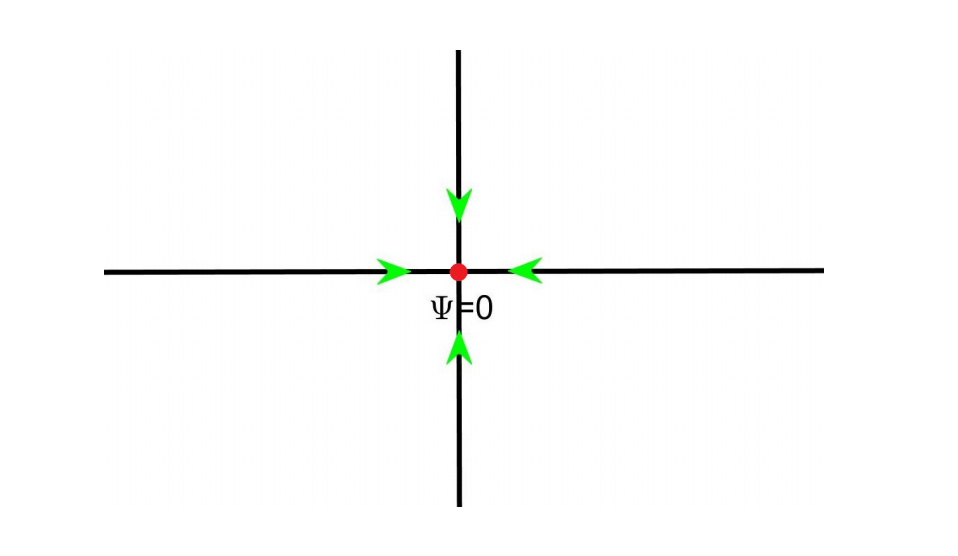}}
		\end{minipage}
		\hfill
		\begin{minipage}[t]{0.49\linewidth}
			\centering
			{\includegraphics[width=2.5in]{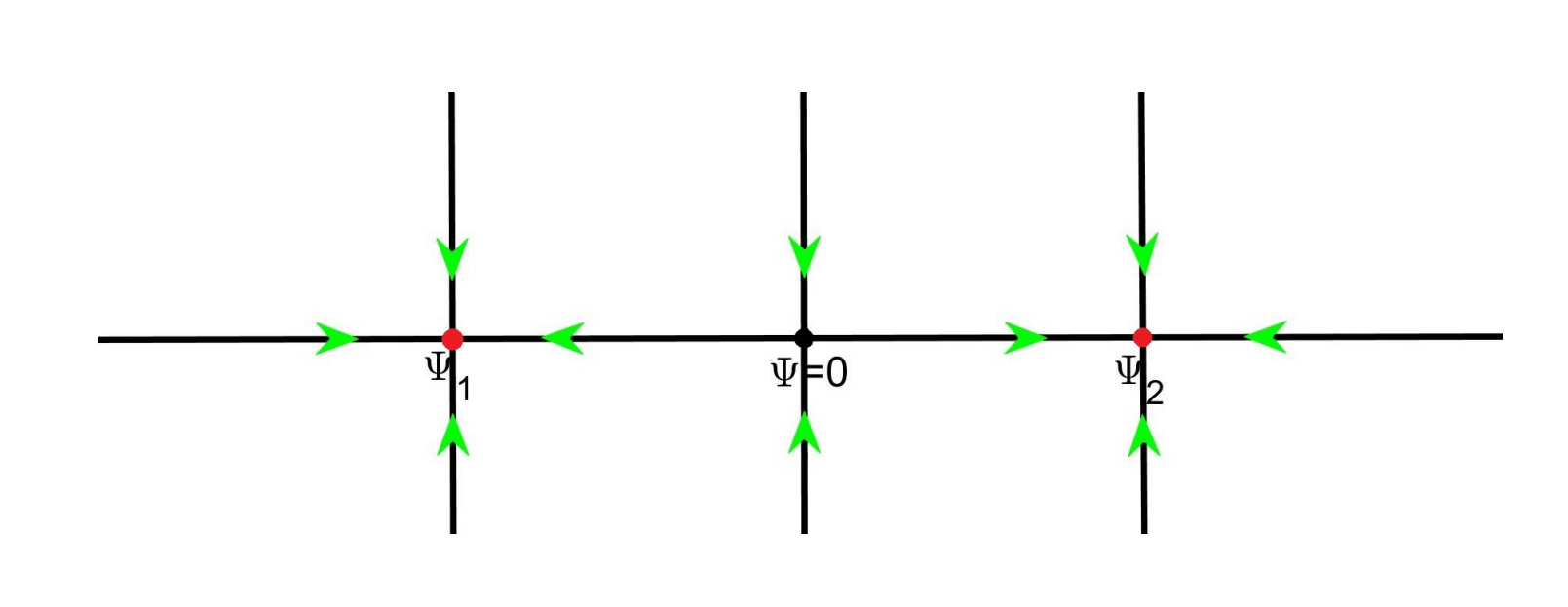}}
		\end{minipage}
		\caption{Topological structure of continuous transition of the system \eqref{abstract1}  when $\delta(\Ra_{c_{1}})<0$: (left)$\Ra< \Ra_{c_{1}}$;
		(right) $\Ra>\Ra_{c_{1}} $.}
			\label{diyigedinglitu1}
	\end{figure}

\begin{figure}[H]
		\begin{minipage}[t]{0.45\linewidth}
			\centering
			{\includegraphics[width=1.5in]{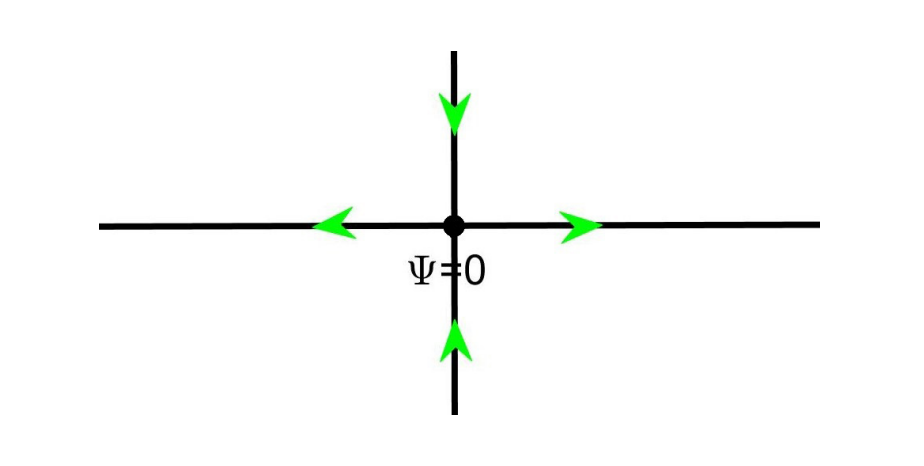}}
		\end{minipage}
		\hfill
		\begin{minipage}[t]{0.49\linewidth}
			\centering
			{\includegraphics[width=2.5in]{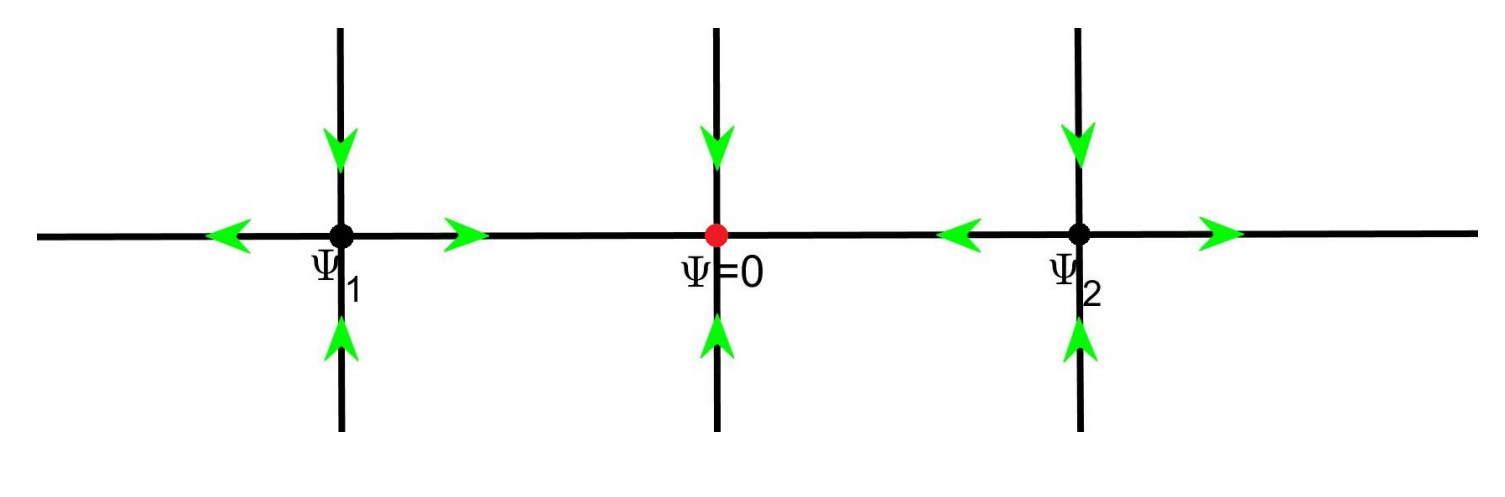}}
		\end{minipage}
		\caption{Topological structure of jump transition of the system \eqref{abstract1}  when $\delta(\Ra_{c_{1}})>0$: (left)$\Ra>\Ra_{c_{1}}$;
		(right) $\Ra<\Ra_{c_{1}} $.}
			\label{diyigedinglitu3}
	\end{figure}

	In this section, we consider the transition in the system \eqref{abstract1} from the simple real eigenvalue $\beta_{J_{0}}^{1}$. This transition is relatively simple, which can be
	completely characterized by introducing a nondimensional number $\delta(\Ra)$ given by
	\begin{align}\label{nondimensionalnumber1}
		\delta(\Ra)=\frac{\pi(2j_{0}\alpha_{1}\Phi_{2}u_{J_{0}}^{21*}+\Phi_{1}\theta_{J_{0}}^{1*}-2k_{0}\alpha_{2}\Phi_{2}u_{J_{0}}^{11*})}{u_{J_{0}}^{11}u_{J_{0}}^{11*}+u_{J_{0}}^{21}u_{J_{0}}^{21*}+1+\theta_{J_{0}}^{1}\theta_{J_{0}}^{1*}},
	\end{align}
	where $\Phi_{1}=-\frac{\Pr \theta_{J_{0}}^{1}}{8\pi}$ and $\Phi_{2}=-\frac{\pi(j_{0}\alpha_{1}u_{J_{0}}^{21}-k_{0}\alpha_{2}u_{J_{0}}^{11})}{4\Q k_{0}^{2}\alpha_{2}^{2}+16\alpha_{j_{0}k_{0}}^4}$.
	More precisely, we have:
	\begin{theorem}\label{theroemoftransition} For the system \eqref{abstract1}, we have the following conclusions:
		\begin{enumerate}
			\item[\rm{(1)}] If $\delta(\Ra_{c_{1}})<0$, it has a continuous transition from
			$(\Psi,\Ra)=(0,\Ra_{c_1})$, and bifurcates on $\Ra>\Ra_{c_{1}}$
			to a local attractor $\Sigma$ which is exactly consist of two steady-state solutions $\Psi_{1}=(\mathbf{u}_{1},\theta_{1})$ and $\Psi_{2}=(\mathbf{u}_{2},\theta_{2})$, as shown in \autoref{diyigedinglitu1}. The two steady-state solutions are approximately given by
	  \begin{align}\label{bifurcatedsolution}
					  \Psi_{m}=(-1)^{m}\sqrt{\frac{\beta_{J_{0}}^{1}(\Ra)}{-\delta(\Ra)}}\Psi_{J_{0}}^{1}+o(|\beta_{J_{0}}^{1}(\Ra)|^{\frac{1}{2}})~(m=1,2).
	 \end{align} 
	 \item[\rm{(2)}] If $\delta(\Ra_{c_{1}})>0$, it has a jump transition from $(\Psi,\Ra)=(0,\Ra_{c_1})$,
	 and bifurcates on $\Ra< \Ra_{c_{1}}$ to exactly two points 
	  $\Psi_{1}=(\mathbf{u}_{1},\theta_{1})$ and $\Psi_{2}=(\mathbf{u}_{2},\theta_{2})$,
	 as shown in \autoref{diyigedinglitu3}.
		\end{enumerate}
	\end{theorem}
	
	\begin{proof}
		According to the decomposition \eqref{decomposition}, the space $\mathbf{H_{1}}$ can be written into the direct sum of two subspaces as follows:
		\begin{align}\label{decompositionofspace}
			\mathbf{H_{1}}=\overline{E_{0}} +\overline{E_{1}},\quad \overline{E_{0}}=\text{span}\{\Psi_{J_{0}}^{1}\},
		\end{align}
	 where $\Psi_{J_{0}}^{1}$ is the eigenvector corresponding to the eigenvalue $\beta_{J_{0}}^{1}$, and $\overline{E_{1}}$ is the space spanned by the rest of eigenvectors.  Note that the PES condition \eqref{PES1} means that there exists a center manifold function $\Phi$ $:\overline{E_{0}}\rightarrow \overline{E_{1}}$ at $\Ra=\Ra_{c_1}$. And, for the solution $\Psi\in \mathbf{H_{1}}$ of \eqref{abstract1},  \eqref{decompositionofspace} infers that it has the following spectral decomposition
		\begin{align}\label{decompositionoffunction}
			\Psi=x\Psi_{J_{0}}^{1}+\Phi(x)=x\Psi_{J_{0}}^{1}+\sum_{(J,s)\neq(J_{0},1)}x_{J}^{s}(x)\Psi_{J}^{s},\quad \Phi(x)\in\overline{E_{1}}.
		\end{align} 
	Making use of \eqref{decompositionoffunction}, the system \eqref{abstract1} is then reduced into 
	 \begin{align}\label{reduceode}
			\frac{dx}{dt}=\beta_{J_{0}}^{1}x+\frac{\left\langle \mathbf{G}(\Psi), \Psi_{J_{0}}^{1*}\right\rangle }{\left\langle \Psi_{J_{0}}^{1}, \Psi_{J_{0}}^{1*}\right\rangle},
		\end{align}
	which is derived by taking inner products of $\Psi_{J_{0}}^{1*}$ with  both sides of \eqref{abstract1}.
	
	Because $G:\mathbf{H_{1}}\times\mathbf{H_{1}}\rightarrow \mathbf{H}$ is a bilinear form, one gets
		\begin{align}\label{bilinearform}
		   \mathbf{G}(\Psi)=G(\Psi,\Psi)=x^2\mathbf{G}(\Psi_{J_{0}}^1)+x [G(\Psi_{J_{0}}^{1},\Phi)+G(\Phi,\Psi_{J_{0}}^{1})]+\mathbf{G}(\Phi).
		\end{align}
	A direct calculation further gives  
		\begin{align}\label{1G11}
			\mathbf{G}(\Psi_{J_{0}}^{1})=\frac{\pi}{2}\begin{pmatrix}
				\mathbf{P}\begin{pmatrix}
					u_{J_{0}}^{11}\sin{2(j_{0}\alpha_{1}x_{1}+k_{0}\alpha_{2}x_{2})}
					\\
					u_{J_{0}}^{21}\sin{2(j_{0}\alpha_{1}x_{1}+k_{0}\alpha_{2}x_{2})}
					\\
					0
				\end{pmatrix}
				\\
				0
			\end{pmatrix}-\frac{\pi}{2}\begin{pmatrix}
				0
				\\
				0
				\\
				0
				\\
				\theta_{J_{0}}^{1}\sin{2\pi x_{3}}
			\end{pmatrix},
		\end{align}
		thus, $\left\langle \mathbf{G}(\Psi_{J_{0}}^{1}), \Psi_{J_{0}}^{1*}\right\rangle=0. $
		Furthermore, the formula (A.1.14) ( \cite{Ma1} ) says that $\Phi$ solves
		\begin{align}\label{centermanifold}
			L\Phi=\sum_{(J,s)\neq(J_{0},1)}x_{J}^{s}\beta_{J}^{s}\Psi_{J}^{s}=-x^2 \mathbf{G}(\Psi_{J_{0}}^{1})+o(2),
		\end{align}
	from which we have
	\[
	x_{J}^{s}=-\frac{x^{2}\left\langle \mathbf{G}(\Psi_{J_{0}}^{1}),\Psi_{J}^{s*}\right\rangle }{\beta_{J}^{s}\left\langle \Psi_{J}^{s},\Psi_{J}^{s*} \right\rangle }.
	\]
		
	Making use of \eqref{1G11}, one can see that for $(J,s)\neq((0,0,2),3)$ and $J\neq (2j_{0},2k_{0},0)$, we have
		$\left\langle \mathbf{G}(\Psi_{J_{0}}^{1}),\Psi_{J}^{s*}\right\rangle=0,$
		and for $(J,s)=((0,0,2),3)$ and $J=2(j_{0},k_{0},0)$, it yields 
		\begin{align}\label{xishu}
			\begin{cases}
				x_{(0,0,2)}^{3}=-\frac{\Pr \theta_{J_{0}}^{1}x^{2}}{8\pi},
				\\
				x_{(2j_{0},2k_{0},0)}=-\frac{\pi(j_{0}\alpha_{1}u_{J_{0}}^{21}-k_{0}\alpha_{2}u_{J_{0}}^{11})x^{2}}{4\Q k_{0}^{2}\alpha_{2}^{2}+16\alpha_{j_{0}k_{0}}^4}.
			\end{cases}
		\end{align}
	Hence,  $\Phi$ has the following expansion 
		\begin{align}
			\Phi=x^{2}(\Phi_{1}\Psi_{(0,0,2)}^{3}+\Phi_{2}\Psi_{(2j_{0},2k_{0},0)})+o({2}),
		\end{align}
		where 
		\[
		\Phi_{1}=-\frac{\Pr \theta_{J_{0}}^{1}}{8\pi},\quad \Phi_{2}=-\frac{\pi(j_{0}\alpha_{1}u_{J_{0}}^{21}-k_{0}\alpha_{2}u_{J_{0}}^{11})}{4\Q k_{0}^{2}\alpha_{2}^{2}+16\alpha_{j_{0}k_{0}}^4}.
		\]
		
		And by computation, we have $
			G(\Phi,\Psi_{J_{0}}^{1})=0$,
		and 
		\begin{align}\label{phi2}
			G(\Psi_{J_{0}}^{1},\Phi)=\pi x^2
			\begin{pmatrix}
				2k_{0}\alpha_{2}\Phi_{2}[\sin{3f_{0102}}-\sin{f_{0102}}]\cos{\pi x_{3}}
				\\
				-2j_{0}\alpha_{1}\Phi_{2}[\sin{3f_{0102}}-\sin{f_{0102}}]\cos{\pi x_{3}}
				\\
				0
				\\
				-\Phi_{1}\cos{f_{0102}(\sin{3\pi x_{3}}-\sin{\pi x_{3}}})
			\end{pmatrix}+o({2}).
		\end{align}
		where $f_{0102}=j_{0}\alpha_{1}x_{1}+k_{0}\alpha_{2}x_{2}$.
		Thus, the reduced system \eqref{reduceode} can be simplified as
		\begin{align}\label{transitionequation}
			\frac{dx}{dt}=\beta_{J_{0}}^{1}x+\delta(\Ra)x^{3}+ o(3),
		\end{align}
		where 
		\begin{align}\label{transitionequationxishu}
			\delta(\Ra)=\frac{\left\langle G(\Psi_{J_{0}}^{1},\Phi),\Psi_{J_{0}}^{1*} \right\rangle }{\left\langle \Psi_{J_{0}}^{1}, \Psi_{J_{0}}^{1*}\right\rangle}=\frac{\pi(2j_{0}\alpha_{1}\Phi_{2}u_{J_{0}}^{21*}+\Phi_{1}\theta_{J_{0}}^{1*}-2k_{0}\alpha_{2}\Phi_{2}u_{J_{0}}^{11*})}{u_{J_{0}}^{11}u_{J_{0}}^{11*}+u_{J_{0}}^{21}u_{J_{0}}^{21*}+1+\theta_{J_{0}}^{1}\theta_{J_{0}}^{1*}}.
		\end{align}
		
	 Finally, by analyzing dynamics of \eqref{transitionequation}
	 at $\Ra=\Ra_c$. one can obtain the theorem.    
		
	\end{proof}

	\section{Transition from first complex eigenvalues}\label{section4}
	In this section, we consider the transition of the system \eqref{abstract1} at $\Ra=\Ra_{c_{2}}$, where a pair of complex conjugate eigenvalues become critical.  The critical index is denoted by $J_{1}=(j_{1}, k_{1}, l_{1})\in X$ which also satisfies $j_{1}k_{1}=0$ and $(j_{1},k_{1})\neq(0,0)$, and the critical conjugate eigenvectors, corresponding to the critical eigenvalues $\beta_{J_{1}}^{1}=\overline{\beta_{J_{1}}^{2}}=\sigma +i\rho$, are $\Psi_{J_{1}}^{1}$ and $\Psi_{J_{1}}^{2}$, respectively. Likewise,  this transition can also be
	completely characterized by the following nondimensional number $a(\Ra)$ 
	\begin{align}\label{nondimensionalnumber2}
		a(\Ra)=3a_{111}(\Ra)+a_{122}(\Ra)+b_{112}(\Ra),
	\end{align}
	 where the exact expressions of $a_{111}$, $a_{122}$ and $b_{122}$ can be found in \eqref{nondimensionalcoefficients}.
	
	\begin{theorem} \label{complex-transition}
		For the system \eqref{abstract1}, we have the following conclusions:
		\begin{enumerate}
			\item[\rm{(1)}] If $a(\Ra_{c_{2}})<0$, it has a continuous transition
	from $(\Psi,\Ra)=(0,\Ra_{c_2})$, and bifurcates on $Ra>Ra_{c_{2}}$ to a stable periodic orbit $\Psi_p$ which is approximately given by 
			\begin{align}\label{zhouqiguidao}
				\begin{aligned}
					&\Psi_p=\left(\frac{4\sigma}{-\pi a(Ra_{c_{2}})}\right)^{\frac{1}{2}}\sin(\rho t)\Psi_{J_{1}}^{1}+\left(\frac{4\sigma}{-\pi a(Ra_{c_{2}})}\right)^{\frac{1}{2}}\cos(\rho t)\Psi_{J_{1}}^{2}+o(\sqrt{|\sigma|}).
				\end{aligned}
			\end{align}
				\item[\rm{(2)}] If $a(\Ra_{c_{2}})>0$, it has a jump transition 
				from $(\Psi,\Ra)=(0,\Ra_{c_2})$, and bifurcates on $\Ra<\Ra_{c_{2}}$ to
				 an unstable periodic orbit $\Psi_p$ with same expression as in
				 \eqref{zhouqiguidao}.
		\end{enumerate}
	\end{theorem}

	\begin{proof} Similarly, the space $\mathbf{H_{1}}$ can be written into the direct sum of two subspaces:
		\begin{align}\label{decompositionofspaces2}
			\mathbf{H_{1}}= \widetilde{E_{1}}+\widetilde{E_{2}}, \quad \widetilde{E_{1}}=\text{span}\{\Psi_{J_{1}}^{1},\Psi_{J_{1}}^{2}\},
		\end{align}
		where $\widetilde{E_{2}}$ is spanned by the rest of eigenvectors. 
		Hence, on the center manifold, the solution $\Psi\in \mathbf{H_{1}}$ of \eqref{abstract1} has the following spectral decomposition
		\begin{align}\label{decompositionoffunction2}
			\begin{aligned}
			\Psi&=x\Psi_{J_{1}}^{1}+y\Psi_{J_{1}}^{2}+\Phi(x,y)\\&=x\Psi_{J_{1}}^{1}+y\Psi_{J_{1}}^{2}+\sum\limits_{(J,s)\neq(J_{1},m)}x_{J}^{s}(x,y)\Psi_{J}^{s},\quad \Phi(x,y)\in\overline{E_{1}},\quad m=1,2.
				\end{aligned} 
		\end{align} 
		
		In order to reduce Eq. \eqref{abstract1} into a system of ODEs, we introduce the following two 
		conjugate eigenvectors $\widetilde{\Psi_{J_{1}}^{1*}}$ and $\widetilde{\Psi_{J_{1}}^{2*}}$:
		\begin{align}\label{conjugateeigenvectors2}
				\widetilde{\Psi_{J_{1}}^{1*}}=\Psi_{J_{1}}^{1*}+\alpha\Psi_{J_{1}}^{2*},\quad
				\widetilde{\Psi_{J_{1}}^{2*}}=-\alpha\Psi_{J_{1}}^{1*}+\Psi_{J_{1}}^{2*}
		\end{align}
		where $\alpha=\frac{\left\langle \Psi_{J_{1}}^{1},\Psi_{J_{1}}^{2*} \right\rangle }{\left\langle \Psi_{J_{1}}^{1},\Psi_{J_{1}}^{1*} \right\rangle }$$=\frac{\Re  u_{J_{1}}^{11} \Im  u_{J_{1}}^{11*}+\Re  u_{J_{1}}^{21}\Im  u_{J_{1}}^{21*}+\Re  \theta_{J_{1}}^{1} \Im  \theta_{J_{1}}^{1*}}{\Re  u_{J_{1}}^{11} \Re  u_{J_{1}}^{11*}+\Re  u_{J_{1}}^{21} \Re  u_{J_{1}}^{21*}+1+\Re  \theta_{J_{1}}^{1} \Re  \theta_{J_{1}}^{1*}}$, $L\Psi_{J_{1}}^{1}=\sigma\Psi_{J_{1}}^{1}-\rho\Psi_{J_{1}}^{2}$, $L\Psi_{J_{1}}^{2}=\rho\Psi_{J_{1}}^{1}+\sigma\Psi_{J_{1}}^{2}$, $L^{*}\Psi_{J_{1}}^{1*}=\sigma\Psi_{J_{1}}^{1*}+\rho\Psi_{J_{1}}^{2*}$ and $L^{*}\Psi_{J_{1}}^{2*}=-\rho\Psi_{J_{1}}^{1*}+\sigma\Psi_{J_{1}}^{2*}$.
	Making use of the expansion \eqref{decompositionoffunction2} and conducting
	some computation, we have
			\begin{align}\label{adjustment}
				\begin{aligned}
				\left\langle \Psi, \widetilde{\Psi_{J_{1}}^{1*}}\right\rangle=x\left\langle \Psi_{J_{1}}^{1},\widetilde{\Psi_{J_{1}}^{1*}}\right\rangle \neq 0,\quad
				\left\langle \Psi, \widetilde{\Psi_{J_{1}}^{2*}}\right\rangle=y\left\langle \Psi_{J_{1}}^{2},\widetilde{\Psi_{J_{1}}^{2*}}\right\rangle \neq 0.
				\end{aligned} 
			\end{align}
			
	For $m=1,2$, taking inner product of $\widetilde{\Psi_{J}^{m}}$ with both sides of Eq.\eqref{abstract1}, we have the following ODEs:
			\begin{align}\label{ODEs2}
				\begin{cases}
					\frac{dx}{dt}=\sigma x+\rho y+\frac{\left\langle \mathbf{G}(\Psi),\widetilde{\Psi_{J_{1}}^{1*}}\right\rangle }{\left\langle \Psi_{J_{1}}^{1},\widetilde{\Psi_{J_{1}}^{1*}}\right\rangle},
					\\
					\frac{dy}{dt}=-\rho x+\sigma y+\frac{\left\langle \mathbf{G}(\Psi),\widetilde{\Psi_{J_{1}}^{2*}}\right\rangle }{\left\langle \Psi_{J_{1}}^{2},\widetilde{\Psi_{J_{1}}^{2*}}\right\rangle},
				\end{cases}
			\end{align}
			where 
			\begin{align}
			   & \begin{aligned}\label{innerproduct}
					\left\langle \Psi_{J_{1}}^{1},\widetilde{\Psi_{J_{1}}^{1*}}\right\rangle =
					 & \pi^2(\alpha_{1}\alpha_{2})^{-2}(\Re u_{J_{1}}^{11}\Re u_{J_{1}}^{11*}+\Re u_{J_{1}}^{21}\Re u_{J_{1}}^{21*}+1 +\Re \theta_{J_{1}}^{1}\Re \theta_{J_{1}}^{1*})
					\\
					 & +\pi^2(\alpha_{1}\alpha_{2})^{-2}\alpha(\Re u_{J_{1}}^{11}\Im u_{J_{1}}^{11*}+\Re u_{J_{1}}^{21}\Im u_{J_{1}}^{21*}+\Re \theta_{J_{1}}^{1}\Im \theta_{J_{1}}^{1*}),
				\end{aligned} \\
			   & \begin{aligned}\label{innerproduct1}
					\left\langle \Psi_{J_{1}}^{2},\widetilde{\Psi_{J_{1}}^{2*}}\right\rangle =
					 & \pi^2(\alpha_{1}\alpha_{2})^{-2}(\Im u_{J_{1}}^{11}\Im u_{J_{1}}^{11*}+\Im u_{J_{1}}^{21}\Im u_{J_{1}}^{21*}+\Im \theta_{J_{1}}^{1}\Im \theta_{J_{1}}^{1*})
					\\
					 & -\pi^2(\alpha_{1}\alpha_{2})^{-2}\alpha(\Im u_{J_{1}}^{11}\Re u_{J_{1}}^{11*}+\Im u_{J_{1}}^{21}\Re u_{J_{1}}^{21*}+\Im \theta_{J_{1}}^{1}\Re \theta_{J_{1}}^{1*}).
				\end{aligned} 
			\end{align}
			
	With the help of the decomposition \eqref{decompositionoffunction2}, we have
	the following expansion
			\begin{align}\label{decomposeG}
				\begin{aligned}
					\mathbf{G}(\Psi)={} & x^{2}G_{11}+xy(G_{12}+G_{21})+y^{2}G_{22}
					 +x[G(\Psi_{J_{1}}^{1},\Phi)+G(\Phi,\Psi_{J_{1}}^{1})]
					\\
							 & +y[G(\Psi_{J_{1}}^{2},\Phi)+G(\Phi,\Psi_{J_{1}}^{2})]+o(2),
				\end{aligned}
			\end{align}
			where $G_{mn}=G(\Psi_{J_{1}}^{m},\Psi_{J_{1}}^{n})~(m,n=1,2).$
	  By simple computation, one gets
			\begin{align}\label{G11}
				&G_{11}=\frac{l_{1}\pi}{2}\begin{pmatrix}
					\mathbf{P}\begin{pmatrix}
						\Re  u_{J_{1}}^{11}\sin{2(j_{1}\alpha_{1}x_{1}+k_{1}\alpha_{2}x_{2})}
						\\
						\Re  u_{J_{1}}^{21}\sin{2(j_{1}\alpha_{1}x_{1}+k_{1}\alpha_{2}x_{2})}
						\\
						0
					\end{pmatrix}
					\\
					0
				\end{pmatrix}
				-\frac{l_{1}\pi\Re  \theta_{J_{1}}^{1}}{2}\begin{pmatrix}
					0
					\\
					0
					\\
					0
					\\
					\sin{2l_{1}\pi x_{3}}
				\end{pmatrix}, \\ \label{G12}
			   & G_{12}=\frac{l_{1}\pi}{2}\begin{pmatrix}
					\Im  u_{J_{1}}^{11}\sin{2(j_{1}\alpha_{1}x_{1}+k_{1}\alpha_{2}x_{2})}
					\\
					\Im  u_{J_{1}}^{21}\sin{2(j_{1}\alpha_{1}x_{1}+k_{1}\alpha_{2}x_{2})}
					\\
					0
					\\
					0
				\end{pmatrix}
				-\frac{l_{1}\pi \Im  \theta_{J_{1}}^{1}}{2}\begin{pmatrix}
					0
					\\
					0
					\\
					0
					\\
					\sin{2l_{1}\pi x_{3}}
				\end{pmatrix},
			\end{align}
			and 
			\begin{align}\label{G2m}
				G_{21}=G_{22}=\mathbf{0}.
			\end{align}
			
			 Utilizing the formula (A.1.19) (Appendix A.1 in \cite{Ma1}), we have 
			\begin{align}\label{centermanifold2}
				\Phi=\Phi_{1}+\Phi_{2}+\Phi_{3}+o(2),
			\end{align}
			where $\Phi_{1}-\Phi_{3}$ can be solved from
			\begin{align}\label{solvecentermanifold}
				\begin{aligned}
				   & L\Phi_{1}={}               - \left(x^2 G_{11}+xyG_{12}\right),
					\\
				   & (L^2+4\rho^2)L\Phi_{2}={}  2\rho^2  \left[\left(x^2-y^2\right)G_{11}+2xyG_{12}\right],
					\\
				   & (L^2+4\rho^2)\Phi_{3}={}   -\rho  \left[\left(y^2-x^{2}\right) G_{12}+2xyG_{11}\right]. 
				\end{aligned}
			\end{align}
	 And,   conducting computation obtains
			\begin{align}\label{phi1phi2phi3}
				\begin{aligned}
					\Phi_{1} & =\Phi_{11}\Psi_{(2j_{1},2k_{1},0)}+\Phi_{12}\Psi_{(0,0,2l_{1})}^{3},
					\\
					\Phi_{2} & =\Phi_{21}\Psi_{(2j_{1},2k_{1},0)}+\Phi_{22}\Psi_{(0,0,2l_{1})}^{3},
					\\
					\Phi_{3} & =\Phi_{31}\Psi_{(2j_{1},2k_{1},0)}+\Phi_{32}\Psi_{(0,0,2l_{1})}^{3},
				\end{aligned}
			\end{align}
			where 
			\begin{align*}
			   & \begin{aligned}
					\Phi_{11} & =\frac{k_{1}\alpha_{2}M_{1}-j_{1}\alpha_{1}M_{2}}{2(\Q k_{1}^{2}\alpha_{2}^{2}+4\alpha_{j_{1}k_{1}}^{4})},                                                                                                                                   & \Phi_{12} & =-\frac{M_{3}\Pr }{4l_{1}^{2}\pi^{2}},
					\\
					\Phi_{21} & =\frac{\frac{(j_{1}\alpha_{1}M_{5}-k_{1}\alpha_{2}M_{4})\alpha_{j_{1}k_{1}}^{2}}{2(\Q k_{1}^{2}\alpha_{2}^{2}+4\alpha_{j_{1}k_{1}}^{4})}}{(\Q k_{1}^{2}\alpha_{2}^{2}+4\alpha_{j_{1}k_{1}}^{4})^{2}+4\alpha_{j_{1}k_{1}}^{2}\rho^{2}}, & \Phi_{22} & =\frac{M_{6}\Pr ^{3}}{16l_{1}^{2}\pi^{2}(4l_{1}^{4}\pi^{4}+\rho^{2}\Pr ^{2})},
					\\
					\Phi_{31} & =\frac{(k_{1}\alpha_{2}M_{7}-j_{1}\alpha_{1}M_{8})\alpha_{j_{1}k_{1}}^{2}}{4[(\Q k_{1}^{2}\alpha_{2}^{2}+4\alpha_{j_{1}k_{1}}^{4})^{2}+4\rho^2 \alpha_{j_{1}k_{1}}^{2}]}           ,                                                         & \Phi_{32} & =\frac{M_{9}\Pr ^{2}}{16l_{1}^2\pi^{2}+4\rho^{2}\Pr ^{2}},
				\end{aligned}
			\end{align*}
			\begin{align*}                   
			   & \begin{aligned}
					M_{1} & =\frac{l_{1}\pi x}{2}(y \Re  u_{J_{1}}^{11}+y\Im  u_{J_{1}}^{11}),~~M_{4}  =l_{1}\pi \rho^{2} \left[  \left( x^{2}-y^{2} \right)\Re  u_{J_{1}}^{11}+2 x y \Im  u_{J_{1}}^{11} \right], \\
					M_{2} & =\frac{l_{1}\pi x}{2}(y \Re  u_{J_{1}}^{21}+y\Im  u_{J_{1}}^{21}),~~ M_{5}  =l_{1}\pi \rho^{2} \left[  \left(x^{2}-y^{2}\right)\Re  u_{J_{1}}^{21}+2 x y \Im  u_{J_{1}}^{21}\right],\\
					M_{3} & =\frac{l_{1}\pi x}{2}(x\Re  \theta_{J_{1}}^{1}+y \Im  \theta_{J_{1}}^{1}),~~M_{6}  =l_{1}\pi \rho^{2}  \left[  \left(x^{2}-y^{2}\right) \Re  \theta_{J_{1}}^{1}+2 x y \Im  \theta_{J_{1}}^{1} \right],
				\end{aligned}\\&
				\begin{aligned}
				   & M_{7}=l_{1}\pi\rho[(x^{2}-y^{2}) \Im  u_{J_{1}}^{11}+2xy\Re  u_{J_{1}}^{11}], ~~
					M_{8}=l_{1}\pi\rho[(x^{2}-y^2)\Im  u_{J_{1}}^{21}+2xy\Re  u_{J_{1}}^{21}], 
					\\  
					&M_{9}=\frac{l_{1}\pi\rho}{2}[(y^2-x^{2})\Im  \theta_{J_{1}}^{1}+2xy\Re  \theta_{J_{1}}^{1}]. 
				\end{aligned}
			\end{align*}
			
	 Upon performing some computation, one can obtain
			\begin{align}\label{reducedequationqianmian}
				\begin{aligned}
					\left\langle G_{11},\widetilde{\Psi_{J_{1}}^{n*}}\right\rangle= \left\langle G_{12},\widetilde{\Psi_{J_{1}}^{n*}}\right\rangle=0\quad  (n=1,2), \\
					G(\Psi_{J_{1}}^{2},\Phi)=G(\Phi,\Psi_{J_{1}}^{1})=G(\Phi,\Psi_{J_{1}}^{2})=0.
				\end{aligned}
			\end{align}
	 and
				\begin{align}\label{202111031}
					 \begin{aligned}
				   & \frac{\left\langle G(\Psi_{J_{1}}^{1},\Phi),\widetilde{\Psi_{J_{1}}^{1*}}\right\rangle}{\left\langle \Psi_{J_{1}}^{1},\widetilde{\Psi_{J_{1}}^{1*}}\right\rangle}
					=\sum_{s=1}^2\sum_{j=1}^3A_{1s}\Phi_{js},\\
				   & \frac{\left\langle G(\Psi_{J_{1}}^{1},\Phi),\widetilde{\Psi_{J_{1}}^{2*}}\right\rangle}{\left\langle \Psi_{J_{1}}^{2},\widetilde{\Psi_{J_{1}}^{2*}}\right\rangle}=\sum_{s=1}^2\sum_{j=1}^3A_{2s}\Phi_{js},
								\end{aligned}
				\end{align}
			where
			\begin{align*}
				\begin{cases}
					A_{11} =\frac{2l_{1}\pi^3}{\alpha_{1}\alpha_{2} \left\langle \Psi_{J_{1}}^{1},\widetilde{\Psi_{J_{1}}^{1*}}\right\rangle}(j_{1}\alpha_{1}\Re  u_{J_{1}}^{21*}+j_{1}\alpha_{1}\alpha \Im  u_{J_{1}}^{21*}-k_{1}\alpha_{2}\Re  u_{J_{1}}^{11*}-k_{1}\alpha_{2}\alpha \Im  u_{J_{1}}^{11*}), 
					\\
					A_{12} =\frac{l_{1}\pi^3}{\alpha_{1}\alpha_{2}\left\langle \Psi_{J_{1}}^{1},\widetilde{\Psi_{J_{1}}^{1*}}\right\rangle}(\Re \theta_{J_{1}}^{1*}+\alpha\Im \theta_{J_{1}}^{1*}),
					\\
					A_{21}=\frac{2l_{1}\pi^3}{\alpha_{1}\alpha_{2}\left\langle \Psi_{J_{1}}^{2},\widetilde{\Psi_{J_{1}}^{2*}}\right\rangle }(j_{1}\alpha_{1}\Im  u_{J_{1}}^{21*}-j_{1}\alpha_{1}\alpha \Re  u_{J_{1}}^{21*}-k_{1}\alpha_{2}\Im  u_{J_{1}}^{11*}+k_{1}\alpha_{2}\alpha \Re  u_{J_{1}}^{11*}), 
					\\
					A_{22}=\frac{l_{1}\pi^3}{\alpha_{1}\alpha_{2}\left\langle \Psi_{J_{1}}^{2},\widetilde{\Psi_{J_{1}}^{2*}}\right\rangle}(\Im \theta_{J_{1}}^{1*}-\alpha\Re \theta_{J_{1}}^{1*}).
				\end{cases}
			\end{align*}
			
		   With the help of \eqref{innerproduct}-\eqref{G2m} and \eqref{202111031}, we simplify the system \eqref{ODEs2} into
			\begin{align}\label{reducedODE2}
				\begin{aligned}
					\frac{dx}{dt} & =\sigma x+\rho y+a_{111}x^{3}+a_{112}x^{2}y+a_{122}xy^{2},
					\\
					\frac{dy}{dt} & =-\rho x+\sigma y+b_{111}x^{3}+b_{112}x^{2}y+b_{122}xy^{2},
				\end{aligned}
			\end{align}
			where 
			\begin{align}\label{nondimensionalcoefficients}
				\begin{aligned}
					a_{111} & =C_{1}A_{11}+C_{2}A_{21}, & a_{112} & =C_{3}A_{11}+C_{4}A_{21}, & a_{122} & =C_{5}A_{11}+C_{6}A_{21},
					\\ 
					b_{111} & =C_{1}A_{12}+C_{2}A_{22}, & b_{112} & =C_{3}A_{12}+C_{4}A_{22}, & b_{122} & =C_{5}A_{12}+C_{6}A_{22},
				\end{aligned}
			\end{align}
			and in which $ C_{1}- C_{5}$ are explicitly given as follows
			\begin{align*}
				C_{1}= {} & \frac{l_{1}\pi(k_{1}\alpha_{2} \Re  u_{J_{1}}^{11}-j_{1}\alpha_{1}\Re  u_{J_{1}}^{21})}{4(\Q k_{1}^{2}\alpha_{2}^{2}+4\alpha_{j_{1}k_{1}}^{4})}+\frac{\frac{l_{1}\pi\rho^2(j_{1}\alpha_{1} \Re  u_{J_{1}}^{21}-k_{1}\alpha_{2}\Re  u_{J_{1}}^{11})\alpha_{j_{1}k_{1}}^{2}}{2(\Q k_{1}^{2}\alpha_{2}^{2}+4\alpha_{j_{1}k_{1}}^{4})}}{(\Q k_{1}^{2}\alpha_{2}^{2}+4\alpha_{j_{1}k_{1}}^{4})^{2}+4\alpha_{j_{1}k_{1}}^2 \rho^2}
				\\
					   & +\frac{l_{1}\pi \rho(k_{1}\alpha_{2}\Im  u_{J_{1}}^{11}-j_{1}\alpha_{1}\Im  u_{J_{1}}^{21})\alpha_{j_{1}k_{1}}^{2}}{4 [(\Q k_{1}^{2}\alpha_{2}^{2}+4\alpha_{j_{1}k_{1}}^{4})^{2}+4\alpha_{j_{1}k_{1}}^{2}\rho^2]},
				\\
				C_{2}= {} & -\frac{\Pr \Re  \theta_{J_{1}}^{1}}{8l_{1}\pi}+\frac{\Pr ^{3}\rho^{2} \Re \theta_{J_{1}}^{1}}{4l_{1}\pi(16l_{1}^{4}\pi^4+4\Pr ^{2}\rho^{2})}-\frac{l_{1}\pi\rho \Pr ^{2}\Im  \theta_{J_{1}}^{1}}{32l_{1}^{2}\pi^{2}+8\Pr ^{2}\rho^{2}},
				\\
				C_{3}= {} & \frac{l_{1}\pi(k_{1}\alpha_{2}\Im  u_{J_{1}}^{11}-j_{1}\alpha_{1} \Im  u_{J_{1}}^{21})}{4(\Q k_{1}^{2}\alpha_{2}^{2}+4\alpha_{j_{1}k_{1}}^{4})}+\frac{\frac{l_{1}\pi\rho^2(j_{1}\alpha_{1} \Im  u_{J_{1}}^{21}-k_{1}\alpha_{2}\Im  u_{J_{1}}^{11})}{\Q k_{1}^{2}\alpha_{2}^{2}+4\alpha_{j_{1}k_{1}}^{4}}}{(\Q k_{1}^{2}\alpha_{2}^{2}+4\alpha_{j_{1}k_{1}}^{4})^{2}+4\alpha_{j_{1}k_{1}}^2 \rho^2}
				\\
					   & +\frac{l_{1}\pi \rho(k_{1}\alpha_{2}\Re  u_{J_{1}}^{11}-j_{1}\alpha_{1}\Re  u_{J_{1}}^{21})}{2 [(\Q k_{1}^{2}\alpha_{2}^{2}+4\alpha_{j_{1}k_{1}}^{4})^{2}+4\alpha_{j_{1}k_{1}}^{2}\rho^2]}, \\
				C_{4}= {} & -\frac{\Pr \Im  \theta_{J_{1}}^{1}}{8l_{1}\pi}+\frac{\Pr ^{3}\rho^{2} \Im \theta_{J_{1}}^{1}}{2l_{1}\pi(16l_{1}^{4}\pi^4+4\Pr ^{2}\rho^{2})}+\frac{l_{1}\pi\rho \Pr ^{2}\Re  \theta_{J_{1}}^{1}}{16l_{1}^{2}\pi^{2}+4\Pr ^{2}\rho^{2}},
				\\
				C_{5}= {} & \frac{\frac{l_{1}\pi\rho^2(k_{1}\alpha_{2} \Re  u_{J_{1}}^{11}-j_{1}\alpha_{1}\Re  u_{J_{1}}^{21})}{2(\Q k_{1}^{2}\alpha_{2}^{2}+4\alpha_{j_{1}k_{1}}^{4})}}{(\Q k_{1}^{2}\alpha_{2}^{2}+4\alpha_{j_{1}k_{1}}^{2})^{2}+4\alpha_{j_{1}k_{1}}^4 \rho^2}+\frac{l_{1}\pi \rho(j_{1}\alpha_1\Im  u_{J_{1}}^{21}-k_{1}\alpha_{2}\Im  u_{J_{1}}^{11})}{4 [(\Q k_{1}^{2}\alpha_{2}^{2}+4\alpha_{j_{1}k_{1}}^{2})^{2}+4\alpha_{j_{1}k_{1}}^{4}\rho^2]},
				\\
				C_{6}= {} & -\frac{\Pr ^{3}\rho^{2} \Re \theta_{J_{1}}^{1}}{4l_{1}\pi(16l_{1}^{4}\pi^4+4\Pr ^{2}\rho^{2})}+\frac{l_{1}\pi \Pr ^{2}\rho \Im  \theta_{J_{1}}^{1}}{32l_{1}^{2}\pi^{2}+8\Pr ^{2}\rho^{2}}.
			\end{align*}
			
			According to the formula (\cite{Wiggins2003} Page 385), the types of transition of the system \eqref{abstract1} can be determined by the sign of the nondimensional number 
			\[
			a(\Ra)=3a_{111}(\Ra)+a_{122}(\Ra)+b_{112}(\Ra).
			\]
	  According to \cite{Ma1}, when $a(Ra_{c_{2}})<0$, the bifurcated periodic solution are given by
			\begin{align}\label{zhouqiguidao1}
				\begin{aligned}
					&\Psi_p=\left(\frac{4\sigma}{-\pi a(Ra_{c_{2}})}\right)^{\frac{1}{2}}\sin(\rho t)\Psi_{J_{1}}^{1}+\left(\frac{4\sigma}{-\pi a(Ra_{c_{2}})}\right)^{\frac{1}{2}}\cos(\rho t)\Psi_{J_{1}}^{2}+o(\sqrt{|\sigma|}).
				\end{aligned}
			\end{align}
	\end{proof}

	\section{Transition from first real eigenvalue with double multiplicity}
\label{section5}
We consider a little more complex situation where 
$\text{card}(X)=2$, $\Ra_c=\Ra_{c_1}$ and 
\[
X=\{J_{2},J_{3}\},\quad J_{2}=(j_{2},k_{2},1),\quad J_{3}=(j_{2},-k_{2},1),
\]
i.e. the first eigenvalues are $\beta_{J_{2}}^{1}$ and $\beta_{J_{3}}^{1}$. One gets  from the equation \eqref{characteristicequation1} that $\beta_{J_{2}}^{1}=\beta_{J_{3}}^{1}$. That is, the first eigenvalue $\beta_{J_{2}}^{1}=\beta_{J_{3}}^{1}$
has multiplicity two.  $\Psi_{J_{2}}^{1}$ and $\Psi_{J_{3}}^{1}$ are the corresponding two eigenvectors,  and $\Psi_{J_{2}}^{1*}$ and $\Psi_{J_{3}}^{1*}$ are the conjugate eigenvectors. A direct verification shows that $
\left\langle \Psi_{J_{2}}^{1},\Psi_{J_{2}}^{1*}\right\rangle=\left\langle \Psi_{J_{3}}^{1},\Psi_{J_{3}}^{1*}\right\rangle$.
Thus, for simplicity, let us denote $\beta=\beta_{J_{2}}^{1}=\beta_{J_{3}}^{1}$ and $q=\left\langle \Psi_{J_{2}}^{1},\Psi_{J_{2}}^{1*}\right\rangle$, we then have the following lemma.

\begin{lemma}\label{reducedodelll}
In the vicinity of $Ra=Ra_{c_{1}}$, the stability and transition of the system \eqref{abstract1} for any small initial condition is equivalent to these of the following ODEs:
\begin{align}\label{reducedODE3}
    \frac{d\mathbf{Y}}{dt}=\beta \mathbf{Y} + g(\mathbf{Y}),\quad  \mathbf{Y}=(y,z)^{T},
\end{align}
where $g(\mathbf{Y})=(y(\Gamma_{1}y^2+\Gamma_{2}z^{2}), z(\Gamma_{3}y^2+\Gamma_{1}z^{2}))^{T}$,
    and $\Gamma_{1}$, $\Gamma_{2}$ and $\Gamma_{3}$ are 
 three transition numbers determining the type of transition, whose expressions are
    given in \eqref{zuida1}.
\end{lemma}
\begin{proof}
    Likewise, the $\mathbf{H_{1}}$ can be decomposed as follows:
    \begin{align}\label{decompositionofspace3}
        \mathbf{H_{1}}=\overline{\overline{E_{2}}}+\overline{\overline{E_{3}}}, 
    \end{align}
    where $\overline{\overline{E_{2}}}=\operatorname*{span}\{\Psi_{J_{2}}^{1},\Psi_{J_{3}}^{1}\}$ and $\overline{\overline{E_{3}}}$ is spanned by the rest of eigenvectors. As before, there exists a center manifold function $\Phi$ $:\overline{\overline{E_{2}}}\rightarrow \overline{\overline{E_{3}}}$ at $Ra=Ra_{c_{1}}$. Hence, on the center manifold, the solution $\Psi\in \mathbf{H_{1}}$ of \eqref{abstract1} has the following spectral decomposition,
    \begin{align}\label{decompositonoffunctions3}
        \Psi=y\Psi_{J_{2}}^{1}+z\Psi_{J_{3}}^{1}+\Phi,~\Phi=\sum\limits_{(J,s)\neq(J_{m},1)}x_{J}^{s}\Psi_{J}^{s}\in \overline{\overline{E_{3}}}~(m=2,3).
    \end{align}

In a similar way, we can get the following ODEs:
    \begin{align}\label{reducedequation3}
        \begin{aligned}
            \frac{dy}{dt}=\beta y+\frac{\left\langle \mathbf{G}(y\Psi_{J_{2}}^{1}+z\Psi_{J_{3}}^{1}+\Phi),\Psi_{J_{2}}^{1*}\right\rangle }{\left\langle \Psi_{J_{2}}^{1}, \Psi_{J_{2}}^{1*}\right\rangle },
            \\
            \frac{dz}{dt}=\beta z+\frac{\left\langle \mathbf{G}(y\Psi_{J_{2}}^{1}+z\Psi_{J_{3}}^{1}+\Phi),\Psi_{J_{3}}^{1*}\right\rangle }{\left\langle \Psi_{J_{3}}^{1}, \Psi_{J_{3}}^{1*}\right\rangle },
        \end{aligned}
    \end{align}
    where 
    \begin{align}\label{innerproduct2}
        q=\left\langle \Psi_{J_{m}}^{1},\Psi_{J_{m}}^{1*} \right\rangle =\frac{\pi^2 (u_{J_{m}}^{11}u_{J_{m}}^{11*}+u_{J_{m}}^{21}u_{J_{m}}^{21*}+1+\theta_{J_{m}}^{1}\theta_{J_{m}}^{1*})}{\alpha_{1}\alpha_{2}}~~~~(m=2,3).
    \end{align}
    
A direct computation gives
    \begin{align}\label{Gbilinearform3}
        \begin{aligned}
            \mathbf{G}(y\Psi_{J_{2}}^{1}+z\Psi_{J_{3}}^{1}+\Phi)  
            ={} & G(\Phi,y\Psi_{J_{2}}^{1}+z\Psi_{J_{3}}^{1})+G(y\Psi_{J_{2}}^{1}+z\Psi_{J_{3}}^{1},\Phi) \\
                & +G(\Phi,y\Psi_{J_{2}}^{1}+z\Psi_{J_{3}}^{1})+G(\Phi,\Phi),
        \end{aligned}
    \end{align}
    and 
    \begin{align}
        \left\langle  \mathbf{G}(y\Psi_{J_{2}}^{1}+z\Psi_{J_{3}}^{1}),\Psi_{J_{m}}^{1*}\right\rangle=0~~(m=2,3). 
    \end{align}
    Furthermore,  utilizing the formula (A.1.14) (\cite{Ma1}), we have 
    \begin{align}\label{centermanifold3}
        L\Phi=-y^{2}G_{22}-yz(G_{23}+G_{32})-z^{2}G_{33}+o,
    \end{align}
    where $G_{mn}=G(\Psi_{J_{m}}^{1},\Psi_{J_{n}}^{1})(m,n=2,3)$ are independent of $y$ and $z$, and 
    \begin{align}\label{G22}
        G_{22}      & =\frac{\pi}{2} \left[  \begin{pmatrix}
                \mathbf{P}\begin{pmatrix}
                    u_{J_{2}}^{11}\sin{2(j_{2}\alpha_{1}x_{1}+k_{2}\alpha_{2}x_{2})}
                    \\
                    u_{J_{2}}^{21}\sin{2(j_{2}\alpha_{1}x_{1}+k_{2}\alpha_{2}x_{2})}
                    \\
                    0
                \end{pmatrix}
                \\
                0
            \end{pmatrix}-
            \begin{pmatrix}
                0
                \\
                0
                \\
                0
                \\
                \theta_{J_{2}}^{1}\sin{2\pi x_{3}}
            \end{pmatrix}
            \right],
        \\ 
        \label{G33}
        ~~~~ G_{33} & =\frac{\pi}{2} \left[  \begin{pmatrix}
                \mathbf{P}\begin{pmatrix}
                    u_{J_{3}}^{11}\sin{2(j_{2}\alpha_{1}x_{1}-k_{2}\alpha_{2}x_{2})}
                    \\
                    u_{J_{3}}^{21}\sin{2(j_{2}\alpha_{1}x_{1}-k_{2}\alpha_{2}x_{2})}
                    \\
                    0
                \end{pmatrix}
                \\
                0
            \end{pmatrix}-
            \begin{pmatrix}
                0
                \\
                0
                \\
                0
                \\
                \theta_{J_{3}}^{1}\sin{2\pi x_{3}}
            \end{pmatrix}
            \right],
    \end{align}
    
    \begin{align}\label{G23}
        \begin{aligned}
            G_{23} = {} & P_{2}\begin{pmatrix}
                \mathbf{P}\begin{pmatrix}
                    u_{J_{3}}^{11}\sin{2j_{2}\alpha_{1}x_{1}}\cos{2\pi x_{3}}
                    \\
                    u_{J_{3}}^{21}\sin{2j_{2}\alpha_{1}x_{1}}\cos{2\pi x_{3}}
                    \\
                    \cos{2j_{2}\alpha_{1}x_{1}}\sin{2\pi x_{3}}
                \end{pmatrix}
                \\
                \theta_{J_{3}}^{1}\cos{2j_{2}\alpha_{1}x_{1}}\sin{2\pi x_{3}}
            \end{pmatrix}
            -P_{1} \begin{pmatrix}
                \mathbf{P}\begin{pmatrix}
                    u_{J_{3}}^{11}\sin{2j_{2}\alpha_{1}x_{1}}
                    \\
                    u_{J_{3}}^{21}\sin{2j_{2}\alpha_{1}x_{1}}
                    \\
                    0
                \end{pmatrix}
                \\
                0
            \end{pmatrix}
            \\        
                     & +P_{2}\begin{pmatrix}
                \mathbf{P}\begin{pmatrix}
                    u_{J_{3}}^{11}\sin{2k_{2}\alpha_{2}x_{2}}
                    \\
                    u_{J_{3}}^{21}\sin{2k_{2}\alpha_{2}x_{2}}
                    \\
                    0
                \end{pmatrix}
                \\
                0
            \end{pmatrix}
            -P_{1}\begin{pmatrix}
                \mathbf{P}\begin{pmatrix}
                    u_{J_{3}}^{11}\sin{2k_{2}\alpha_{2}x_{2}}\cos{2\pi x_{3}}
                    \\
                    u_{J_{3}}^{21}\sin{2k_{2}\alpha_{2}x_{2}}\cos{2\pi x_{3}}
                    \\
                    -\cos{2k_{2}\alpha_{2}x_{2}}\sin{2\pi x_{3}}
                \end{pmatrix}
                \\
                -\theta_{J_{3}}^{1}\cos{2k_{2}\alpha_{2}x_{2}}\sin{2\pi x_{3}}
            \end{pmatrix},
        \end{aligned}
    \end{align}
    
    \begin{align}\label{G32}
        \begin{aligned}
            G_{32}  = {} & \widetilde{P_{1}} \begin{pmatrix}
                \mathbf{P}\begin{pmatrix}
                    u_{J_{2}}^{11}\sin{2k_{2}\alpha_{2}x_{2}}\cos{2\pi x_{3}}
                    \\
                    u_{J_{2}}^{21}\sin{2k_{2}\alpha_{2}x_{2}}\cos{2\pi x_{3}}
                    \\
                    \cos{2k_{2}\alpha_{2}x_{2}}\sin{2\pi x_{3}}
                \end{pmatrix}
                \\
                \theta_{J_{2}}^{1}\cos{2k_{2}\alpha_{2}x_{2}}\sin{2\pi x_{3}}
            \end{pmatrix}
            -\widetilde{P_{1}}  \begin{pmatrix}
                \mathbf{P}\begin{pmatrix}
                    u_{J_{2}}^{11}\sin{2j_{2}\alpha_{1}x_{1}}
                    \\
                    u_{J_{2}}^{21}\sin{2j_{2}\alpha_{1}x_{1}}
                    \\
                    0
                \end{pmatrix}
                \\
                0
            \end{pmatrix}
            \\
                      & +\widetilde{P_{2}} \begin{pmatrix}
                \mathbf{P}\begin{pmatrix}
                    u_{J_{2}}^{11}\sin{2k_{2}\alpha_{2}x_{2}}
                    \\
                    u_{J_{2}}^{21}\sin{2k_{2}\alpha_{2}x_{2}}
                    \\
                    0
                \end{pmatrix}
                \\
                0
            \end{pmatrix}
            -\widetilde{P_{2}} \begin{pmatrix}
                \mathbf{P}\begin{pmatrix}
                    u_{J_{2}}^{11}\sin{2j_{2}\alpha_{1}x_{1}}\cos{2\pi x_{3}}
                    \\
                    u_{J_{2}}^{21}\sin{2j_{2}\alpha_{1}x_{1}}\cos{2\pi x_{3}}
                    \\
                    \cos{2j_{2}\alpha_{1}x_{1}}\sin{2\pi x_{3}}
                \end{pmatrix}
                \\
                \theta_{J_{2}}^{1}\cos{2j_{2}\alpha_{1}x_{1}}\sin{2\pi x_{3}}
            \end{pmatrix},
        \end{aligned}
    \end{align}
in which
    \begin{align*}
        \begin{aligned}
            P_{1}=\frac{j_{2}\alpha_{1}u_{J_{2}}^{11}}{2},~~~~P_{2}=\frac{k_{2}\alpha_{2}u_{J_{2}}^{21}}{2},~~~
            \widetilde{P_{1}} =\frac{j_{2}\alpha_{1}u_{J_{3}}^{11}}{2},~~~~\widetilde{P_{2}} =\frac{k_{2}\alpha_{2}u_{J_{3}}^{21}}{2}.
        \end{aligned}
    \end{align*}

  According to \eqref{G22}-\eqref{G32} and  for the convenience of computing $x_{J}^{s}$, we introduce the following notations:
    \begin{align}\label{notations}
        \begin{aligned}
             & K_{1}=(0,0,2),~K_{2}=(2j_{2},0,0),~K_{3}=(0,2k_{2},0),~K_{4}=(2j_{2},2k_{2},0),
            \\
             & K_{5}=(2j_{2},-2k_{2},0),~K_{6}=(2j_{2},0,2),~K_{7}=(0,2k_{2},2).
        \end{aligned}
    \end{align}
It is known from the section \ref{section2} that $\beta_{K_{1}}^{3}$ and $\beta_{K_{m}}^{s}~(m=2,\cdots,5;~s=1,2,3)$ are real. There are two cases in the process of computing the center manifold function $\Phi$. 

Case i: $\beta_{K_{m}}^{s}(m=6,7;~s=1,2,3)$ are all real.
Comparing the coefficients on both sides of \eqref{centermanifold3}, we have
    \begin{align}\label{computation1}
        x_{J}^{s}=-\frac{\left\langle y^{2}G_{22}+yz(G_{23}+G_{32})+z^{2}G_{33},\Psi_{J}^{s*}\right\rangle }{\beta_{J}^{s}\left\langle \Psi_{J}^{s}, \Psi_{J}^{s*}\right\rangle},~~J\in \{K_{1},\cdots,K_{7}\}.
    \end{align}
By some calculation we have:
    \begin{align}\label{xishu3}
        \begin{cases}
            x_{K_{1}}^{3}=\frac{\pi\theta_{J_{2}}^{1}y^{2}+\pi\theta_{J_{3}}^{1}z^{2}}{2\beta_{K_{1}}^{3}},~~x_{K_{2}}=-\frac{(P_{1}u_{J_{3}}^{21}+\widetilde{P_{1}} u_{J_{2}}^{21})yz}{2j_{2}\alpha_{1}\beta_{K_{2}}},~~x_{K_{3}}=-\frac{(P_{2}u_{J_{3}}^{11}+\widetilde{P_{2}}u_{J_{2}}^{11})yz}{2k_{2}\alpha_{2}\beta_{K_{3}}},
            \\
            x_{K_{4}}=\frac{\pi(j_{2}\alpha_{1}u_{J_{2}}^{21}-k_{2}\alpha_{2}u_{J_{2}}^{11})y^{2}}{4\alpha_{j_{2}k_{2}}^{2}\beta_{K_{4}}},~~x_{K_{5}}=\frac{\pi (j_{2}\alpha_{1}u_{J_{3}}^{21}+k_{2}\alpha_{2}u_{J_{3}}^{11})z^{2}}{4\alpha_{j_{2}k_{2}}^{2}\beta_{K_{5}}},
            \\
            x_{K_{6}}^{s}=-\frac{ [(P_{2}u_{J_{3}}^{11}-\widetilde{P_{2}}u_{J_{2}}^{11})u_{k_{6}}^{1s*}+(P_{2}u_{J_{3}}^{21}-\widetilde{P_{2}}u_{J_{2}}^{21})u_{K_{6}}^{2s*}+(P_{2}-\widetilde{P_{2}})+(P_{2}\theta_{J_{3}}^{1}-\widetilde{P_{2}}\theta_{J_{2}}^{1})\theta_{K_{6}}^{s*}]yz}{(u_{K_{6}}^{1s}u_{K_{6}}^{1s*}+u_{K_{6}}^{2s}u_{K_{6}}^{2s*}+1+\theta_{K_{6}}^{s}\theta_{K_{6}}^{s*})\beta_{K_{6}}^{s}},
            \\
            x_{K_{7}}^{s}=\frac{ [(P_{1}u_{J_{3}}^{11}-\widetilde{P_{1}}u_{J_{2}}^{11})u_{k_{7}}^{1s*}+(P_{1}u_{J_{3}}^{21}-\widetilde{P_{1}}u_{J_{2}}^{21})u_{K_{7}}^{2s*}-(P_{1}+\widetilde{P_{1}})-(P_{1}\theta_{J_{2}}^{1}+\widetilde{P_{1}}\theta_{J_{3}}^{1})\theta_{K_{7}}^{s*}]yz}{(u_{K_{7}}^{1s}u_{K_{7}}^{1s*}+u_{K_{7}}^{2s}u_{K_{7}}^{2s*}+1+\theta_{K_{7}}^{s}\theta_{K_{7}}^{s*})\beta_{K_{7}}^{s}}.
        \end{cases}
    \end{align}
   Therefore, the center manifold function $\Phi$ has the following expansion
    \begin{align}\label{centermanifoldfunction}
        \begin{aligned}
            \Phi= {} & (G_{K_{1}}^{3}\Psi_{K_{1}}^{3}+G_{K_{4}}\Psi_{K_{4}})y^{2}+(G_{K_{1}}^{3}\Psi_{K_{1}}^{3}+G_{K_{5}}\Psi_{K_{5}})z^2
            \\
                  & - [(G_{K_{2}}\Psi_{K_{2}}+G_{K_{3}}\Psi_{K_{3}})+\sum_{s=1}^{3}(G_{K_{6}}^{s}\Psi_{K_{6}}^{s}-G_{K_{7}}^{s}\Psi_{K_{7}}^{s})]yz+o,
        \end{aligned}
    \end{align}
    where for $s=1,2,3$
    \begin{align}\label{2021110303}
        \begin{cases}
            G_{K_{1}}^{3}=\frac{\pi\theta_{J_{2}}^{1}}{2\beta_{K_{1}}^{3}},~~G_{K_{4}}=\frac{\pi(j_{2}\alpha_{1}u_{J_{2}}^{21}-k_{2}\alpha_{2}u_{J_{2}}^{11})}{4\alpha_{j_{2}k_{2}}^{2}\beta_{K_{4}}},
            \\
            G_{K_{5}}=\frac{\pi(j_{2}\alpha_{1}u_{J_{3}}^{21}+k_{2}\alpha_{2}u_{J_{3}}^{11})}{4\alpha_{j_{2}k_{2}}^{2}
            \beta_{K_{5}}},~~G_{K_{2}}=\frac{u_{J_{2}}^{11}u_{J_{3}}^{21}+u_{J_{2}}^{21}u_{J_{3}}^{11}}{4\beta_{K_{2}}},~~G_{K_{3}}=\frac{u_{J_{2}}^{21}u_{J_{3}}^{11}+u_{J_{2}}^{11}u_{J_{3}}^{21}}{4\beta_{K_{3}}},
            \\
            G_{K_{6}}^{s}=\frac{(P_{2}u_{J_{3}}^{11}-\widetilde{P_{2}}u_{J_{2}}^{11})u_{K_{6}}^{1s*}+(P_{2}u_{J_{3}}^{21}-\widetilde{P_{2}}u_{J_{2}}^{21})u_{K_{6}}^{2s*}+(P_{2}-\widetilde{P_{2}})+(P_{2}\theta_{J_{3}}^{1}-\widetilde{P_{2}}\theta_{J_{2}}^{1})\theta_{K_{6}}^{s*}}{(u_{K_{6}}^{1s}u_{K_{6}}^{1s*}+u_{K_{6}}^{2s}u_{K_{6}}^{2s*}+1+\theta^{s}_{K_{6}}\theta_{K_{6}}^{s*})\beta_{K_{6}}^{s}},
            \\
            G_{K_{7}}^{s}=\frac{(P_{1}u_{J_{3}}^{11}-\widetilde{P_{1}}u_{J_{2}}^{11})u_{K_{7}}^{1s*}+(P_{1}u_{J_{3}}^{21}-\widetilde{P_{1}}u_{J_{2}}^{21})u_{K_{7}}^{2s*}-(P_{1}+\widetilde{P_{1}})-(P_{1}\theta_{J_{2}}^{1}+\widetilde{P_{1}}\theta_{J_{3}}^{1})\theta_{K_{7}}^{s*}}{(u_{K_{7}}^{1s}u_{K_{7}}^{1s*}+u_{K_{7}}^{2s}u_{K_{7}}^{2s*}+1+\theta_{K_{7}}^{s}\theta_{K_{7}}^{s*})\beta_{K_{7}}^{s}}.
    \end{cases}
    \end{align}

Case ii: there exists one $m\in\{6,7\}$ such that $\beta_{K_{m}}^{1}$$=$$\overline{\beta_{K_{m}}^{2}}$$=$$\sigma_{m}+i\rho_{m}$ with $\rho_{m}\neq 0$. We can compute the $x_{K_{m}}^{s}(s=1,2)$ by employing the following method. Firstly, we let
    \begin{align}\label{weilehuajian}
        \begin{aligned}
        &\widetilde{\Psi_{K_{m}}^{1*}}=\Psi_{K_{m}}^{1*}+\widetilde{\alpha}  \Psi_{K_{m}}^{2*}
        \\
        &\widetilde{\Psi_{K_{m}}^{1*}}=-\widetilde{\alpha} \Psi_{K_{m}}^{1*}+\Psi_{K_{m}}^{2*}
        \end{aligned} ,
    \end{align} 
    where $\widetilde{\alpha} =\frac{\left\langle \Psi_{K_{m}}^{1},\Psi_{K_{m}}^{2*}\right\rangle }{\left\langle \Psi_{K_{m}}^{1},\Psi_{K_{m}}^{1*}\right\rangle }$, $\Psi_{K_{m}}^{s*}(s=1,2)$ and $\Psi_{K_{m}}^{s}(s=1,2)$  are in \eqref{conjugateeigenvectorI11} and \eqref{EigenvectorI11}, respectively. For $p,q=1,2$, we have
    \begin{align}\label{zhengjiaoguanxi}
        \left\langle\Psi_{K_{m}}^{p}, \widetilde{\Psi_{K_{m}}^{q*}}\right\rangle 
        \begin{cases}
            =0 ~~(p\neq q)
            \\
            \neq 0~~(p=q)
        \end{cases}.
    \end{align}
    Conducting some computation, we have 
     \begin{align}\label{jutisuan}
        \begin{cases}
            x_{K_{m}}^{1}= [\frac{\rho_{m}\left\langle G_{23}+G_{32},\widetilde{\Psi_{K_{m}}^{2*}} \right\rangle }{|\beta_{K_{m}}^{1}|^2\left\langle \Psi_{K_{m}}^{2},\widetilde{\Psi_{K_{m}}^{2*}}\right\rangle}-\frac{\sigma_{m}\left\langle G_{23}+G_{32},\widetilde{\Psi_{K_{m}}^{1*}} \right\rangle}{|\beta_{K_{m}}^{1}|^2\left\langle \Psi_{K_{m}}^{1},\widetilde{\Psi_{K_{m}}^{1*}}\right\rangle}]yz,
            \\
            x_{K_{m}}^{2}= [-\frac{\rho_{m}\left\langle G_{23}+G_{32},\widetilde{\Psi_{K_{m}}^{1*}} \right\rangle }{|\beta_{K_{m}}^{1}|^2\left\langle \Psi_{K_{m}}^{1},\widetilde{\Psi_{K_{m}}^{1*}}\right\rangle}-\frac{\sigma_{m}\left\langle G_{23}+G_{32},\widetilde{\Psi_{K_{m}}^{2*}} \right\rangle}{|\beta_{K_{m}}^{1}|^2\left\langle \Psi_{K_{m}}^{2},\widetilde{\Psi_{K_{m}}^{2*}}\right\rangle}]yz,
        \end{cases}
     \end{align}
by which one has
    \begin{align}\label{jisuan}
        \begin{cases}
            G_{K_{m}}^{1}=\frac{\rho_{m}\left\langle G_{23}+G_{32},\widetilde{\Psi_{K_{m}}^{2*}} \right\rangle }{ |\beta_{K_{m}}^{1}|^2\left\langle \Psi_{K_{m}}^{2},\widetilde{\Psi_{K_{m}}^{2*}}\right\rangle}-\frac{\sigma_{m}\left\langle G_{23}+G_{32},\widetilde{\Psi_{K_{m}}^{1*}} \right\rangle}{ |\beta_{K_{m}}^{1}|^2\left\langle \Psi_{K_{m}}^{1},\widetilde{\Psi_{K_{m}}^{1*}}\right\rangle},
            \\
            G_{K_{m}}^{2}=-\frac{\rho_{m}\left\langle G_{23}+G_{32},\widetilde{\Psi_{K_{m}}^{1*}} \right\rangle }{ |\beta_{K_{m}}^{1}|^2\left\langle \Psi_{K_{m}}^{1},\widetilde{\Psi_{K_{m}}^{1*}}\right\rangle}-\frac{\sigma_{m}\left\langle G_{23}+G_{32},\widetilde{\Psi_{K_{m}}^{2*}} \right\rangle}{ |\beta_{K_{m}}^{1}|^2\left\langle \Psi_{K_{m}}^{2},\widetilde{\Psi_{K_{m}}^{2*}}\right\rangle},
        \end{cases}
    \end{align}
    and the others in \eqref{centermanifoldfunction} remain the same, given in \eqref    
  {2021110303}.

    To sum up, when $\beta_{K_{m}}^{s}(m=6,7;~s=1,2,3)$ are all real, we take \eqref{2021110303} into \eqref{Gbilinearform3}; when one of $\beta_{K_{m}}^{1}(m=6,7)$ is complex, we substitute \eqref{jisuan} and the others given in 
    \eqref{2021110303}
     into \eqref{Gbilinearform3}. Finally, the system \eqref{reducedequation3} is simplified into 
    \begin{align}\label{reducechenweiode}
        \frac{d\mathbf{Y}}{dt}=\beta \mathbf{Y} + g(\mathbf{Y}),
    \end{align}
    where 
    \begin{align*}
             & \mathbf{Y}=(y,z)^{T},\quad g(\mathbf{Y})=(y(\Gamma_{1}y^2+\Gamma_{2}z^{2}), z(\Gamma_{3}y^2+\Gamma_{1}z^{2}))^{T},
    \end{align*}
    \begin{align}\label{zuida1}
       & \begin{aligned}
             & \Gamma_{1}=\frac{\pi^2 [\pi G_{K_{1}}^{3}\theta_{J_{2}}^{1*}+2\pi G_{K_{4}}(j_{2}\alpha_{1}u_{J_{2}}^{21*}-k_{2}\alpha_{2}u_{J_{2}}^{11*})]}{q\alpha_{1}\alpha_{2}},
            \\
             & 
             \begin{aligned}\Gamma_{2}={}& \frac{\pi^2}{q\alpha_{1}\alpha_{2}} \Big[G_{K_{1}}^{3}\pi\theta_{J_{2}}^{1*}+G_{K_{2}K_{3}1}+G_{K_{2}K_{3}2}  \\
                &  -\sum_{s=1}^{3} ((G_{K_{61}}^{s}+G_{K_{62}}^{s})G_{K_{6}}^{s}+(G_{K_{71}}^{s}+G_{K_{72}}^{s})G_{K_{7}}^{s}) \Big], 
             \end{aligned}
            \\
             & 
             \begin{aligned}\Gamma_{3}= {}&\frac{\pi^2}{q\alpha_{1}\alpha_{2}} \Big[G_{K_{1}}^{3}\pi\theta_{J_{3}}^{1*}+G_{K_{2}K_{3}3}-G_{K_{2}K_{3}4}
                \\
                &+\sum_{s=1}^{3} ((G_{K_{63}}^{s}+G_{K_{64}}^{s})G_{K_{6}}^{s}+(G_{K_{73}}^{s}+G_{K_{74}}^{s})G_{K_{7}}^{s}) \Big],
             \end{aligned}
        \end{aligned}
    \end{align}
    and 
    \begin{align}\label{110}
        &\begin{aligned}
         & G_{K_{2}K_{3}1}=2j_{2}^{2}\alpha_{1}^{2}u_{J_{3}}^{11}u_{J_{2}}^{21*}G_{K_{2}}+2k_{2}^{2}\alpha_{2}^{2}u_{J_{3}}^{21}u_{J_{2}}^{11*}G_{K_{3}},
        \\
         & G_{K_{2}K_{3}3}=2j_{2}^{2}\alpha_{1}^{2}u_{J_{2}}^{11}u_{J_{3}}^{21*}G_{K_{2}}+2k_{2}^{2}\alpha_{2}^{2}u_{J_{2}}^{21}u_{J_{3}}^{11*}G_{K_{3}},
        \end{aligned}\\ \label{112}
       & \begin{aligned}
        &\begin{aligned}
         & G_{K_{62}}^{s}=2^{-1} k_{2}\alpha_{2}u_{J_{3}}^{21} (u_{K_{6}}^{1s}u_{J_{2}}^{11*}+u_{K_{6}}^{2s}u_{J_{2}}^{21*}+1+\theta_{K_{6}}^{s}\theta_{J_{2}}^{1*}),
        \\
         & G_{K_{63}}^{s}=2^{-1} k_{2}\alpha_{2}u_{J_{2}}^{21} (u_{K_{6}}^{1s}u_{J_{3}}^{11*}+u_{K_{6}}^{2s}u_{J_{3}}^{21*}+1+\theta_{K_{6}}^{s}\theta_{J_{3}}^{1*}),
        \end{aligned}
  \\   
        &\begin{aligned}
         & G_{K_{71}}^{s}=2^{-1}j_{2}\alpha_{1}u_{J_{3}}^{11} (u_{K_{7}}^{1s}u_{J_{2}}^{11*}+u_{K_{7}}^{2s}u_{J_{2}}^{21*}+1+\theta_{K_{7}}^{s}\theta_{J_{2}}^{1*}),
        \\
         & G_{K_{73}}^{s}=2^{-1}j_{2}\alpha_{1}u_{J_{2}}^{11} (u_{K_{7}}^{1s}u_{J_{3}}^{11*}+u_{K_{7}}^{2s}u_{J_{3}}^{21*}-1-\theta_{K_{7}}^{s}\theta_{J_{3}}^{1*}),
        \end{aligned}
    \end{aligned}\\
    &\label{111}
        \begin{aligned}
         &\begin{aligned} 
         G_{K_{2}K_{3}2} &= {} j_{2}\alpha_{1}k_{2}\alpha_{2}G_{K_{2}} (u_{J_{3}}^{11}u_{J_{2}}^{11*}+u_{J_{3}}^{21}u_{J_{2}}^{21*}  -1  - \theta_{J_{3}}\theta_{J_{2}}^{1*})
         \\
         &+j_{2}\alpha_{1}k_{2}\alpha_{2}G_{K_{3}} (u_{J_{3}}^{11}u_{J_{2}}^{11*}+u_{J_{3}}^{21}u_{J_{2}}^{21*}+1+\theta_{J_{3}}\theta_{J_{2}}^{1*}),
         \end{aligned}
         \\
         &\begin{aligned}
          G_{K_{2}K_{3}4}&={}j_{2}\alpha_{1}k_{2}\alpha_{2}G_{K_{2}} (u_{J_{2}}^{11}u_{J_{3}}^{11*}+u_{J_{2}}^{21}u_{J_{3}}^{21*}-1-\theta_{J_{2}}\theta_{J_{3}}^{1*})
         \\
         &+j_{2}\alpha_{1}k_{2}\alpha_{2}G_{K_{3}} (u_{J_{2}}^{11}u_{J_{3}}^{11*}+u_{J_{2}}^{21}u_{J_{3}}^{21*}+1+\theta_{J_{2}}\theta_{J_{3}}^{1*}),
         \end{aligned}
        \end{aligned}\\   
        \label{113}
      &  \begin{aligned}
       & \begin{aligned}
         & G_{K_{61}}^{s} = - 4^{-1}(j_{2}\alpha_{1}u_{K_{6}}^{1s}-k_{2}\alpha_{2}u_{K_{6}}^{2s}+\pi) (u_{J_{3}}^{11}u_{J_{2}}^{11*}+u_{J_{3}}^{21}u_{J_{2}}^{21*}+1+\theta_{J_{3}}^{1}\theta_{J_{2}}^{1*}),
        \\
         & G_{K_{64}}^{s}=4^{-1} (j_{2}\alpha_{1}u_{K_{6}}^{1s}+k_{2}\alpha_{2}u_{K_{6}}^{2s}+\pi) (u_{J_{2}}^{11}u_{J_{3}}^{11*}+u_{J_{2}}^{21}u_{J_{3}}^{21*}+1+\theta_{J_{2}}^{1}\theta_{J_{3}}^{1*}),
        \end{aligned}
  \\ 
       & \begin{aligned}
         & G_{K_{72}}^{s}= 4^{-1}(j_{2}\alpha_{1}u_{K_{7}}^{1s}-k_{2}\alpha_{2}u_{K_{7}}^{2s}-\pi) (u_{J_{3}}^{11}u_{J_{2}}^{11*}+u_{J_{3}}^{21}u_{J_{2}}^{21*}-1-\theta_{J_{3}}^{1}\theta_{J_{2}}^{1*}),
        \\
         & G_{K_{74}}^{s}=4^{-1}(j_{2}\alpha_{1}u_{K_{7}}^{1s}+k_{2}\alpha_{2}u_{K_{7}}^{2s}+\pi) (u_{J_{2}}^{11}u_{J_{3}}^{11*}+u_{J_{2}}^{21}u_{J_{3}}^{21*}-1-\theta_{J_{2}}^{1}\theta_{J_{3}}^{1*}).
        \end{aligned}
       \end{aligned}
    \end{align}
\end{proof}

 The possible non-zero equilibrium points of the ODEs \eqref{reducedODE3} are as follows:
\begin{align}\label{Odekenengpinghengdian}
    \mathbf{Y}= (0,\pm\sqrt{\gamma}),\ (\pm\sqrt{-\gamma},0),\ (\pm\sqrt{\xi},\pm\sqrt{\eta})
\end{align}
where
\[
\gamma=-\frac{\beta}{\Gamma_{1}},\quad 
\xi=\frac{\Gamma_{2}\beta-\Gamma_{1}\beta}{\Gamma_{1}^{2}-\Gamma_{2}\Gamma_{3}} ,\quad
\eta=\frac{\Gamma_{3}\beta-\Gamma_{1}\beta}{\Gamma_{1}^{2}-\Gamma_{2}\Gamma_{3}},
\quad 
\]
and which are assumed to be nonnegative.

In order to narrate easily, we make the following notations:
\begin{align}\label{biaohao}
    \begin{aligned}
     & &\mathbf{Y}_{1} &= -\mathbf{Y}_{2}= (0,\sqrt{\gamma}),\quad & & \mathbf{Y}_{5}   =-\mathbf{Y}_{8}= (\sqrt{\xi},\sqrt{\eta}),
     \\
    & & \mathbf{Y}_{3}& =-\mathbf{Y}_{4}=  (\sqrt{\gamma},0),\quad ,\quad & &\mathbf{Y}_{6} =-\mathbf{Y}_{7}   = (\sqrt{\xi},-\sqrt{\eta}),
    \end{aligned}
\end{align}
and 
\begin{equation}\label{pdeyuehuahoudejie}
    \begin{alignedat}{2}
        \Psi_{1} &=-  \Psi_{2} =\sqrt{\gamma}\Psi_{J_{3}}^1,\quad & \Psi_{5} &=-\Psi_{8}=\sqrt{\xi}\Psi_{J_{2}}^{1}+\sqrt{\eta}\Psi_{J_{3}}^{1}, \\
        \Psi_{3} &=-  \Psi_{4} =\sqrt{\gamma}\Psi_{J_{2}}^1,\quad &\Psi_{6}&=-\Psi_{7} =\sqrt{\xi}\Psi_{J_{2}}^{1}-\sqrt{\eta}\Psi_{J_{3}}^{1},
       \end{alignedat}
\end{equation}

\begin{figure}[H]
	\centering
	\begin{tabular}{cc}
		\includegraphics[width=.3\linewidth,height=.2\linewidth]{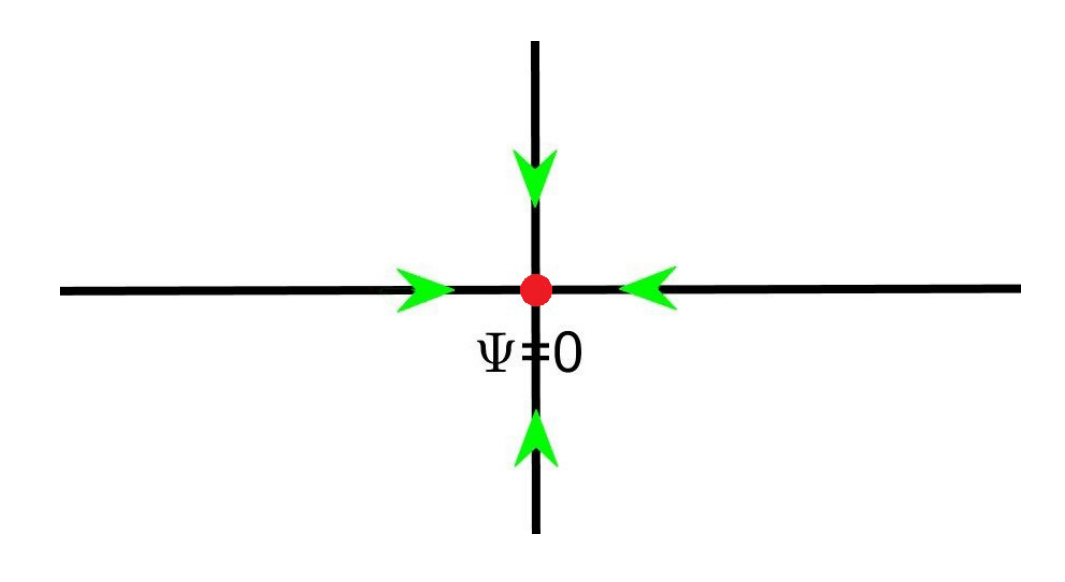}	~~~& \includegraphics[width=.3\linewidth,height=.2\linewidth]{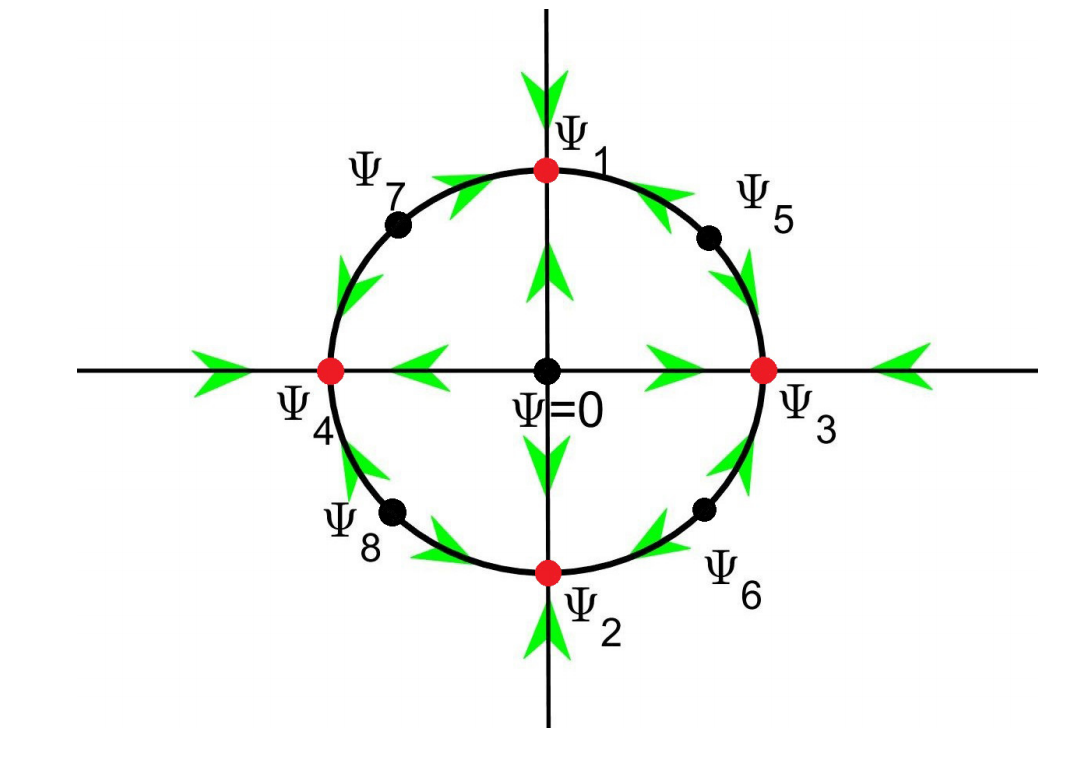}\\
				$(a)$~~~~& $(b)$
	\end{tabular}
	\caption{Topological structure of continuous transition of \eqref{abstract1}
	when $\Gamma_{1}<0$, $\Gamma_{2}\Gamma_{3}>\Gamma_{1}^{2}>0$, $\Gamma_{1}>\Gamma_{2}$ and $\Gamma_{1}>\Gamma_{3}$:
	(a) $\Ra< \Ra_{c_{1}}$; (b)  $\Ra> \Ra_{c_{1}}$.}  
	  \label{diyigedinglitu6dd}
\end{figure}

\begin{figure}[H]
	\centering
	\begin{tabular}{cc}
		\includegraphics[width=.3\linewidth,height=.2\linewidth]{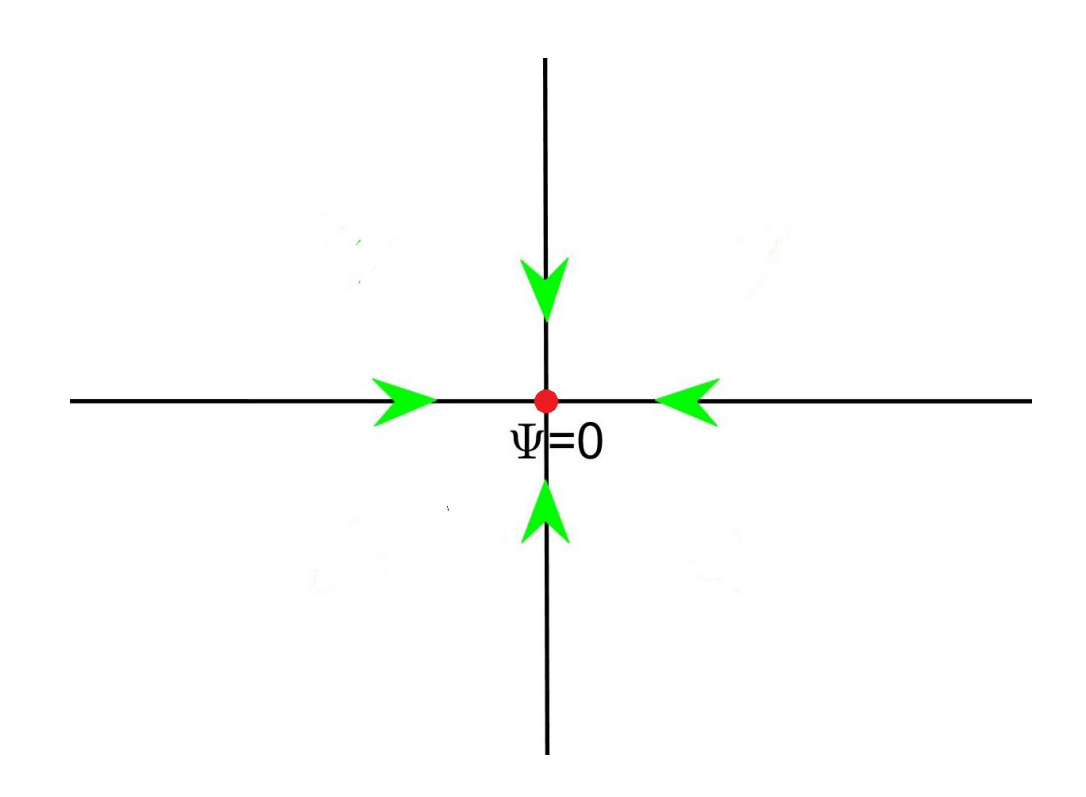}	~~~& \includegraphics[width=.3\linewidth,height=.2\linewidth]{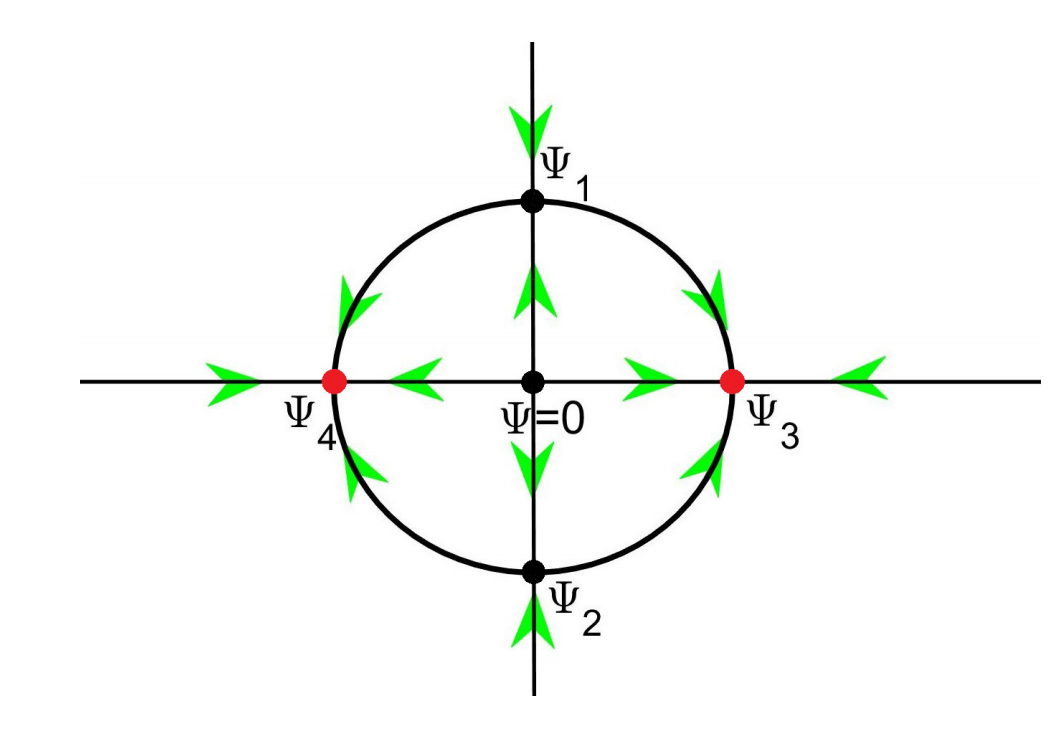}\\
				$(a)$~~~~& $(b)$
	\end{tabular}
	\caption{Topological structure of continuous transition of \eqref{abstract1}
	when $\Gamma_{1}<0$, $\Gamma_{1}<\Gamma_{2}$ and $\Gamma_{1}>\Gamma_{3}$:
 (a) $\Ra< \Ra_{c_{1}}$; (b)  $\Ra> \Ra_{c_{1}}$.}  
	  \label{diyigedinglitu156}
\end{figure}

\begin{figure}[H]
	\centering
	\begin{tabular}{cc}
		\includegraphics[width=.3\linewidth,height=.2\linewidth]{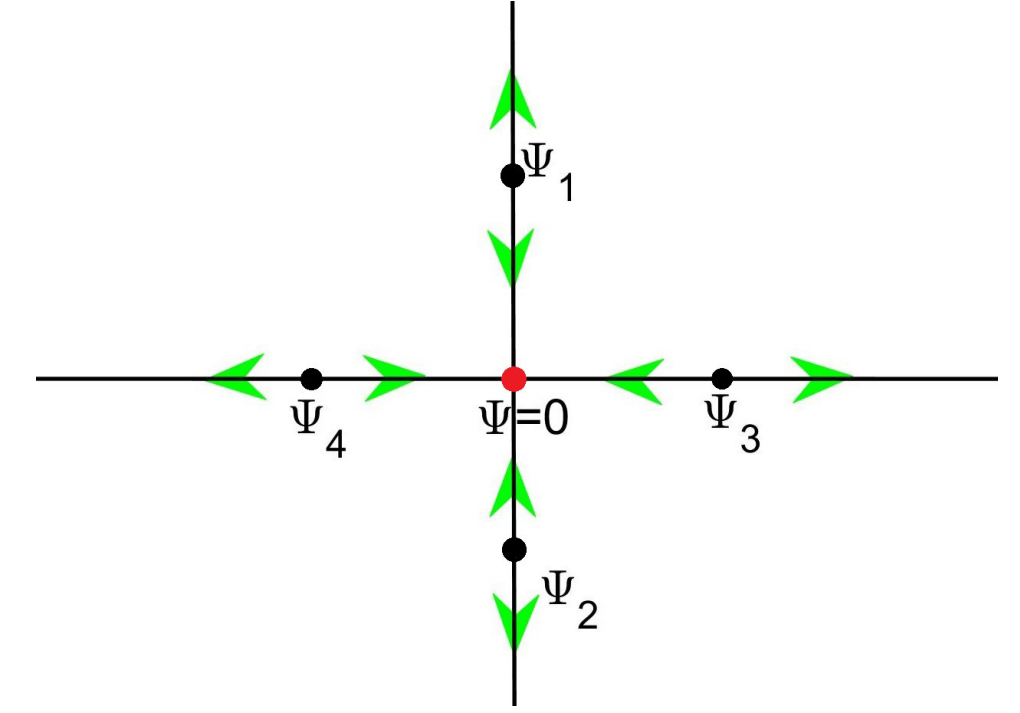}	~~~& \includegraphics[width=.3\linewidth,height=.2\linewidth]{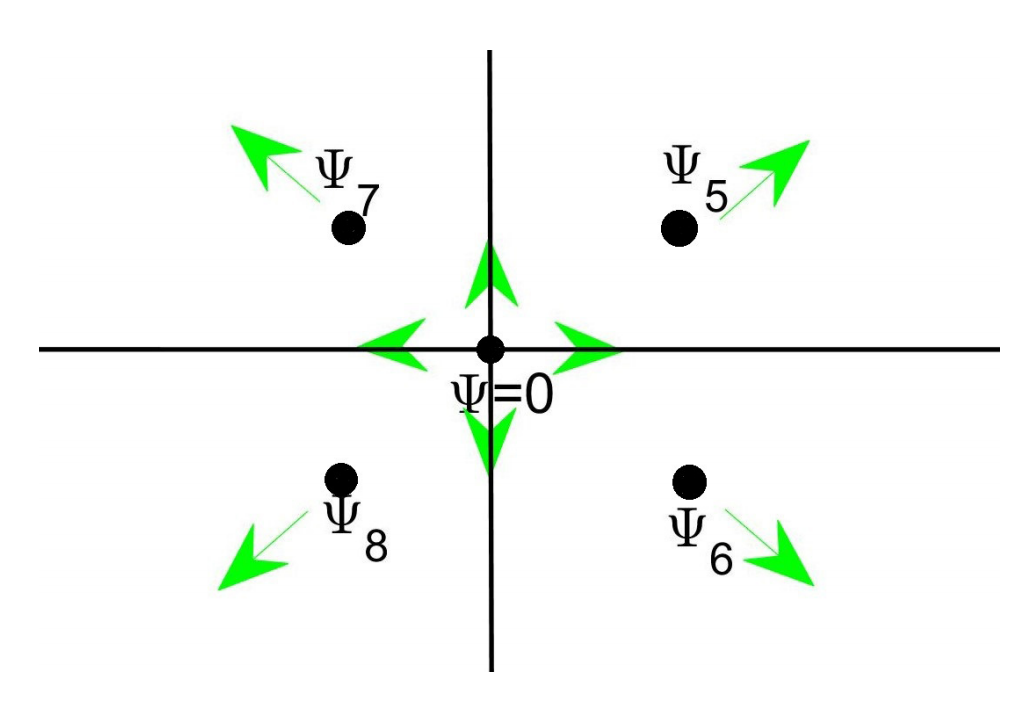}\\
				$(a)$~~~~& $(b)$
	\end{tabular}
	\caption{Topological structure of jump transition of \eqref{abstract1}
	when  $\Gamma_{1}>0$, $(\Gamma_{1}-\Gamma_{2})(\Gamma_{2}\Gamma_{3}-\Gamma_{1}^{2})>0$ and $(\Gamma_{1}-\Gamma_{3})(\Gamma_{2}\Gamma_{3}-\Gamma_{1}^{2})>0$:
 (a) $\Ra< \Ra_{c_{1}}$; (b)  $\Ra> \Ra_{c_{1}}$.}  
	  \label{diyigedinglitu1993}
\end{figure}

\begin{figure}[H]
	\centering
	\begin{tabular}{cc}
		\includegraphics[width=.3\linewidth,height=.2\linewidth]{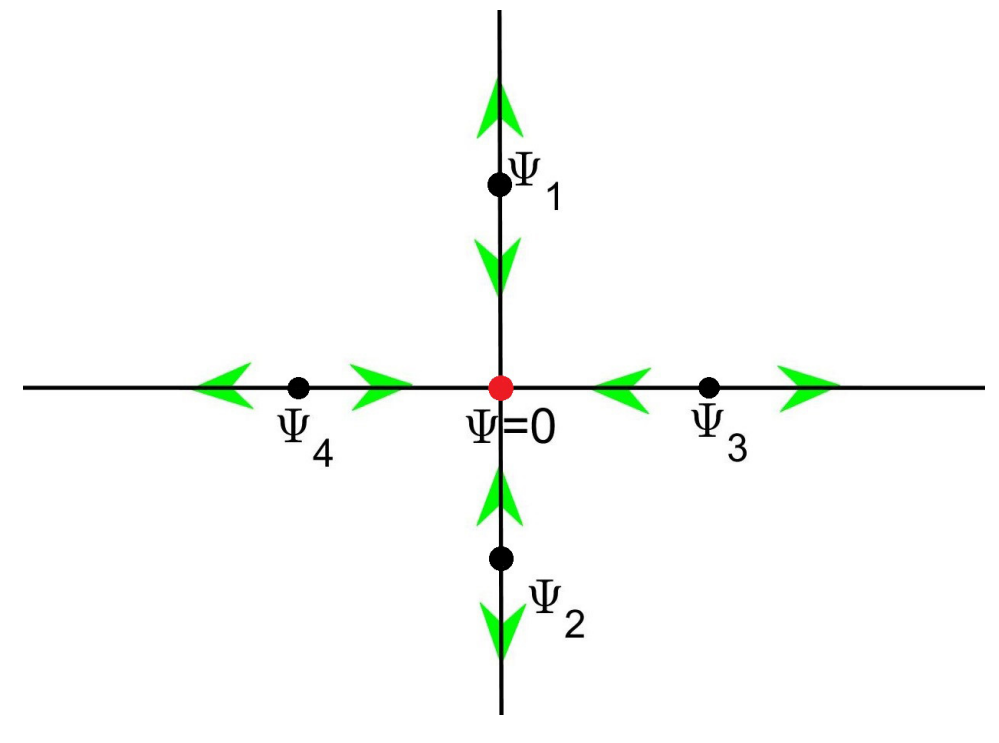}	~~~& \includegraphics[width=.3\linewidth,height=.2\linewidth]{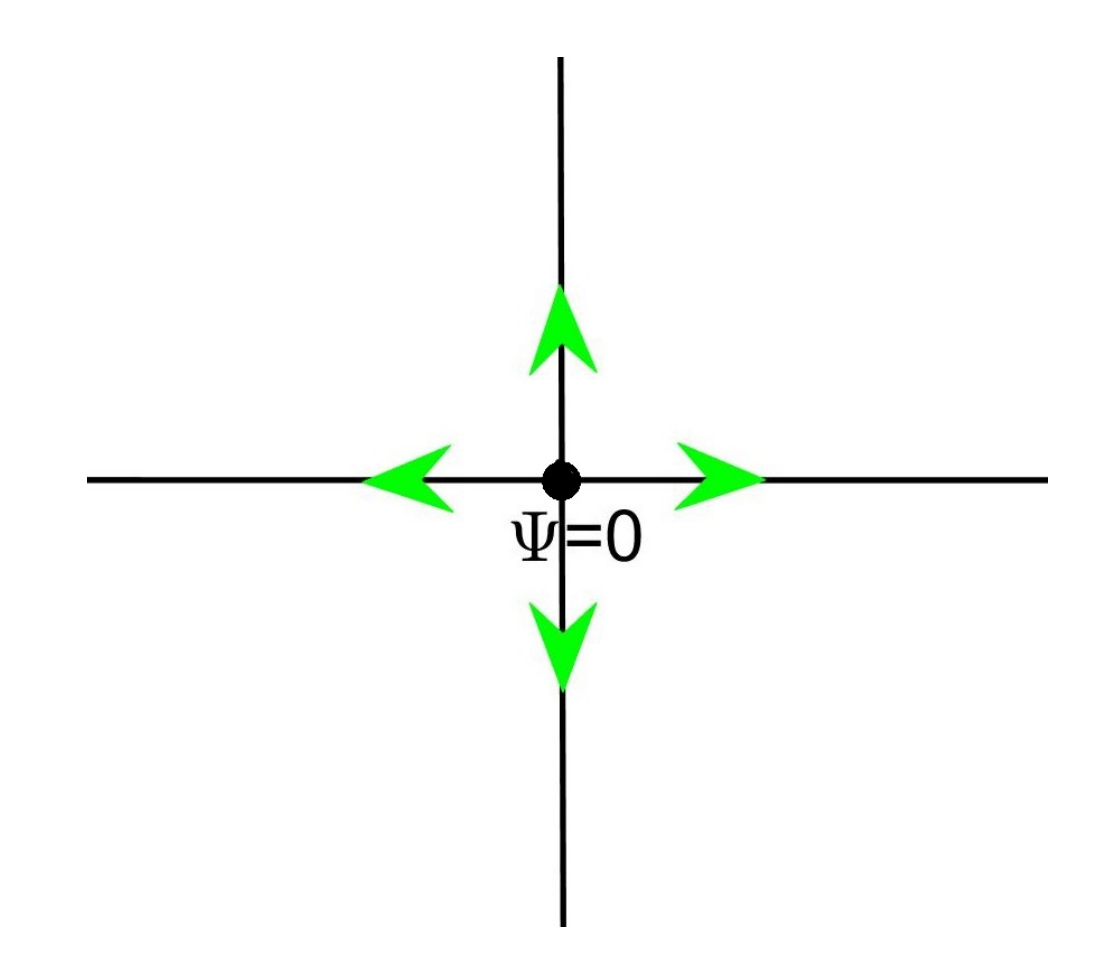}\\
				$(a)$~~~~& $(b)$
	\end{tabular}
	\caption{Topological structure of jump transition of \eqref{abstract1}
	when $\Gamma_{1}>0$ and $(\Gamma_{1}-\Gamma_{2})(\Gamma_{1}-\Gamma_{3})<0$:
 (a) $\Ra< \Ra_{c_{1}}$; (b)  $\Ra> \Ra_{c_{1}}$.}  
	  \label{diyigedin1993}
\end{figure}

\begin{figure}[H]
	\centering
	\begin{tabular}{cc}
		\includegraphics[width=.3\linewidth,height=.2\linewidth]{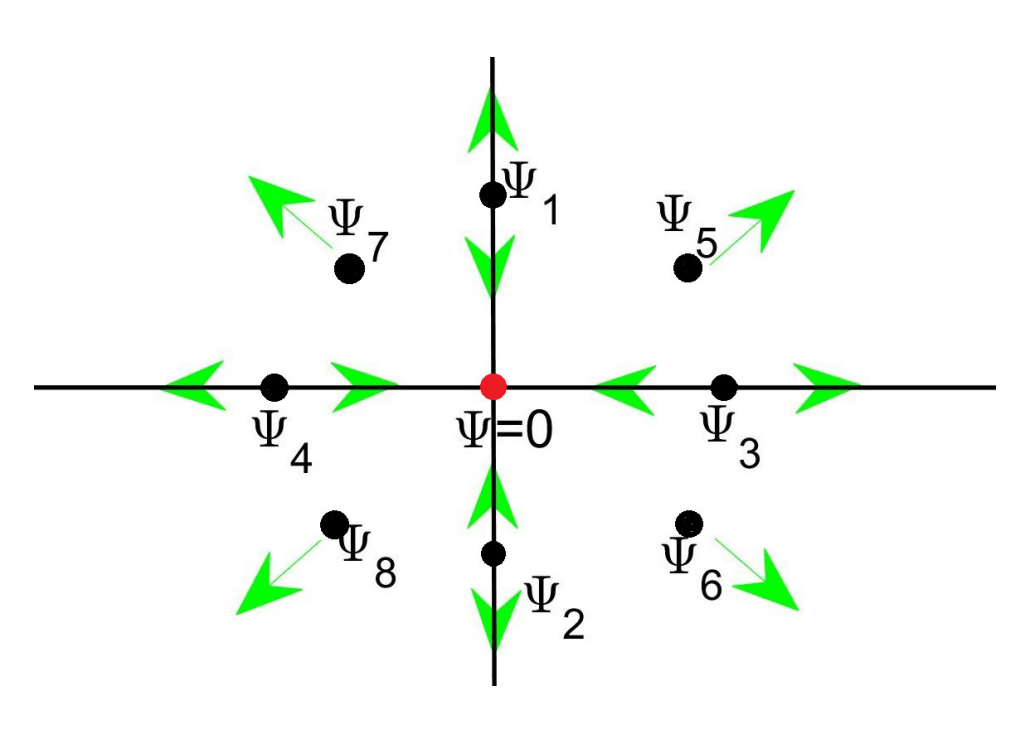}	~~~& \includegraphics[width=.3\linewidth,height=.2\linewidth]{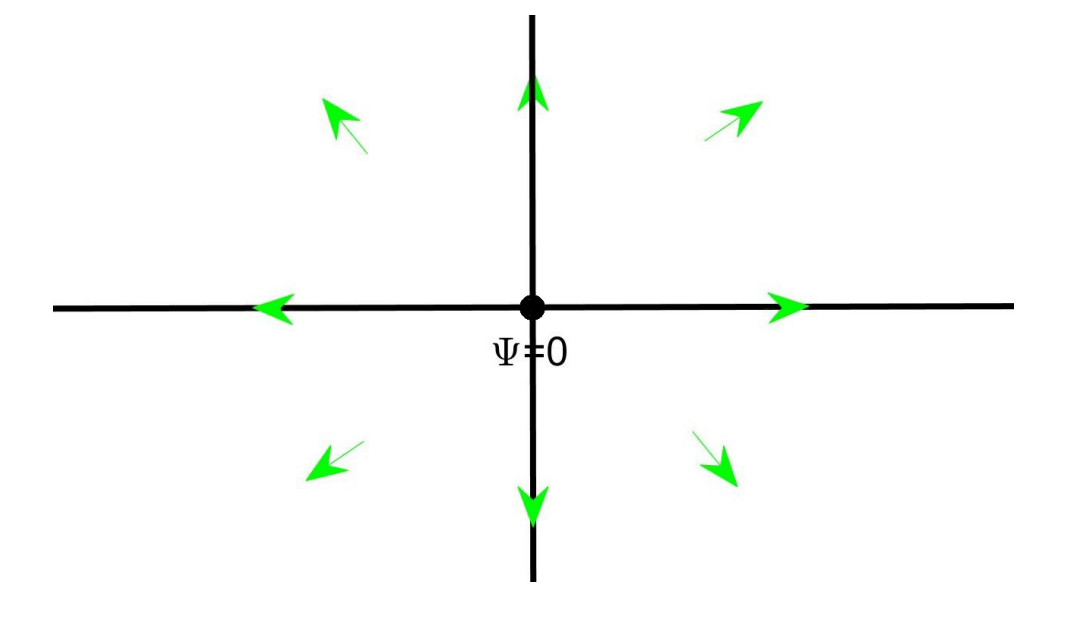}\\
				$(a)$~~~~& $(b)$
	\end{tabular}
	\caption{Topological structure of jump transition of \eqref{abstract1}
	when $\Gamma_{1}>0$, $(\Gamma_{1}-\Gamma_{2})(\Gamma_{2}\Gamma_{3}-\Gamma_{1}^{2})<0$ and $(\Gamma_{1}-\Gamma_{3})(\Gamma_{2}\Gamma_{3}-\Gamma_{1}^{2})<0$:
 (a) $\Ra< \Ra_{c_{1}}$; (b)  $\Ra> \Ra_{c_{1}}$.}  
	  \label{diyigedinglitu7}
\end{figure}

Note that $\Gamma_{s}(s=1,2,3)$ are continuously dependent on $\Ra$. However, 
the transition from the first real eigenvalue with multiplicity two at $\Ra=Ra_{{c_{1}}}$
are determined by the values of $\Gamma_{s}(Ra_{{c_{1}}})~(s=1,2,3)$. Thus, in the following, we alway use $\Gamma_{s}~(s=1,2,3)$ to represent the $\Gamma_{s}(Ra_{{c_{1}}})$. Based on the ODEs \eqref{reducedODE3}, we have the following theorem:
\begin{theorem}
    \label{double-real-theorem} For the system \eqref{abstract1}, we have the following conclusions:
    \begin{enumerate}
        \item [\rm{(1)}]If $\Gamma_{1}<0$, $\Gamma_{2}\Gamma_{3}>\Gamma_{1}^{2}$, $\Gamma_{1}>\Gamma_{2}$ and $\Gamma_{1}>\Gamma_{3}$, it has a continuous transition from
        $(\Psi,\Ra)=(0,\Ra_{c_1})$, and bifurcates on $Ra>Ra_{c_{1}}$ to an attractor $\mathcal{A}$ which exactly contains eight non-degenerate equilibrium points $\Psi_m~(m=1,\cdots,8)$ and is homeomorphic to the one-dimensional sphere $S^1$, as shown in
\autoref{diyigedinglitu6dd}. Among them, $\Psi_m~(m=1,\cdots,4)$
are stable while $\Psi_m~(m=5,\cdots,8)$ are unstable.
        \item[\rm{(2)}] If  $\Gamma_{1}<0$, $\Gamma_{1}<\Gamma_{2}$ and $\Gamma_{1}>\Gamma_{3}$, it has a continuous transition from $(\Psi,\Ra)=(0,\Ra_{c_1})$, and bifurcates on $Ra>Ra_{c_{1}}$ to an attractor $\mathcal{A} $  which exactly contains four non-degenerate equilibrium points $\Psi_m~(m=1,\cdots,4)$ and is homeomorphic to the one-dimensional sphere $S^1$, as shown in \autoref{diyigedinglitu156}. 
        Among them, $\Psi_m~(m=3,4)$
are stable while $\Psi_m~(m=1,2)$ are unstable.
        \item[\rm{(3)}] If $\Gamma_{1}>0$, $(\Gamma_{1}-\Gamma_{2})(\Gamma_{2}\Gamma_{3}-\Gamma_{1}^{2})>0$ and $(\Gamma_{1}-\Gamma_{3})(\Gamma_{2}\Gamma_{3}-\Gamma_{1}^{2})>0$, it has a jump transition from $(\Psi,\Ra)=(0,\Ra_{c_1})$, and bifurcates on both sides of $Ra=Ra_{c_{1}}$ to
eight unstable non-degenerate points $\Psi_m~(m=1,\cdots,8)$, as shown in \autoref{diyigedinglitu1993}.
        \item[\rm{(4)}] If $\Gamma_{1}>0$ and $(\Gamma_{1}-\Gamma_{2})(\Gamma_{1}-\Gamma_{3})<0$, it has a jump transition from $(\Psi,\Ra)=(0,\Ra_{c_1})$, and bifurcates on $Ra<Ra_{c_{1}}$ to four unstable non-degenerate equilibrium points $\Psi_m~(m=1,\cdots,4)$, as shown in \autoref{diyigedin1993} .
        \item[\rm{(5)}] If $\Gamma_{1}>0$, $(\Gamma_{1}-\Gamma_{2})(\Gamma_{2}\Gamma_{3}-\Gamma_{1}^{2})<0$ and $(\Gamma_{1}-\Gamma_{3})(\Gamma_{2}\Gamma_{3}-\Gamma_{1}^{2})<0$,  it has a jump transition from $(\Psi,\Ra)=(0,\Ra_{c_1})$, and bifurcates on $Ra<Ra_{c_{1}}$ to eight unstable non-degenerate equilibrium points $\Psi_m~(m=1,\cdots,4)$, as shown in \autoref{diyigedinglitu7}.
    \end{enumerate}
\end{theorem}

\begin{figure}[H]
	\centering
	\begin{tabular}{cc}
		\includegraphics[width=.3\linewidth,height=.2\linewidth]{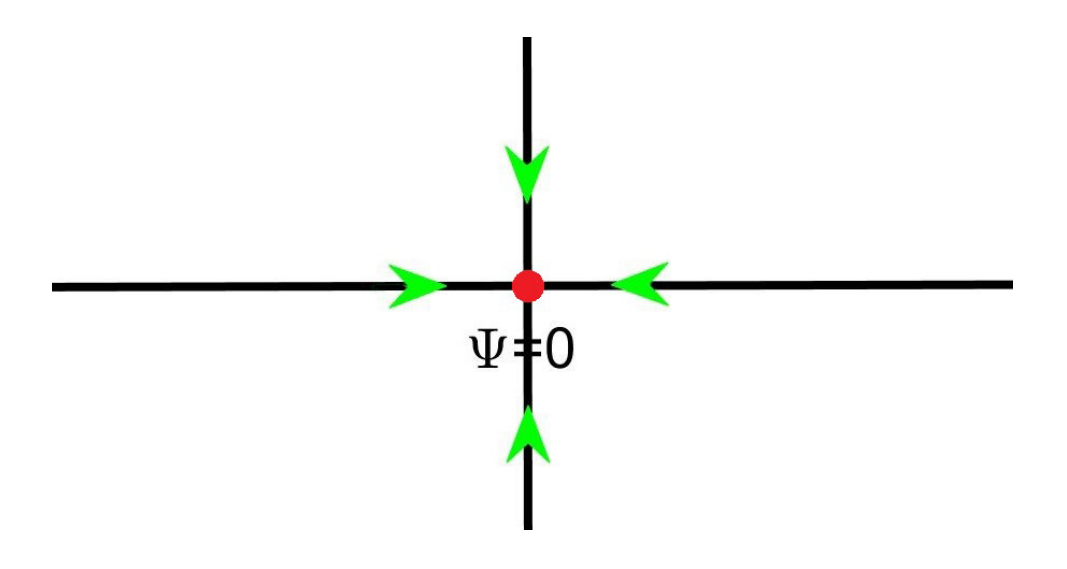}	~~~& \includegraphics[width=.3\linewidth,height=.2\linewidth]{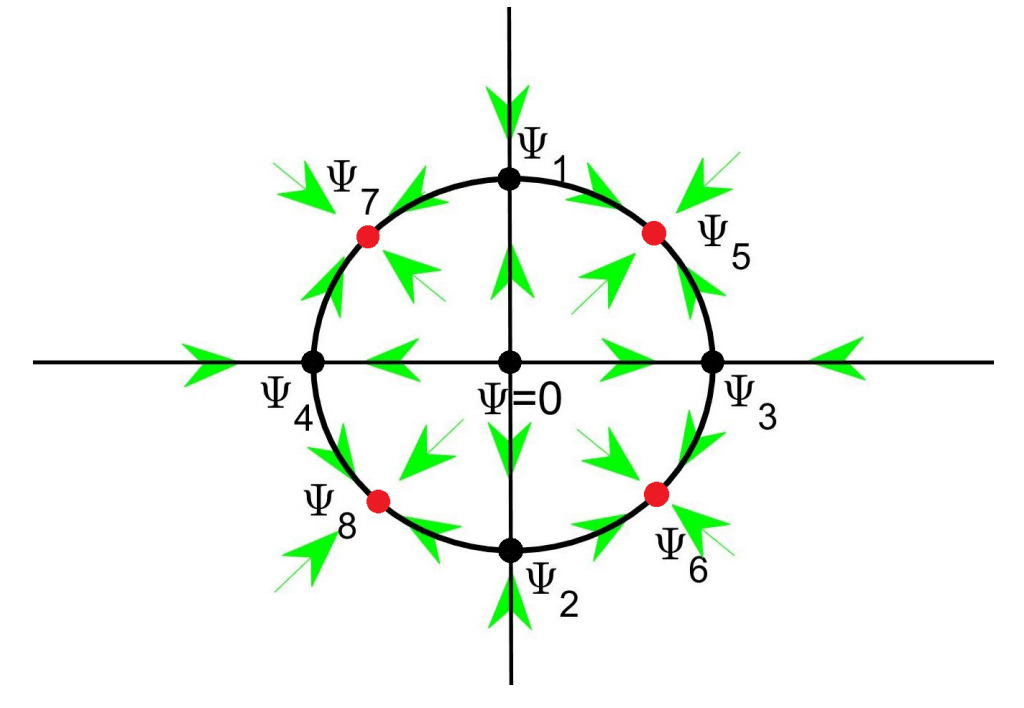}\\
				$(a)$~~~~& $(b)$
	\end{tabular}
	\caption{Topological structure of continuous transition of \eqref{abstract1}
	when $\Gamma_{1}<0$, $\Gamma_{2}\Gamma_{3}<\Gamma_{1}^{2}$, $\Gamma_{1}<\Gamma_{2}$ and $\Gamma_{1}<\Gamma_{3}$:
 (a) $\Ra< \Ra_{c_{1}}$; (b)  $\Ra> \Ra_{c_{1}}$.}  
	  \label{diyigedinglitu6}
\end{figure}

\begin{figure}[H]
	\centering
	\begin{tabular}{cc}
		\includegraphics[width=.3\linewidth,height=.2\linewidth]{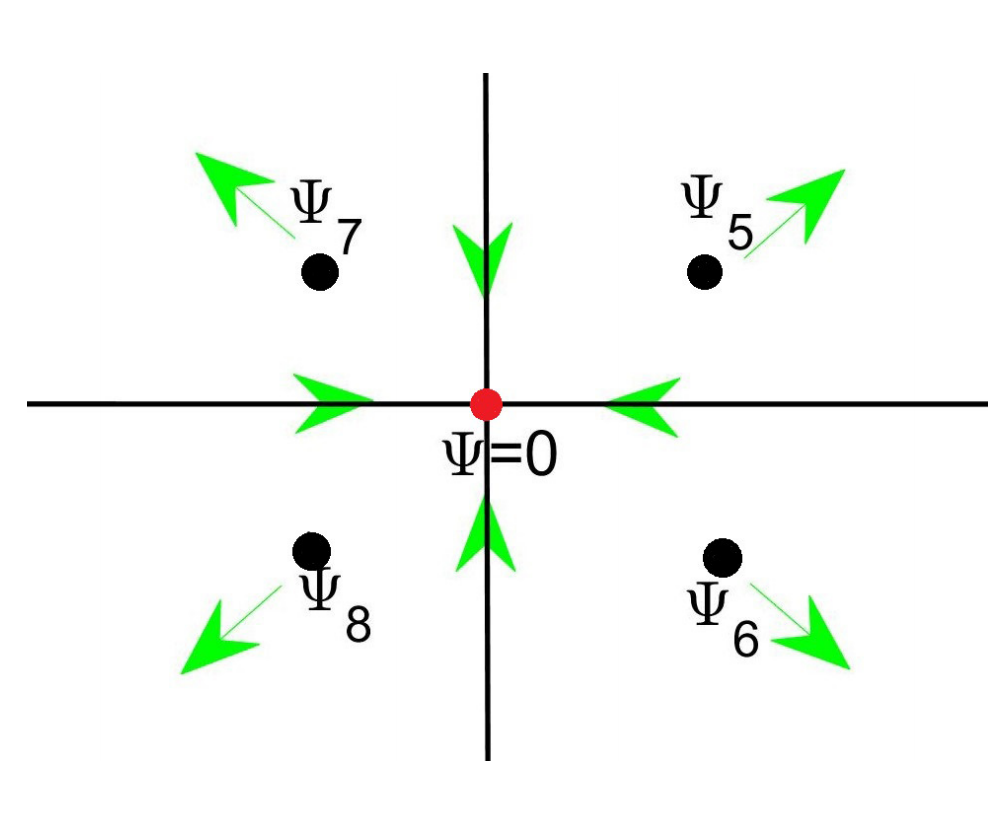}	~~~& \includegraphics[width=.3\linewidth,height=.2\linewidth]{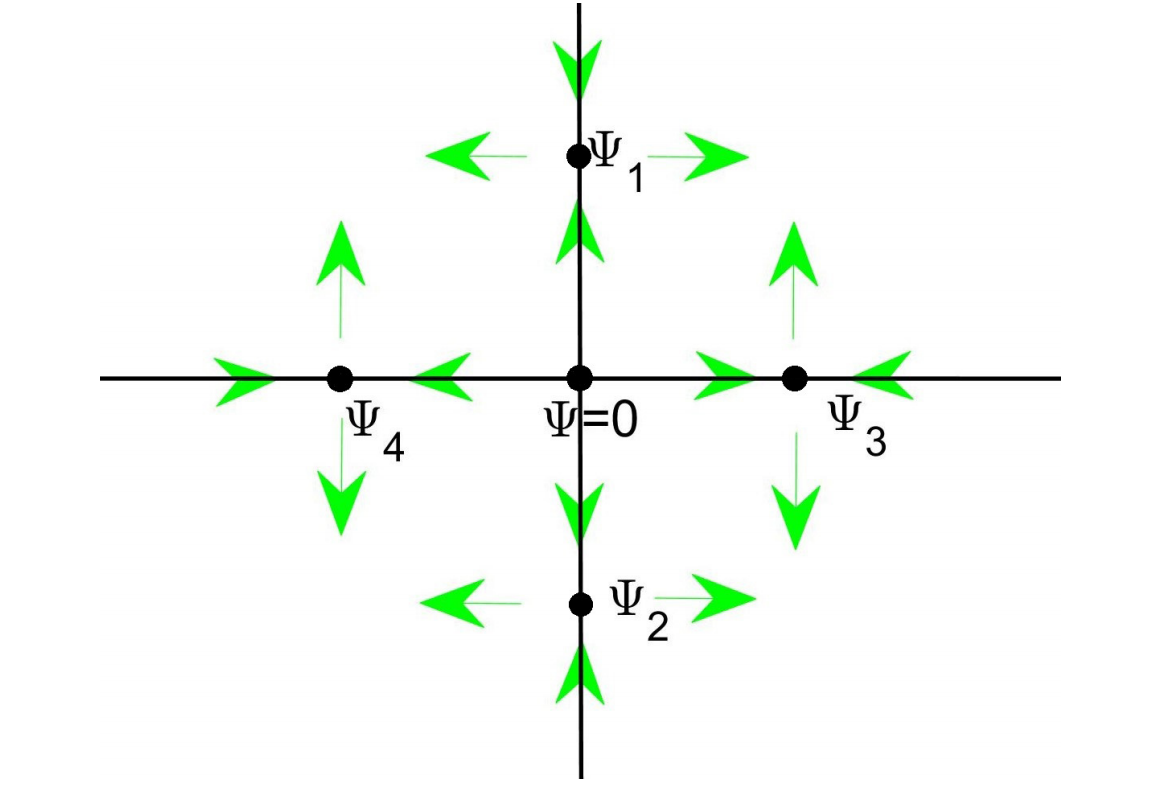}\\
				$(a)$~~~~& $(b)$
	\end{tabular}
	\caption{Topological structure of jump transition of \eqref{abstract1}
	when $\Gamma_{1}<0$, $\Gamma_{2}\Gamma_{3}>\Gamma_{1}^{2}$, $\Gamma_{1}<\Gamma_{2}$ and $\Gamma_{1}<\Gamma_{3}$:
 (a) $\Ra< \Ra_{c_{1}}$; (b)  $\Ra> \Ra_{c_{1}}$.}  
	  \label{diyigedinglitu1323}
\end{figure}

\begin{figure}[H]
	\centering
	\begin{tabular}{cc}
		\includegraphics[width=.3\linewidth,height=.2\linewidth]{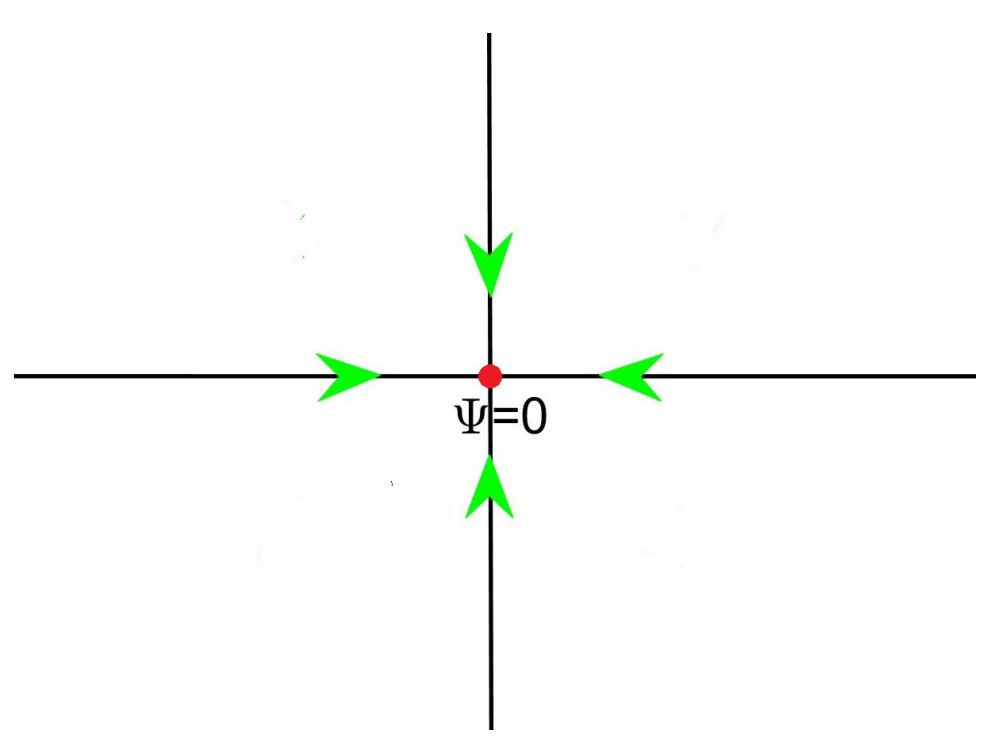}	~~~& \includegraphics[width=.3\linewidth,height=.2\linewidth]{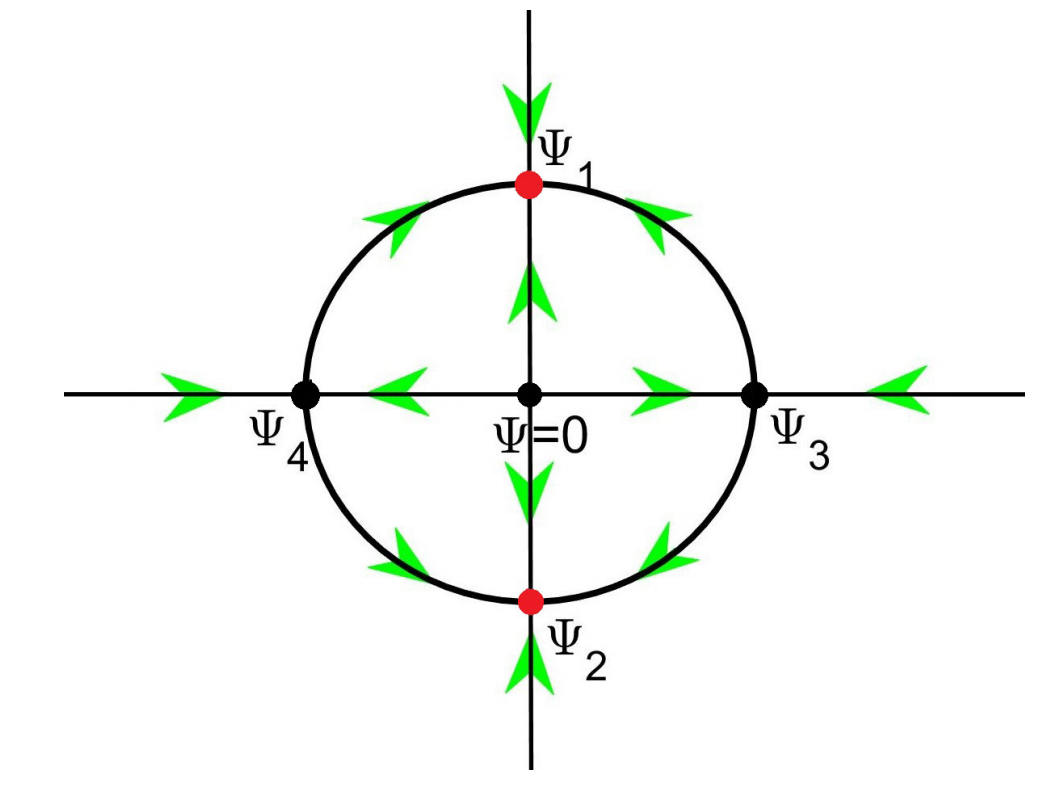}\\
				$(a)$~~~~& $(b)$
	\end{tabular}
	\caption{Topological structure of jump transition of \eqref{abstract1}
	when $\Gamma_{1}<0$, $\Gamma_{1}>\Gamma_{2}$ and $\Gamma_{1}<\Gamma_{3}$:
 (a) $\Ra< \Ra_{c_{1}}$; (b)  $\Ra> \Ra_{c_{1}}$.}  
	  \label{diyigedinglitu1377}
\end{figure}

\begin{remark} Actually, in theory, there are three other situations as follows:
    \begin{enumerate}
    \item [\rm{(1)}]If $\Gamma_{1}<0$, $\Gamma_{2}\Gamma_{3}<\Gamma_{1}^{2}$, $\Gamma_{1}<\Gamma_{2}$ and $\Gamma_{1}<\Gamma_{3}$, it has a continuous transition from
        $(\Psi,\Ra)=(0,\Ra_{c_1})$, and bifurcates on $Ra>Ra_{c_{1}}$ to an attractor $\mathcal{A}$ which exactly contains eight non-degenerate equilibrium points $\Psi_m~(m=1,\cdots,8)$ and is homeomorphic to the one-dimensional sphere $S^1$, as shown in
 \autoref{diyigedinglitu6}. Among them, $\Psi_m~(m=1,\cdots,4)$
are unstable while $\Psi_m~(m=5,\cdots,8)$ are stable.
    \item [\rm{(2)}]If $\Gamma_{1}<0$, $\Gamma_{2}\Gamma_{3}>\Gamma_{1}^{2}$, $\Gamma_{1}<\Gamma_{2}$ and $\Gamma_{1}<\Gamma_{3}$, it has a jump transition from $(\Psi,\Ra)=(0,\Ra_{c_1})$, and bifurcates on both sides of $Ra=Ra_{c_{1}}$ to
eight unstable non-degenerate points $\Psi_m~(m=1,\cdots,8)$, as shown in
    \autoref{diyigedinglitu1323}.
    \item [\rm{(3)}]If $\Gamma_{1}<0$, $\Gamma_{1}>\Gamma_{2}$ and $\Gamma_{1}<\Gamma_{3}$, it has a continuous transition from $(\Psi,\Ra)=(0,\Ra_{c_1})$, and bifurcates on $Ra>Ra_{c_{1}}$ to an attractor $\mathcal{A} $  which exactly contains four non-degenerate equilibrium points $\Psi_m~(m=1,\cdots,4)$ and is homeomorphic to the one-dimensional sphere $S^1$, as shown in \autoref{diyigedinglitu1377}. 
        Among them, $\Psi_m~(m=1,2)$
are stable while $\Psi_m~(m=3,4)$ are unstable.
    \end{enumerate}
However, we do not find the corresponding examples in numerical stimulation.
\end{remark}

\begin{proof} 
    We only prove the first conclusion, the others can be proven in a similar way.
    
 For $\beta>0$, i.e. $\Ra>\Ra_{c_{1}}$, we can infer from  $\Gamma_{1}<0$, $\Gamma_{2}\Gamma_{3}>\Gamma_{1}^{2}$, $\Gamma_{1}>\Gamma_{2}$ and $\Gamma_{1}>\Gamma_{3}$ that $\mathbf{Y}_{n}~(n=1,2,...8)$ all exist.  We then discuss the stability of each $\mathbf{Y}_{n}$. Linearizing the ODEs \eqref{reducechenweiode} at $\mathbf{Y}=\mathbf{Y}_{n}$, the corresponding matrix reads
    \begin{align}\label{xianxingjuzhen}
        Ma=\begin{pmatrix}
        \beta+3\Gamma_{1} y^2+\Gamma_{2} z^2   &2\Gamma_{2} y z\\
        2\Gamma_{3} y z                         &\beta+\Gamma_{3}y^2+3\Gamma_{1}z^2
        \end{pmatrix},
    \end{align}
where $(y,z)=\mathbf{Y}_{n}~(n=1,2,...8)$. Then through judging the sign of eigenvalues of $Ma$ at each $\mathbf{Y_{n}},$ we know $\mathbf{Y}_{n}(n=1,2,3,4)$ are stable while $\mathbf{Y}_{n}(n=5,6,7,8)$ are unstable.

For $\beta<0$, i.e. $\Ra<\Ra_{c_{1}}$, the unique zero equilibrium point is asymptotic stable.

\end{proof}
 
\section{Numerical investigations}\label{section6}

In the preceding sections, the stabilities and transition of the system \eqref{abstract1} have been theoretically studied in various scenarios. In this section, some numerical examples are offered to show the specific transition types.  Numerical investigation in the following is carried out within the parameter ranges \(750 \leq  \Ta \leq  3000\), \(0 < \Q \leq  1000\) and \( 0 < \Pr \leq 1  \). 

We numerically analyze the impact of the Chandrasekahr number \( \Q \)  on the values of the two critical parameters \( \Ra_{c_1} \) and \( \Ra_{c_2} \).  The influence of  \( \Q \) on the two critical parameters is very subtle, see \autoref{Ta2700Pr05L11L212}. 
If the critical index is \( (j , 0 , l) \), \( \Q \) has no effect on the critical parameters.  
 If the critical index is not \( (j , 0 , l) \), 
 the effect of \( \Q \) on critical parameters can always be canceled by decreasing of \( \abs{k} \) where $k$ is the second component of the critical index.  The reason lies in the two expressions in \eqref{critical1} from which we see that the number \( \Q \) is always binding with the second component \( k \).
 
\begin{figure}[H]
    \begin{minipage}[t]{0.45\linewidth}
        \centering
        \includegraphics[width=1.8in]{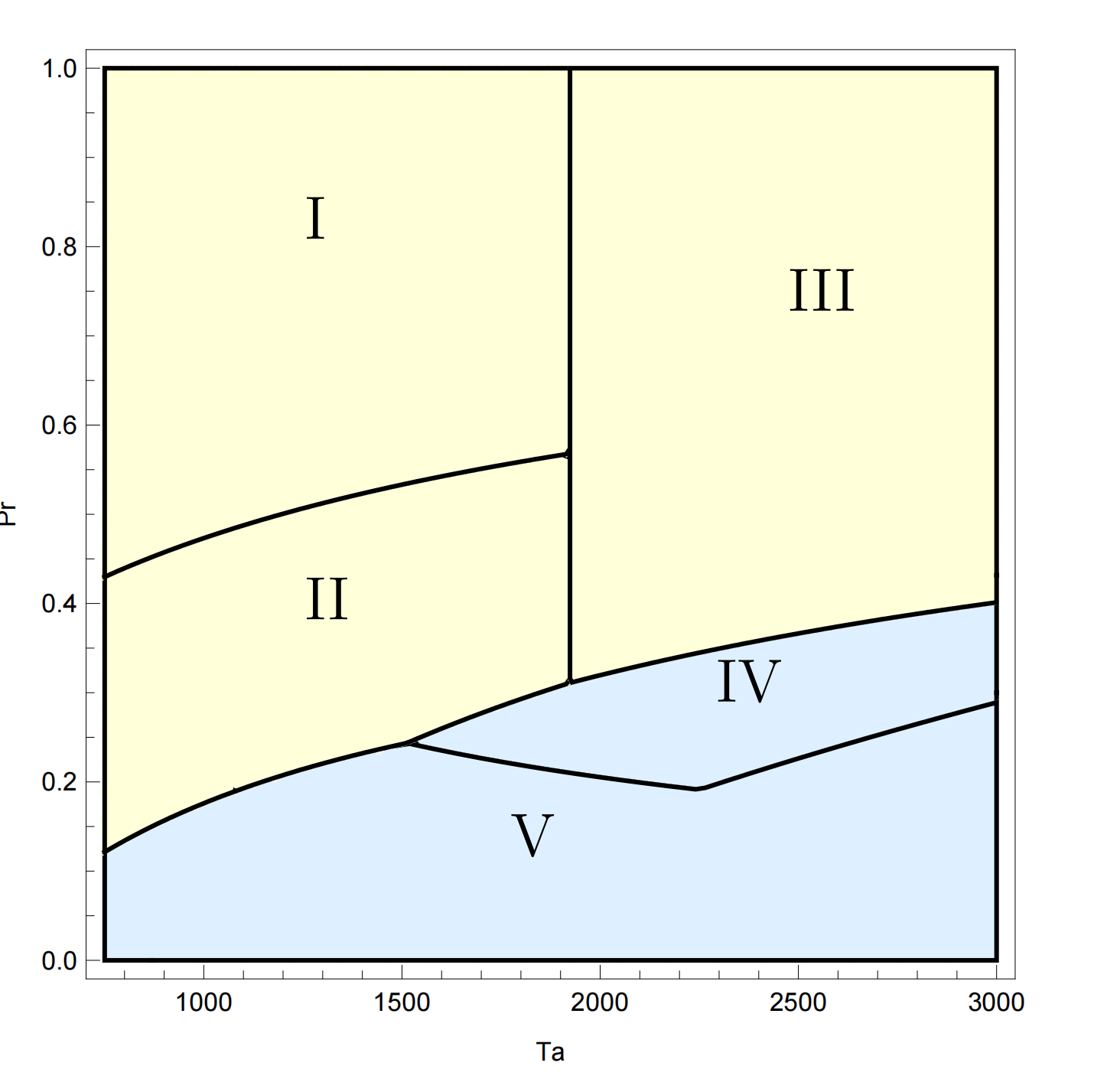}
          
    \end{minipage}
    \hfill
    \begin{minipage}[t]{0.6\linewidth}
        \centering
        \includegraphics[width=1.8in]{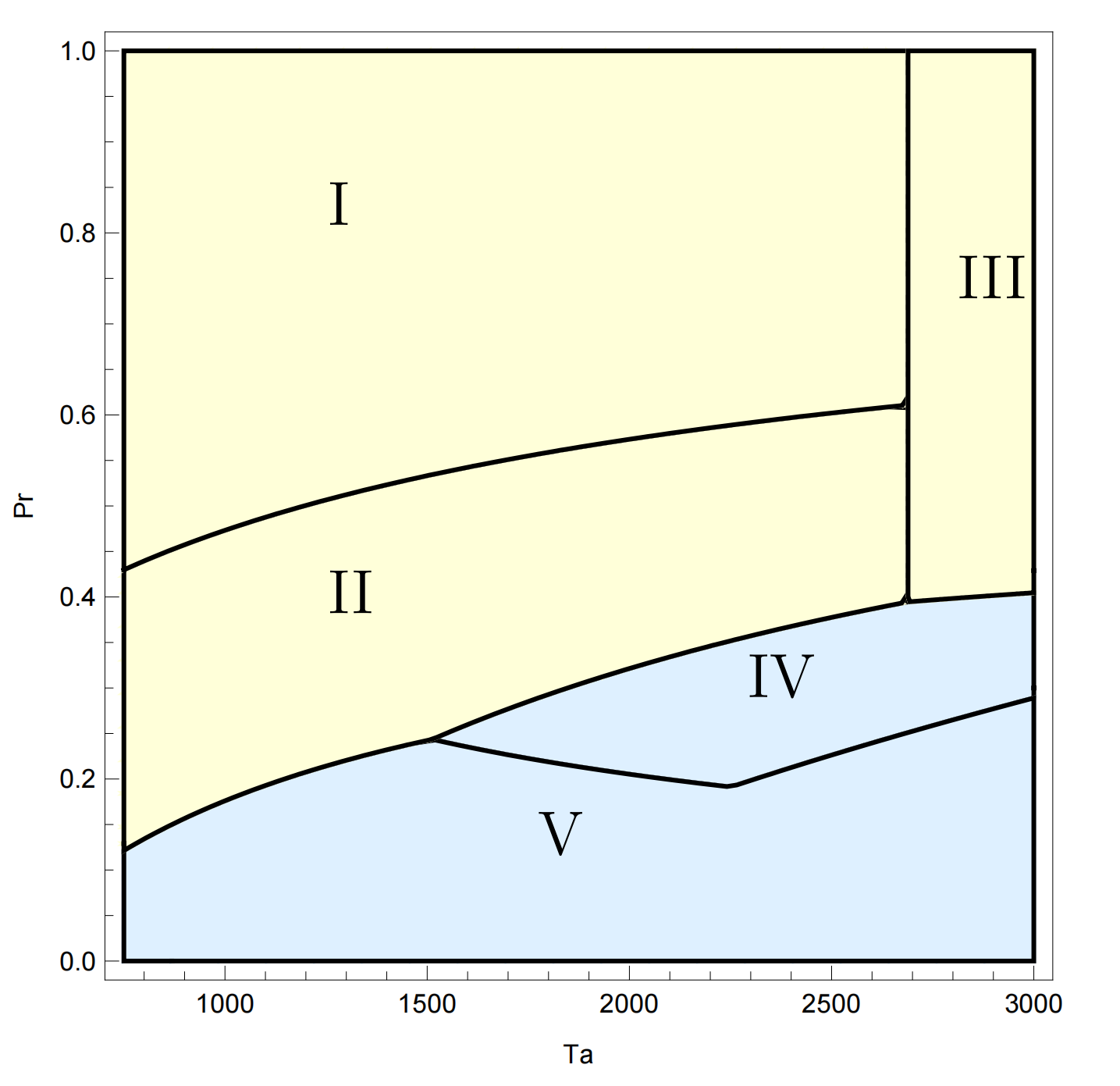}
    \end{minipage}
    \caption{The \( \Ta -\Pr \) plane showing the relation between critical numbers \( \Ra_{c_1} \) and \( \Ra_{c_2} \), for \( L_{1} = 1 \), \( L_2 = 1 . 2 \) and \( \Q = 100 \) (left) \( \Q = 500 \) (right). 
    \MakeUppercase{\romannumeral1}, \MakeUppercase{\romannumeral2} and \MakeUppercase{\romannumeral3} are the domains where \( \Ra_{c_1} < \Ra_{c_2} \) while \MakeUppercase{\romannumeral4} and \MakeUppercase{\romannumeral5} are the domains where \( \Ra_{c_1} > \Ra_{c_2} \). 
    }
    \label{Rac1Rac1Relation}
\end{figure}

We also perform some numerical analysis to examine the influences of the Taylor number \( \Ta \) and the Prandtl number \( \Pr \) on the values of critical parameters  \( \Ra_{c_1} \) and \( \Ra_{c_2} \), the type and multiplicity of first eigenvalues and the signs of transition numbers. In \autoref{Rac1Rac1Relation}, \MakeUppercase{\romannumeral1}, \MakeUppercase{\romannumeral2} and \MakeUppercase{\romannumeral3} are the parameter domains for $(\Ta,\Pr)$ in which \( \Ra_{c_1} <\Ra_{c_2}\)  while  \( \Ra_{c_1} >\Ra_{c_2}\) for $(\Ta,\Pr)$ in domains \MakeUppercase{\romannumeral4} and \MakeUppercase{\romannumeral5}. Besides, the first component $j$ of critical indexes of all cases are always non-zero, see \autoref{Ta2700Q100Pr06L1052L2052}.  The symmetry of expressions in \eqref{critical1} with respect to the first two components of index allow us to use the second component k of critical index to judge the multiplicity of the first eigenvalue. That is, for the critical index $(j,k,l)\in X$, if $k\neq 0$, then the multiplicity of first eigenvalue is two because $(j,-k,l)\in X$, otherwise it is one.
In fact, our numerical investigations show that there is no eigenvalue whose multiplicity is greater than two. Hence, we roughly say that for $(\Ta,\Pr)$ in \MakeUppercase{\romannumeral1} and \MakeUppercase{\romannumeral2}, the multiplicity of first real eigenvalue is one, and its multiplicity is two for $(\Ta,\Pr)$ in \MakeUppercase{\romannumeral3} on the parameter plane.

For the parameter regions shown in
 \autoref{Rac1Rac1Relation}, we also estimate the numerical values of the
 transition number $\delta(\Ra_{c_1})$ and $a(\Ra_{c_2})$, finding that
  \[
 \delta(\Ra_{c_1})
 \begin{cases}
 <0,\quad( \Ta, \Pr )\in \MakeUppercase{\romannumeral1}\\
 >0, \quad (\Ta, \Pr )\in \MakeUppercase{\romannumeral2} 
 \end{cases}\quad\text{and}\quad
a(\Ra_{c_1})
 \begin{cases}
 <0,\quad( \Ta, \Pr )\in \MakeUppercase{\romannumeral4}\\
 >0, \quad (\Ta, \Pr )\in \MakeUppercase{\romannumeral5} 
 \end{cases}.
 \]
 Hence, both scenarios of \autoref{theroemoftransition} and  \autoref{complex-transition} are realized. Here, we give two specific examples to illustrate the stable state after continuous transitions.
 
         \begin{figure}[h]
        \begin{minipage}[t]{0.5\linewidth}
            \centering
            {\includegraphics[width=.7\linewidth,height=.45\linewidth]{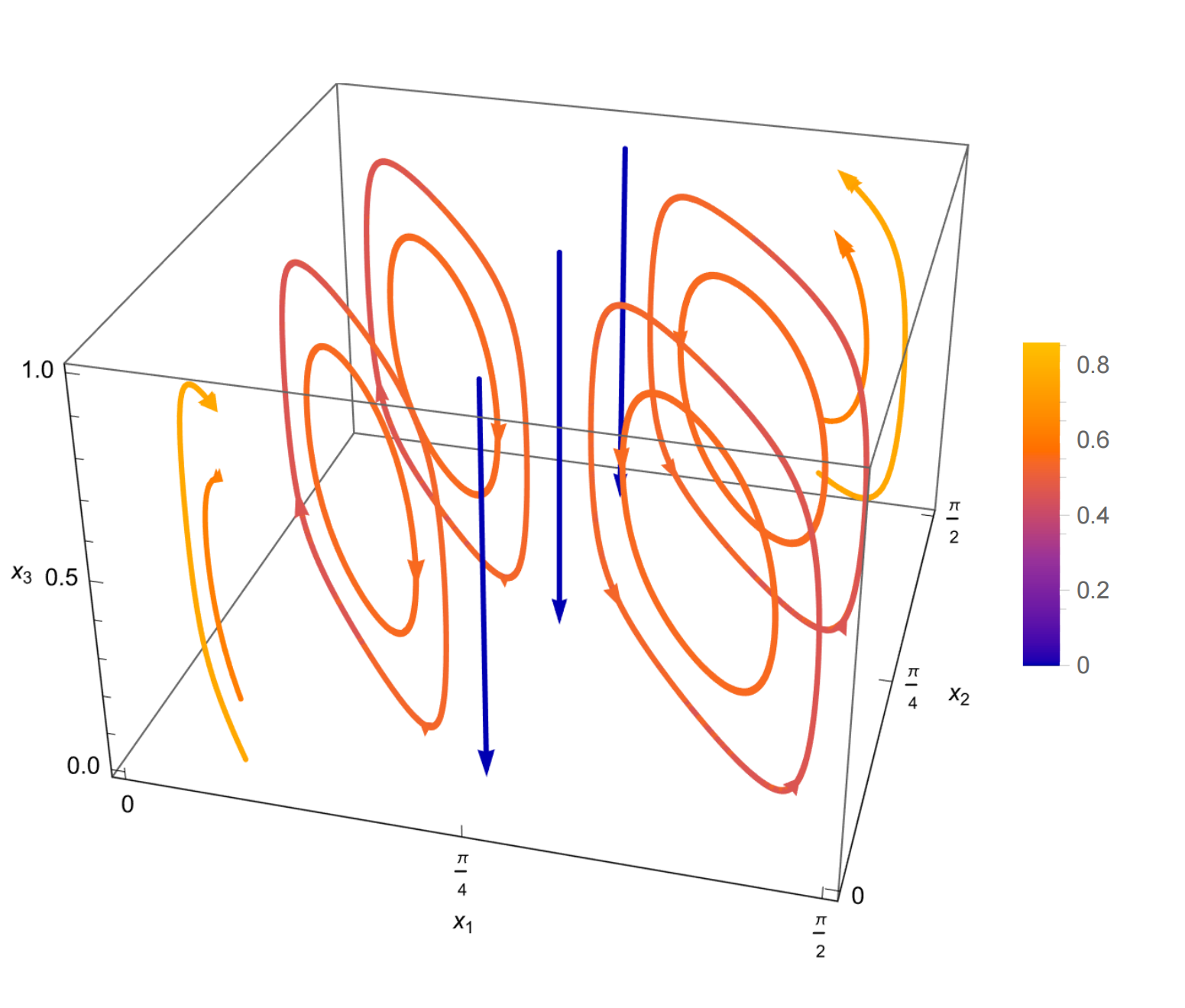}}
        \end{minipage}
        \hfill
        \begin{minipage}[t]{0.5\linewidth}
            \centering
            {\includegraphics[width=.7\linewidth,height=.45\linewidth]{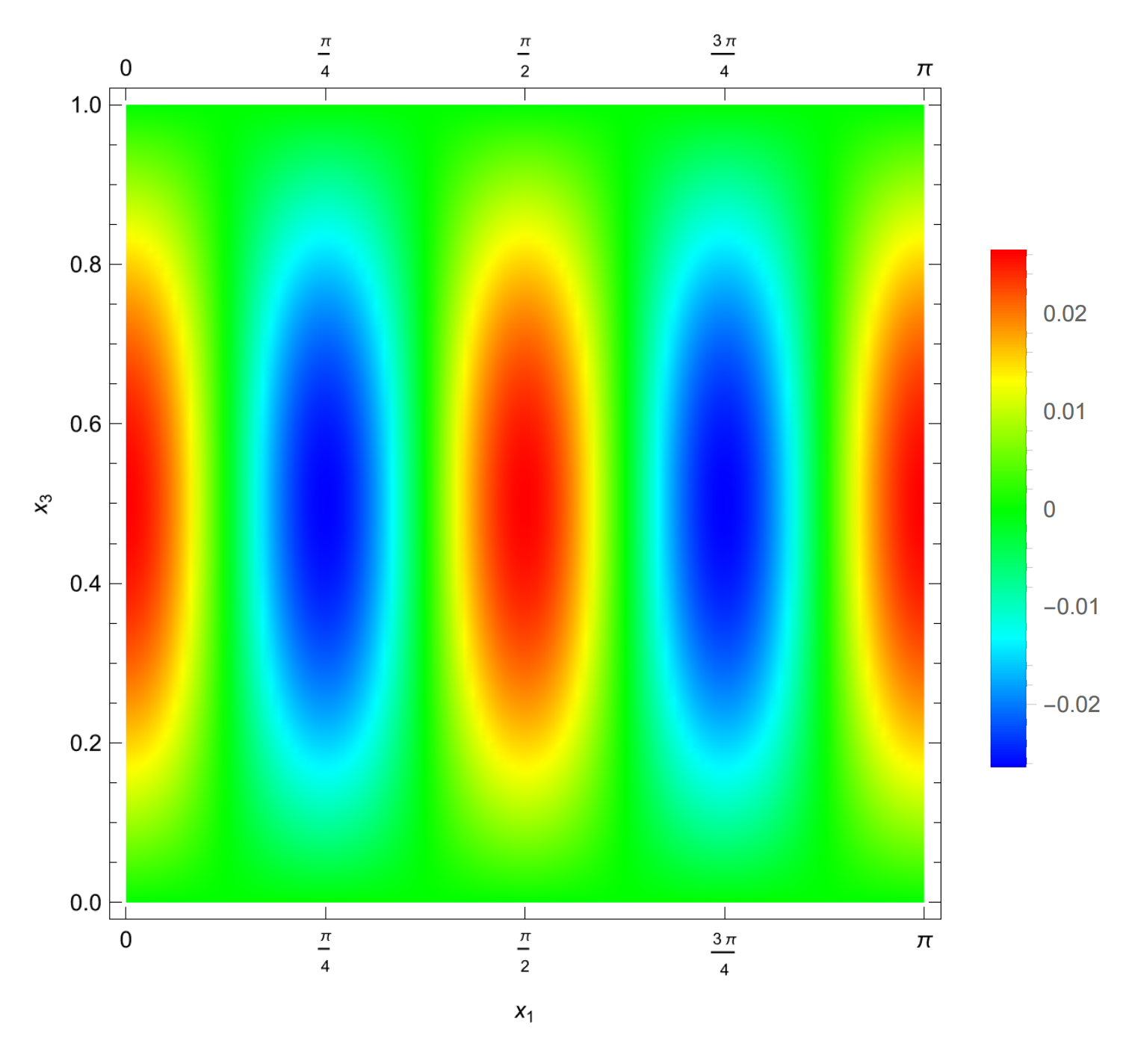}}
        \end{minipage}
        \caption{The approximate bifurcated solutions – stream plot (left) and temperature (right), where \( \Ta = 1100\), \(\Q =100\), \(\Pr = 0.7\), \(L_1 =1\), \(L_2 = 1.5 \) and \( \Ra = 1761.59>Ra_{c_{1}}\ \). }
        \label{Ta1100Q100Pr07L11L215velocityandtempreture}
    \end{figure}

Setting  \( (\Ta, \Q, \Pr, L_1,L_2) = (1100,100,0.7,1,1.5)\), we have \( \Ra_c=\Ra_{c_1 }  \approx 1760.59 \), critical index \( J_0= (j_0,k_0,1) = (4,0,1) \) and \(\delta(\Ra_{c_1}) \approx  -0.033936 \).  Namely, the system \eqref{abstract1} has a continuous transition from $(\Psi,\Ra)=(0,1760.59 )$. After the continuous transition, there are two new states $ \Psi_{1}$ and $\Psi_{2}$ given by
\[
  \Psi_{m}=(-1)^{m}\left(\frac{\beta_{J_{0}}^{1}(\Ra)}{-\delta(\Ra)}\right)^{\frac{1}{2}}\Psi_{J_{0}}^{1}
  +o(|\beta_{J_{0}}^{1}(\Ra)|^{\frac{1}{2}})~(m=1,2).
\]
Here, we plot $\Psi_{1}$ at $\Ra=1761.59$ by using its leading term, see  \autoref{Ta1100Q100Pr07L11L215velocityandtempreture}.

\begin{figure}[htb]
	\centering
	\begin{tabular}{cc}
		\includegraphics[width=.3\linewidth,height=.2\linewidth]{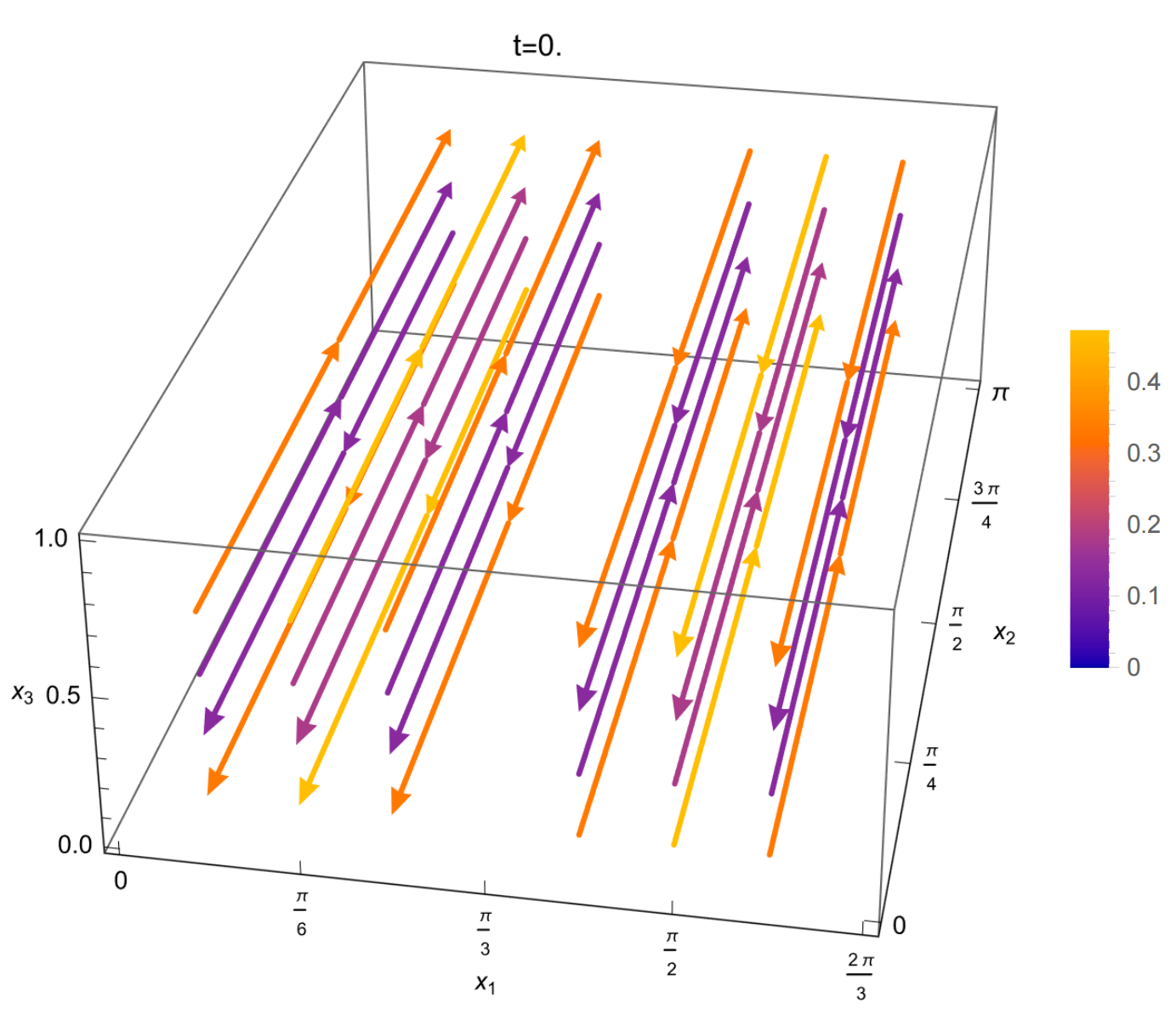}& \includegraphics[width=.3\linewidth,height=.2\linewidth]{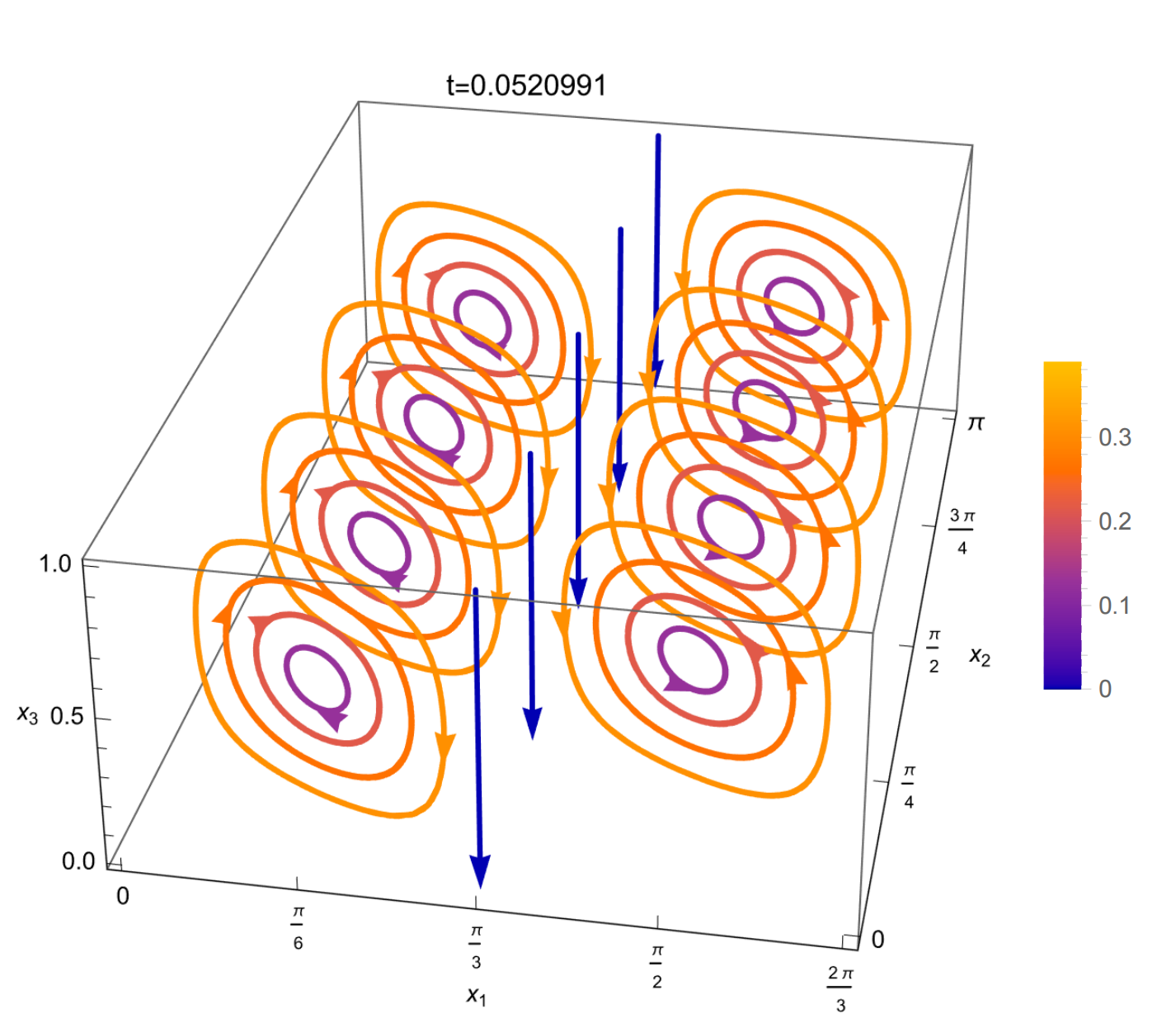}
		\\
		\includegraphics[width=.3\linewidth,height=.2\linewidth]{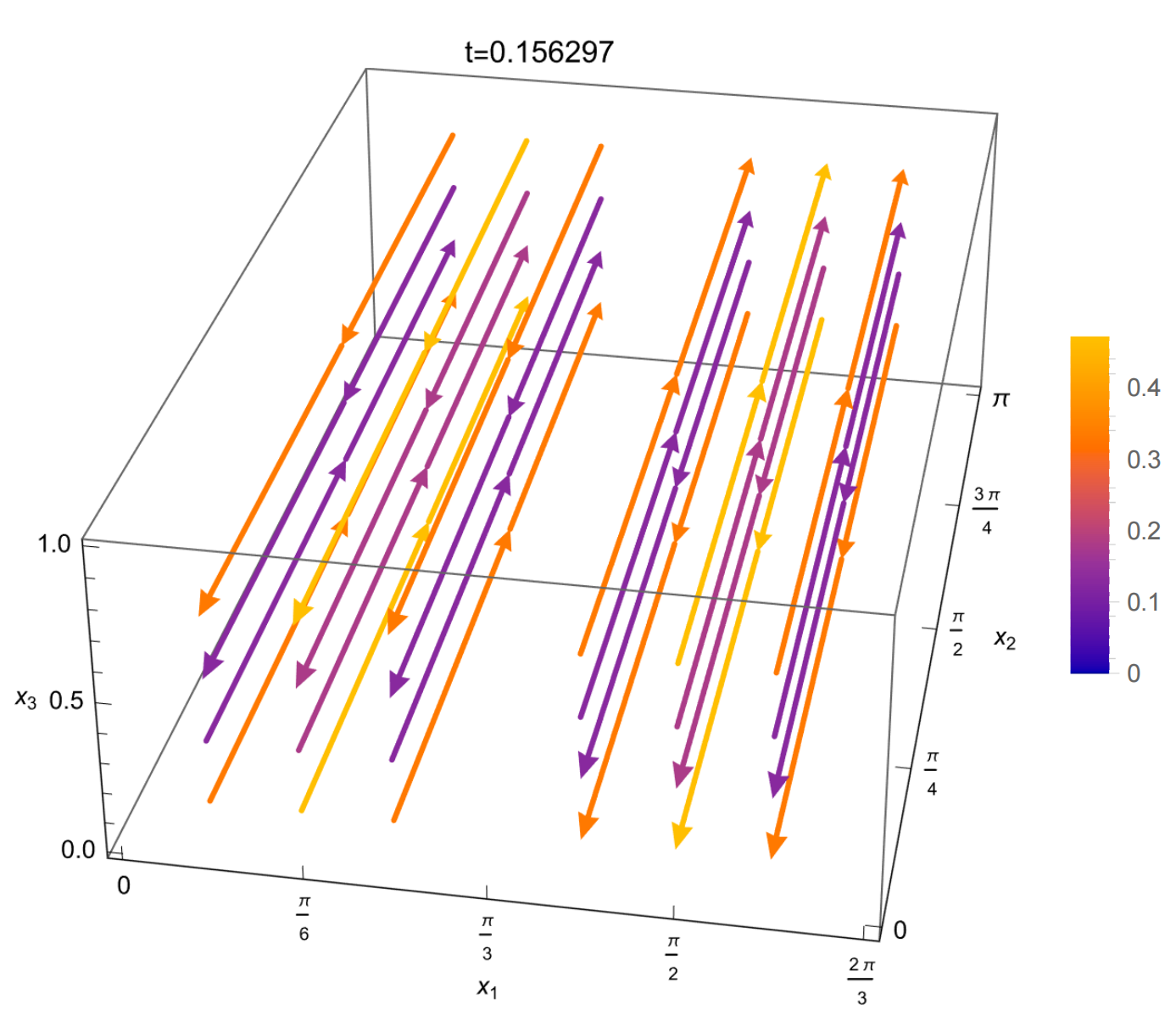}& \includegraphics[width=.3\linewidth,height=.2\linewidth]{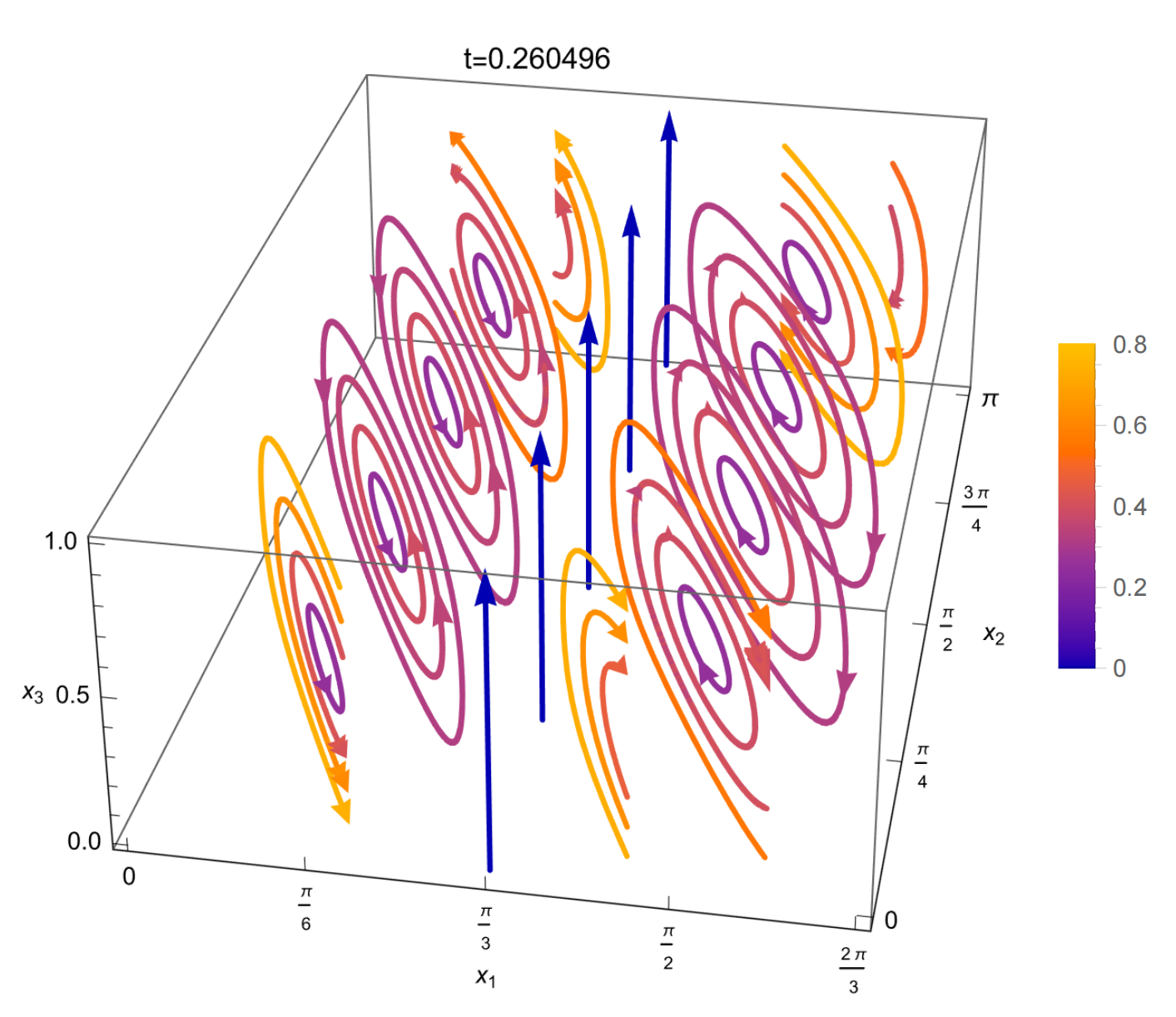}
	\end{tabular}
	\caption{Stream plot in a period for parameter configuration \( \Ta = 2700\), \(\Q = 500\), \(\Pr = 0.3\), \(L_1 =1\), \(L_2 = 1.2 \). The approximate period is \( 0.312595 \).}
    \label{Ta2700Q500Pr01L11L212stream}
\end{figure}

\begin{figure}[htb]
	\centering
	\begin{tabular}{cc}
		\includegraphics[width=.3\linewidth,height=.2\linewidth]{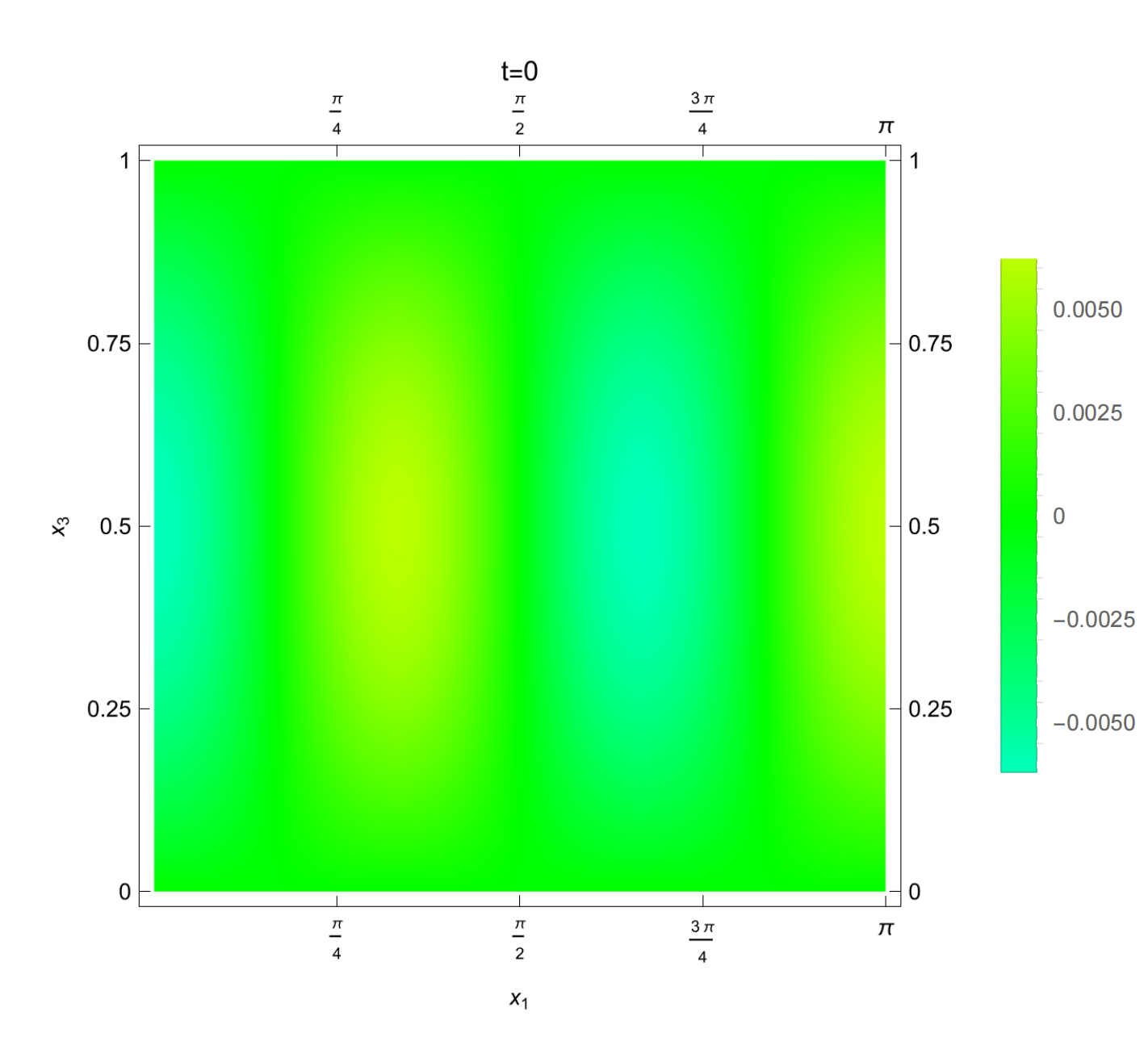}& \includegraphics[width=.3\linewidth,height=.2\linewidth]{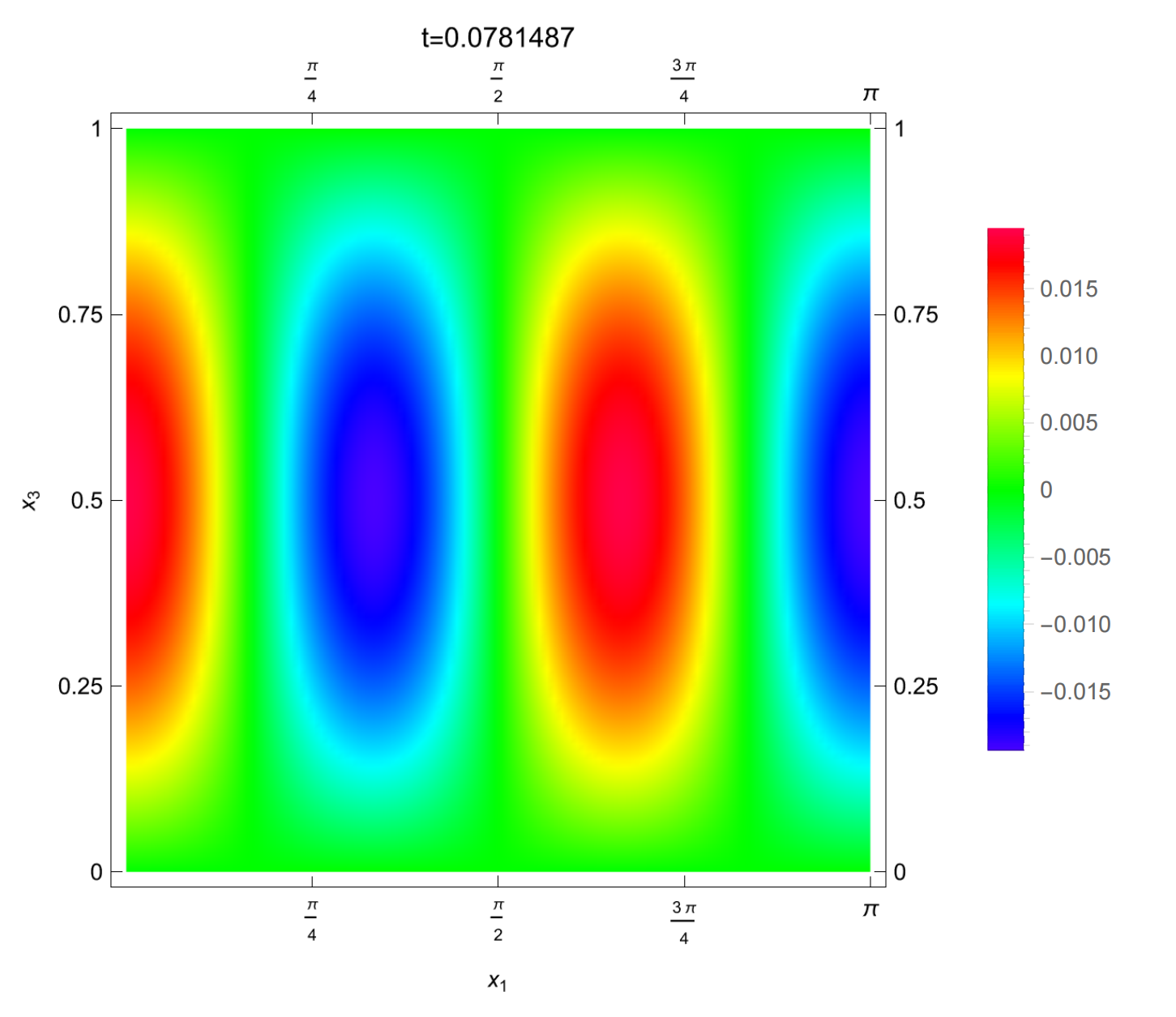}
		\\
		\includegraphics[width=.3\linewidth,height=.2\linewidth]{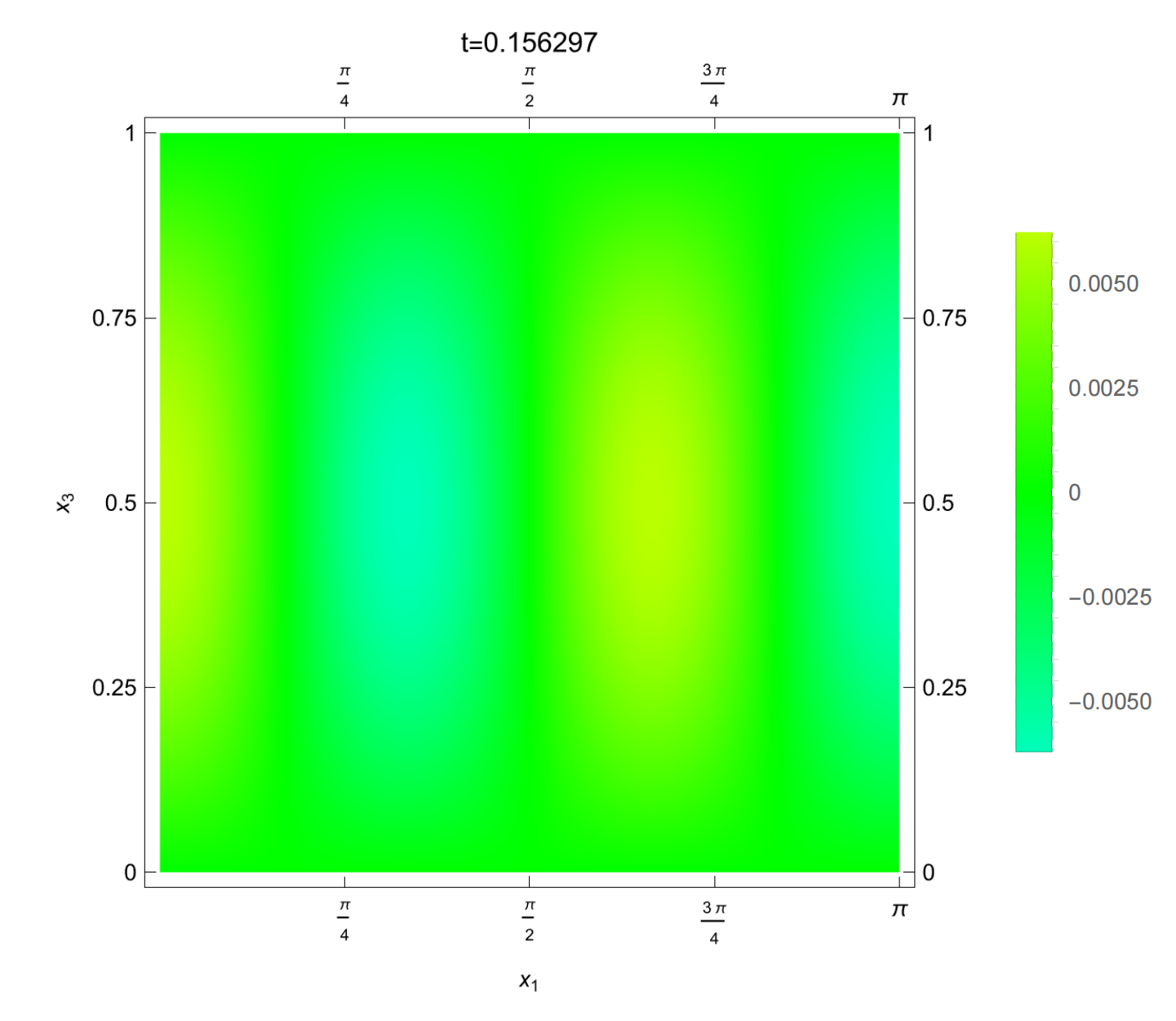}& \includegraphics[width=.3\linewidth,height=.2\linewidth]{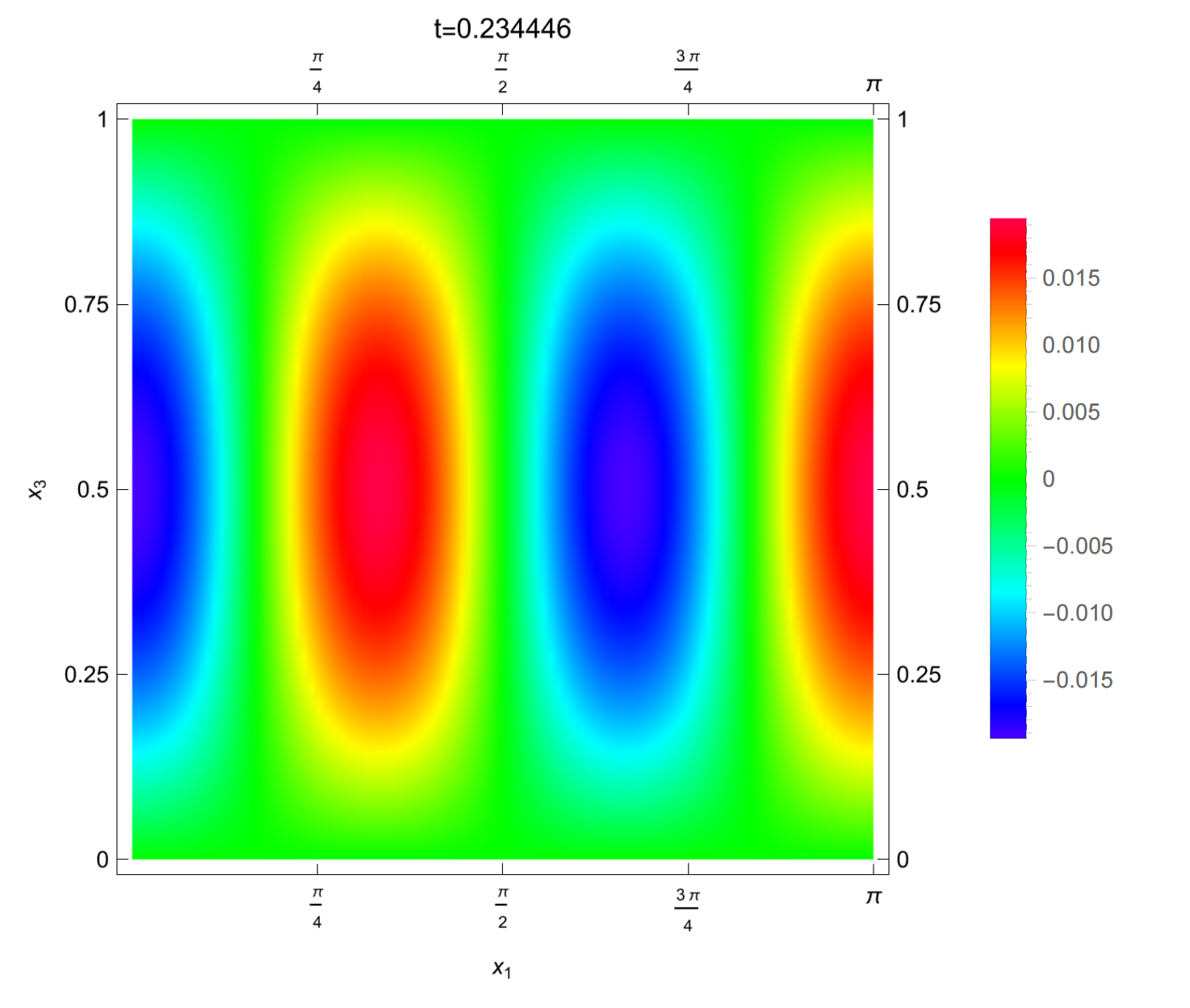}
	\end{tabular}
	\caption{Temperature plot in a period for parameter configuration \( \Ta = 2700\), \(\Q = 500\), \(\Pr = 0.3\), \(L_1 =1\), \(L_2 = 1.2 \). The approximate period is \( 0.312595 \).}
    \label{Ta2700Q500Pr01L11L212thermal}
\end{figure}

  Setting \( (\Ta, \Q, \Pr, L_1,L_2) = (2700,500,0.3,1,1.2)\), we have \( \Ra_c=\Ra_{c_2}  \approx 2350.94 \),  critical index \( J_1=(j_1,k_1,l_1) = (3,0,1) \), and \(a(Ra)= a(\Ra_{c_2})  \approx  -0.0480634\). Namely, the system \eqref{abstract1} undergoes a continuous transition from $(\Psi,\Ra)=(0,2350.94 )$.
After the continuous transition, the new states $\Psi_{p}$
is a periodic solution, given by 
        \begin{align*}
            \begin{aligned}
                &\Psi_p=\left(\frac{4\sigma}{-\pi a(Ra_{c_{2}})}\right)^{\frac{1}{2}}\sin(\rho t)\Psi_{J_{1}}^{1}+\left(\frac{4\sigma}{-\pi a(Ra_{c_{2}})}\right)^{\frac{1}{2}}\cos(\rho t)\Psi_{J_{1}}^{2}+o(\sqrt{\abs{\sigma}}).
            \end{aligned}
        \end{align*}
Here, we plot $\Psi_{p}$ at $\Ra=2351.94$ by using its leading term, see \autoref{Ta2700Q500Pr01L11L212stream}-\autoref{Ta2700Q500Pr01L11L212thermal}.

For $( \Ta, \Pr )\in \text{\MakeUppercase{\romannumeral3}}$, there exists transition from
 real eigenvalue with multiplicity two. Although in theory there are eight
  scenarios involved the transitions, only the five of which given in \autoref{double-real-theorem} are realized
  as we show below through five group of specified values of $(\Ta,\Q,\Pr,L_1,L_2)$,
  see \autoref{tab1}. In fact, the $k$-group of specified value of $(\Ta,\Q,\Pr,L_1,L_2)$ given in the $k$-row of \autoref{tab1} are these parameters at which the condition for the $k$-scenario of \autoref{double-real-theorem} is satisfied.

   \begin{table*}[htbp]
    \centering
    \begin{tabular}{|c|c|c|c|}
        \hline
               \cline{1-4}
       Specified parameter &Critical index&Critical value&Transition number\\
       $ (\Ta,\Q,\Pr,L_1,L_2)$&$(j,k,l)$&$\Ra_c=\Ra_{c_1}$&$(\Gamma_{1},\Gamma_{2},\Gamma_{3})$\\
        \cline{1-4}
        $(2800,500,0.75,1,1.2)$&$(4,\pm1,1)$&2784.75&(-0.004,-0.3,-4)\\
        \cline{1-4}
         $(2700,500,0.9,1,1.2)$&$(4,\pm1,1)$&2745.12&(-0.03,0.1,-2)\\
        \cline{1-4}
         $(2000,100,0.6,1,1.5)$&$(4,\pm1,1)$&2307.09&(2.07,-9.6,-1.9)\\
        \cline{1-4}
         $(2700,100,0.9,1,1.2)$&$(4,\pm1,1)$&2686.52&(0.8,1.5,-10.3)\\
        \cline{1-4}
         $(2900,100,0.7,1,1.2)$&$(4,\pm1,1)$&2794.17& (1.2,-0.3,-3.7)\\
        \hline
    \end{tabular}
    \caption{Five group of specified values of control parameters
    }\label{tab1}
\end{table*}

\begin{figure}[H]
    \begin{minipage}[t]{0.45\linewidth}
        \centering
        {\includegraphics[width=1.5in]{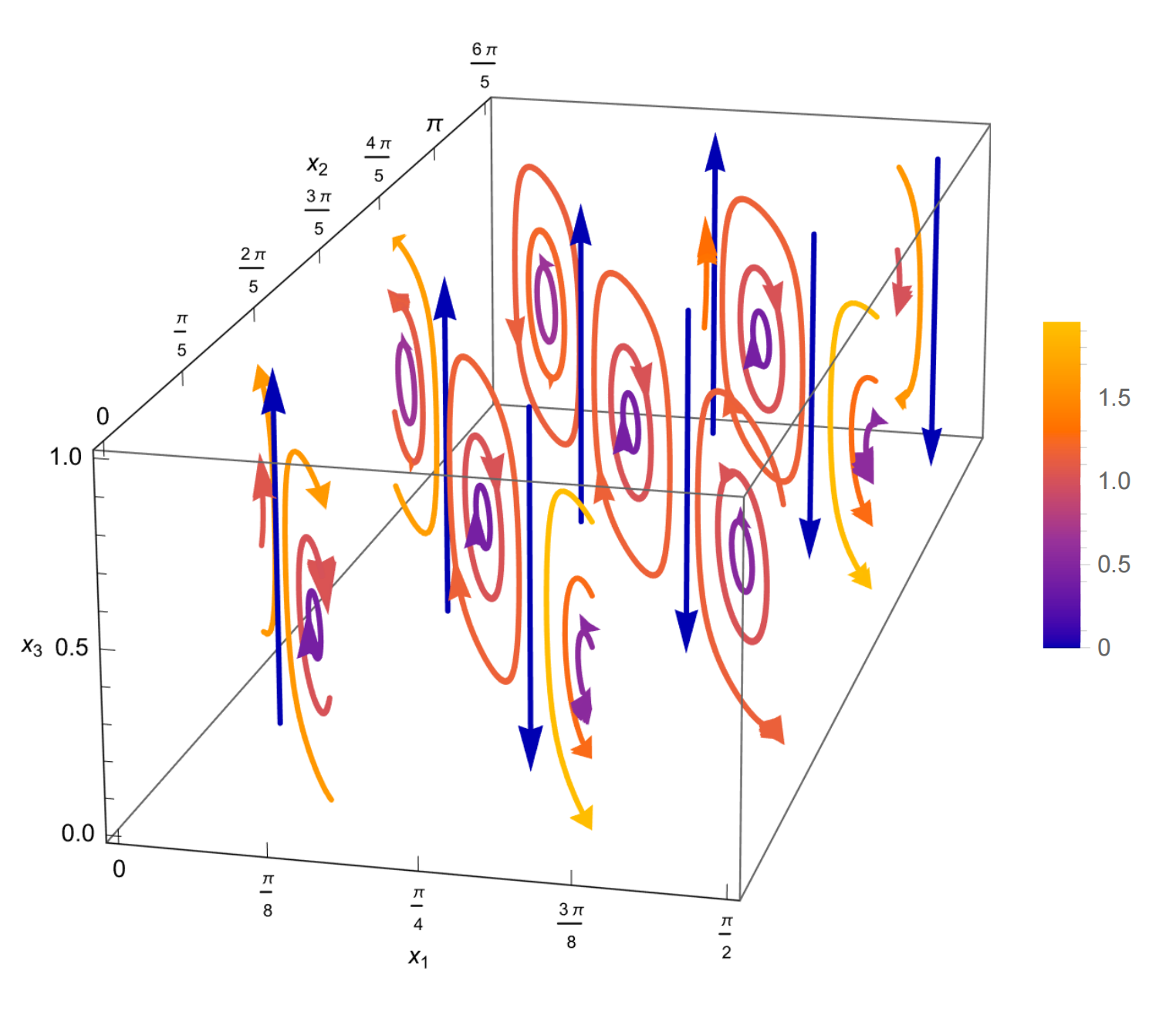}}
    \end{minipage}
    \hfill
    \begin{minipage}[t]{0.45\linewidth}
        \centering
        {\includegraphics[width=1.5in]{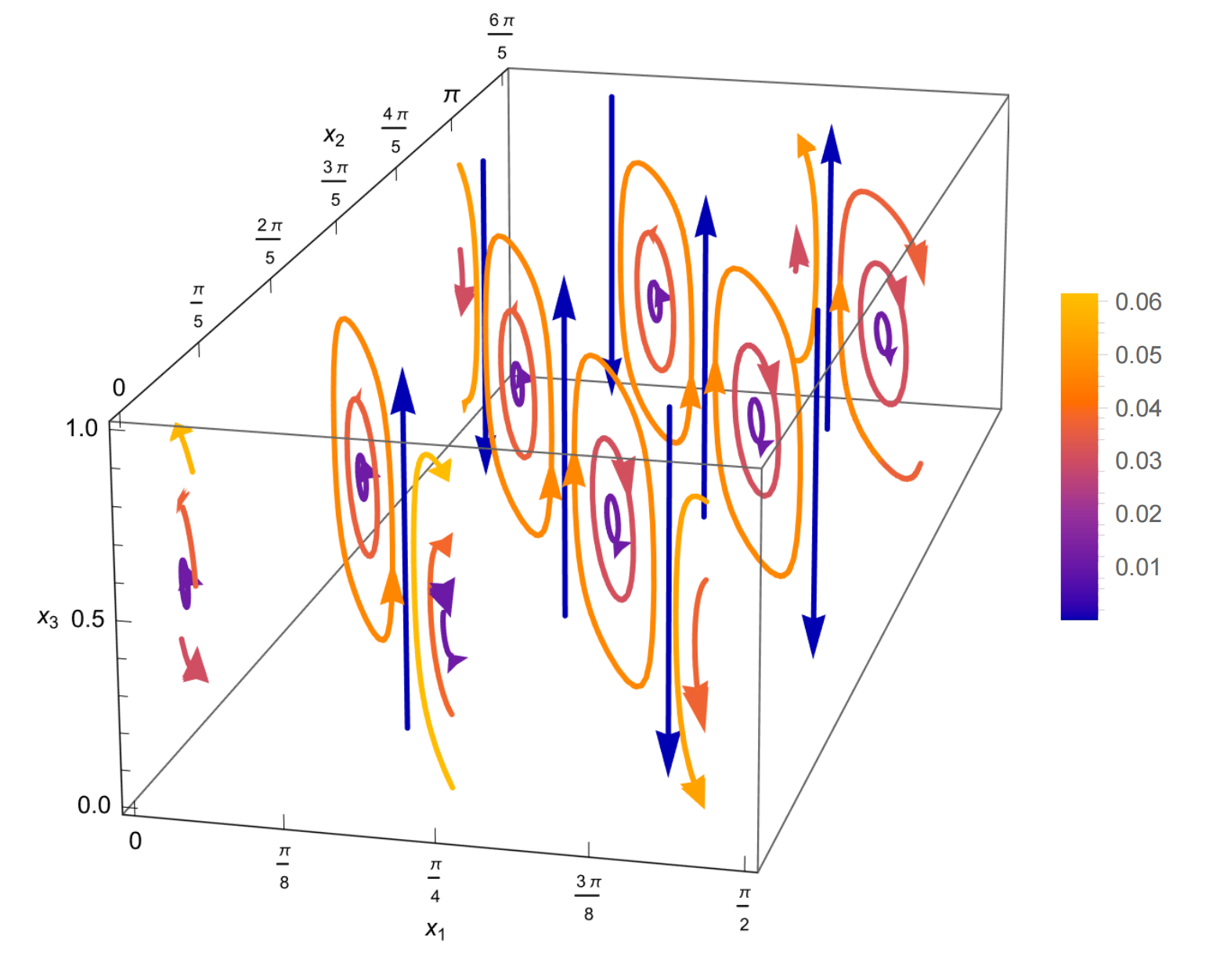}}
    \end{minipage}
    \caption{The approximate bifurcated stable solution \( \Psi_{1} \) – stream (left) and magnetic induction (right). In both figures, half of the period in \( x_2 \) is drawn, where \( \Ta = 2800\), \(\Q = 500\), \(\Pr = 0.75 \), \(L_1 =1\), \(L_2 = 1.2 \) and \( \Ra =2785.75> \Ra_{c_1} \).}
    \label{Ta2800Q500Pr075L11L212Psi1velocityandmag}
\end{figure}

\begin{figure}[H]
    \centering
    \includegraphics[height=1.5in]{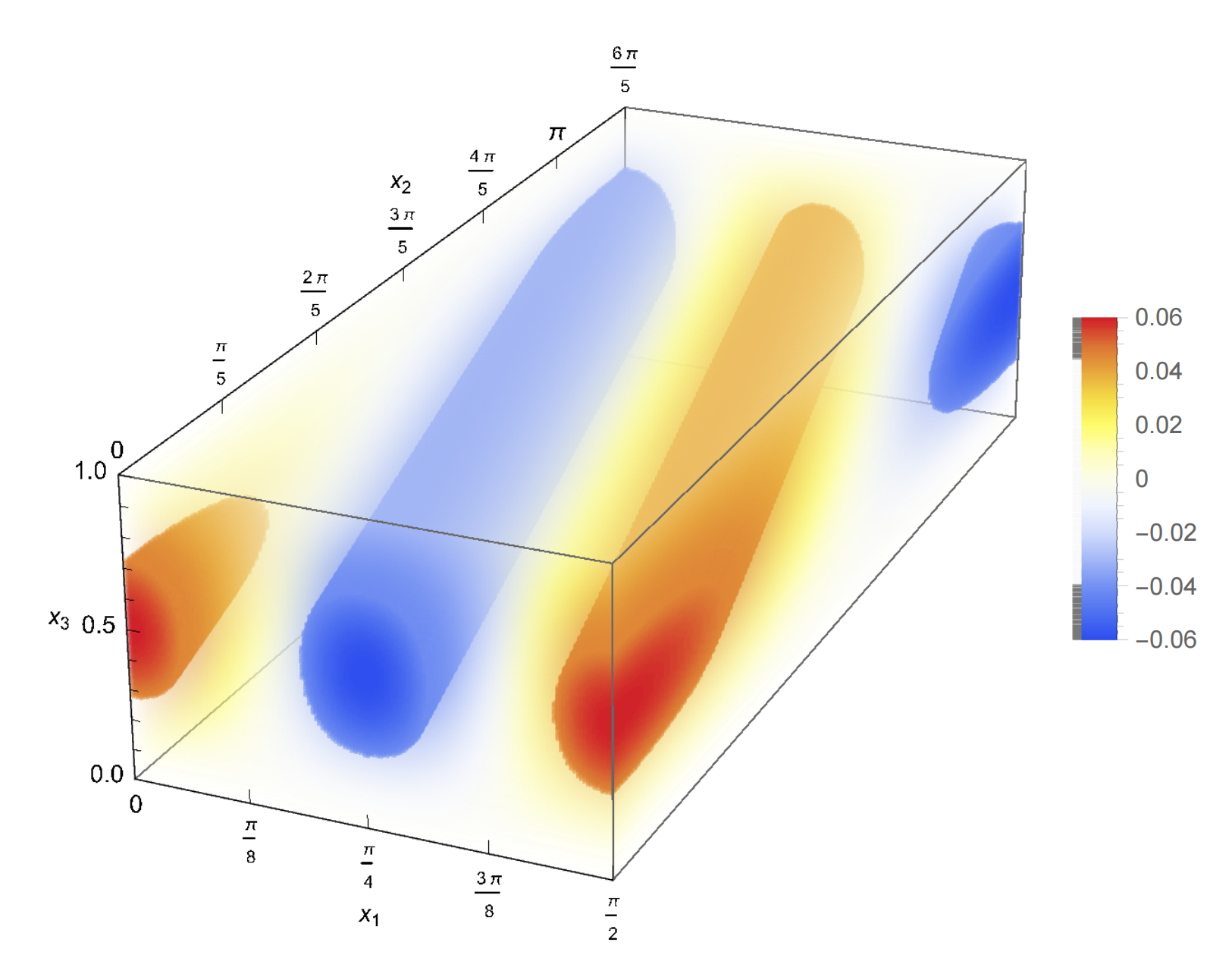}
    \caption{Temperature of the approximate bifurcated solution \( \Psi_{1} \) (half period in \( x_2 \)), where \( \Ta = 2800\), \(\Q = 500\), \(\Pr = 0.75 \), \(L_1 =1\), \(L_2 = 1.2 \) and \( \Ra = 2785.75>\Ra_{c_1} \).}
    \label{Ta2800Q500Pr075L11L212Psi1}
\end{figure}

\begin{figure}[H]
    \begin{minipage}[t]{0.45\linewidth}
        \centering
        {\includegraphics[width=1.5in]{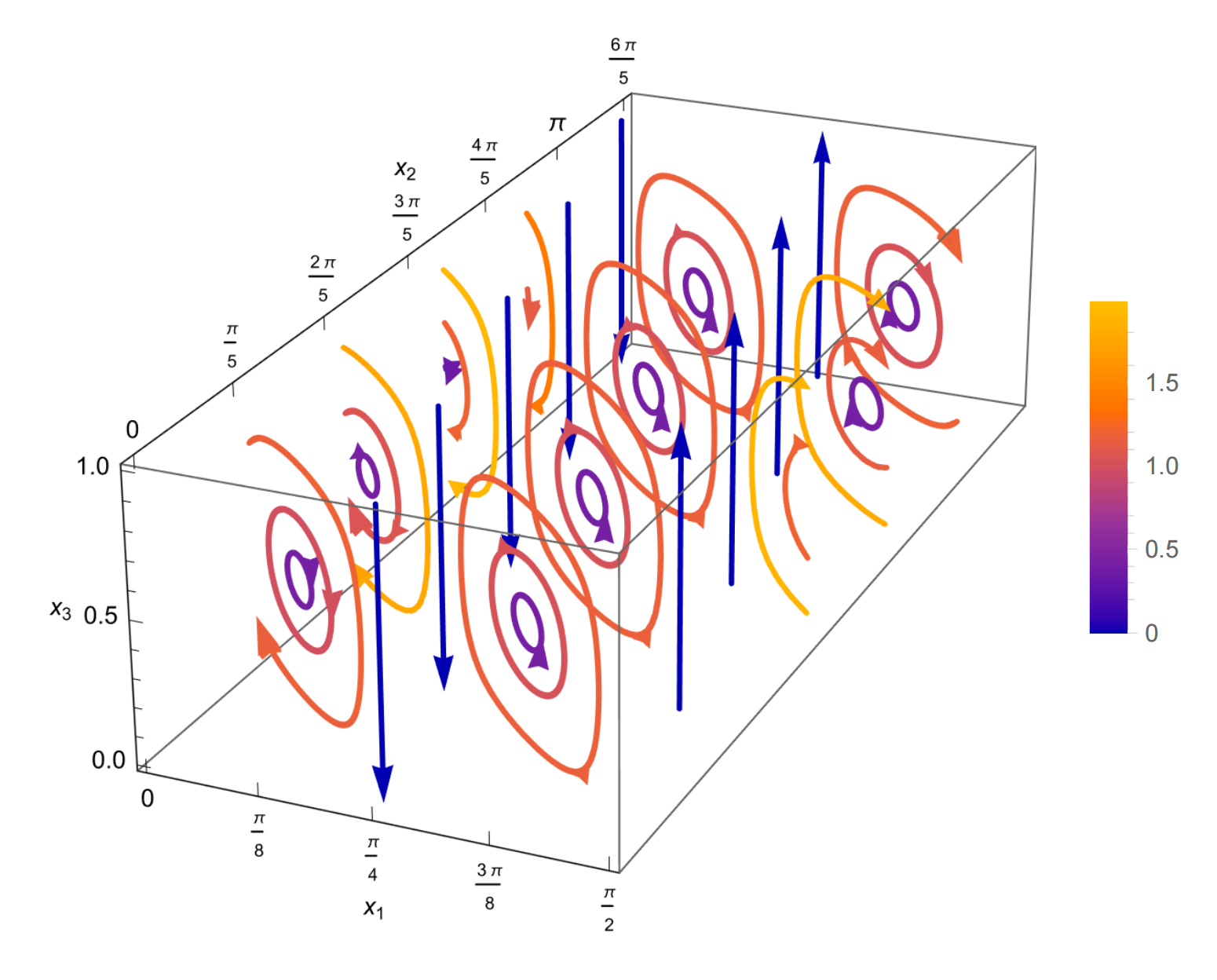}}
    \end{minipage}
    \hfill
    \begin{minipage}[t]{0.45\linewidth}
        \centering
        {\includegraphics[width=1.5in]{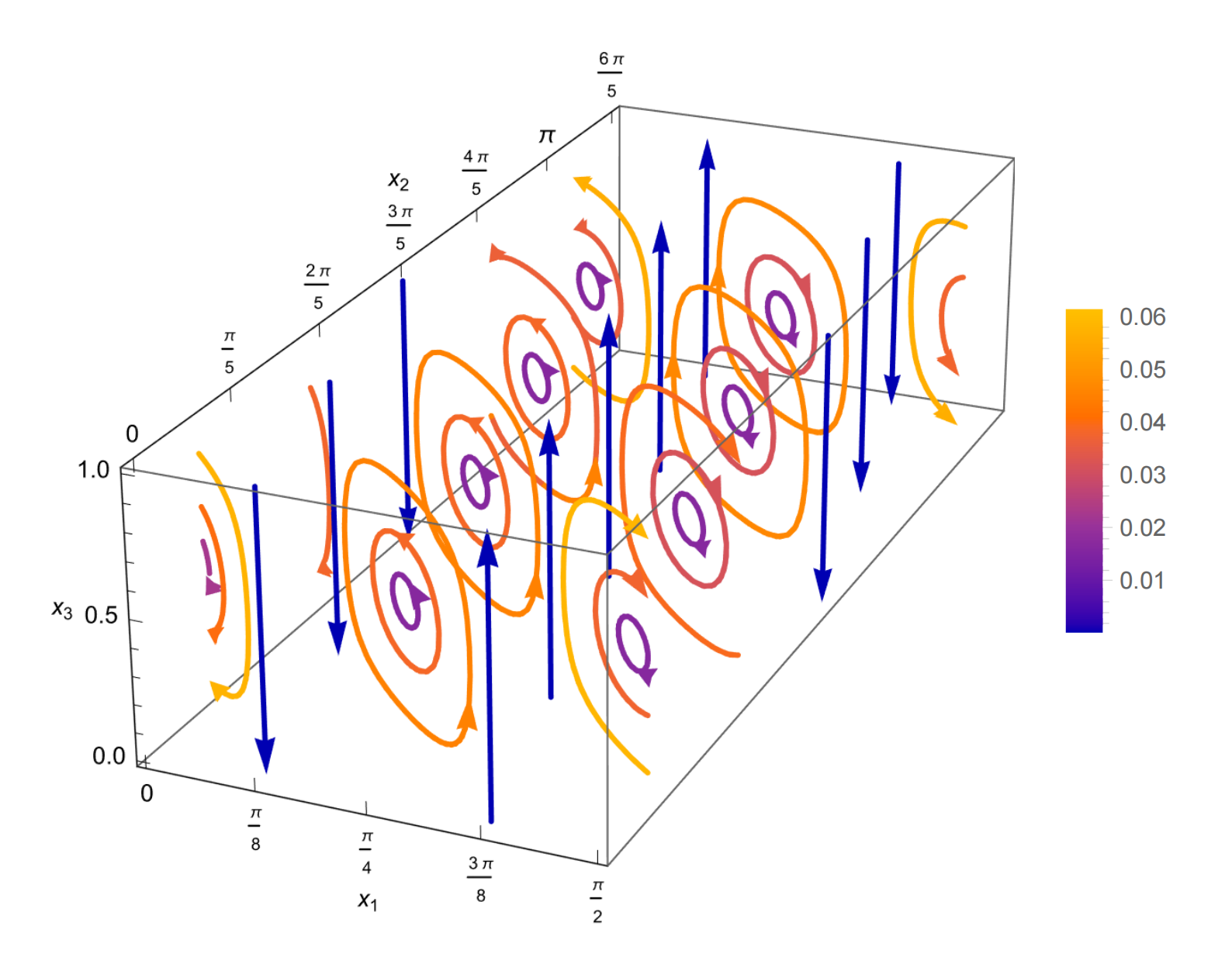}}
    \end{minipage}
    \caption{The approximate bifurcated stable solution \( \Psi_{3} \) – stream (left) and magnetic induction (right). In both figures, half of the period in \( x_2 \) is drawn, where \( \Ta = 2800\), \(\Q = 500\), \(\Pr = 0.75 \), \(L_1 =1\), \(L_2 = 1.2 \) and \( \Ra =2785.75> \Ra_{c_1}\).}
    \label{Ta2800Q500Pr075L11L212Psi3velocityandmag}
\end{figure}

\begin{figure}[H]
    \centering
    \includegraphics[height=1.5in]{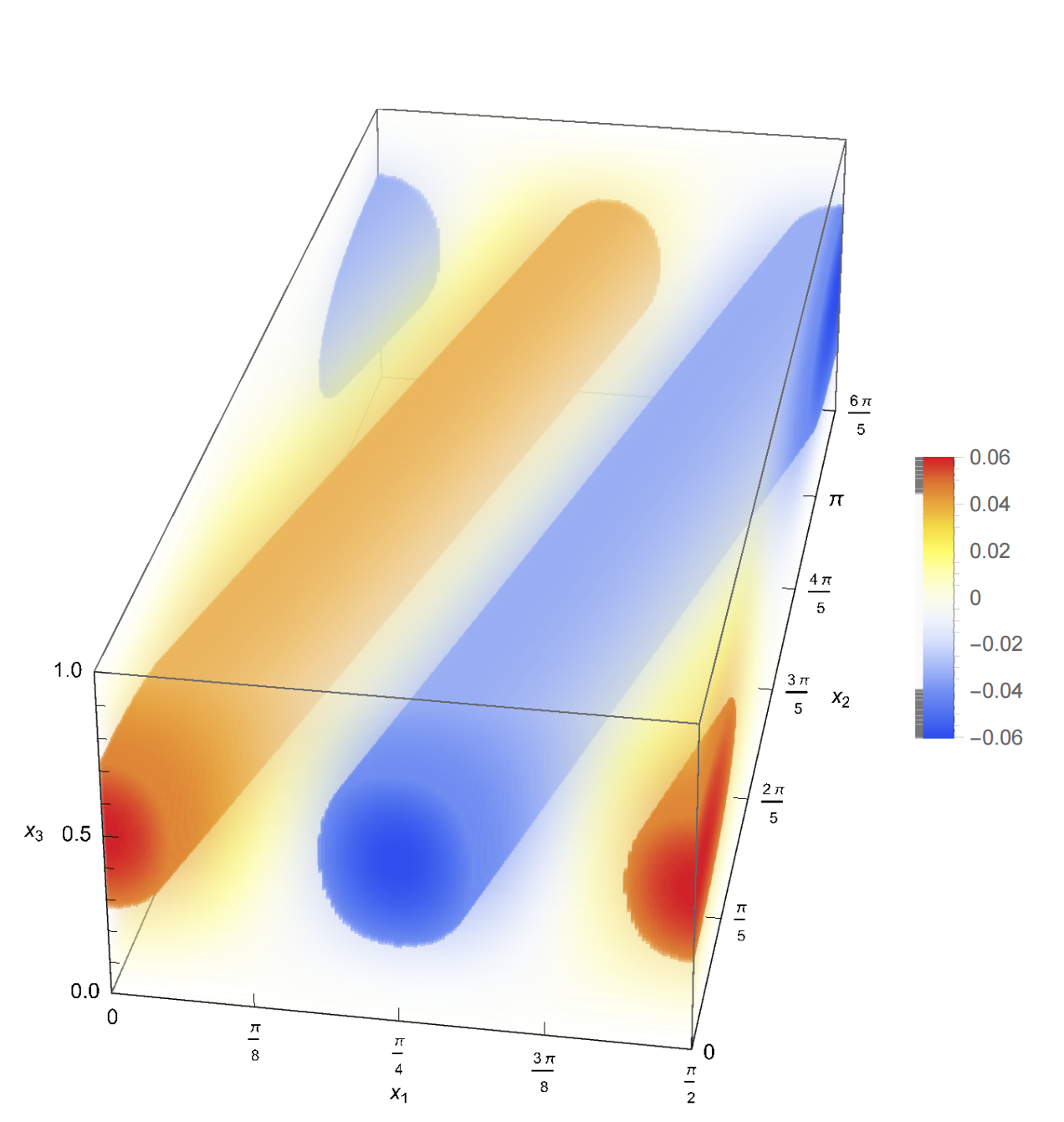}
    \caption{Temperature of the approximate bifurcated solution \( \Psi_{3} \) (half period in \( x_2 \)), where \( \Ta = 2800\), \(\Q = 500\), \(\Pr = 0.75 \), \(L_1 =1\), \(L_2 = 1.2 \) and \( \Ra =2785.75> \Ra_{c_1} \).}
    \label{Ta2800Q500Pr075L11L212Psi3}
\end{figure}

\begin{figure}[H]
    \begin{minipage}[t]{0.45\linewidth}
        \centering
        {\includegraphics[width=1.5in]{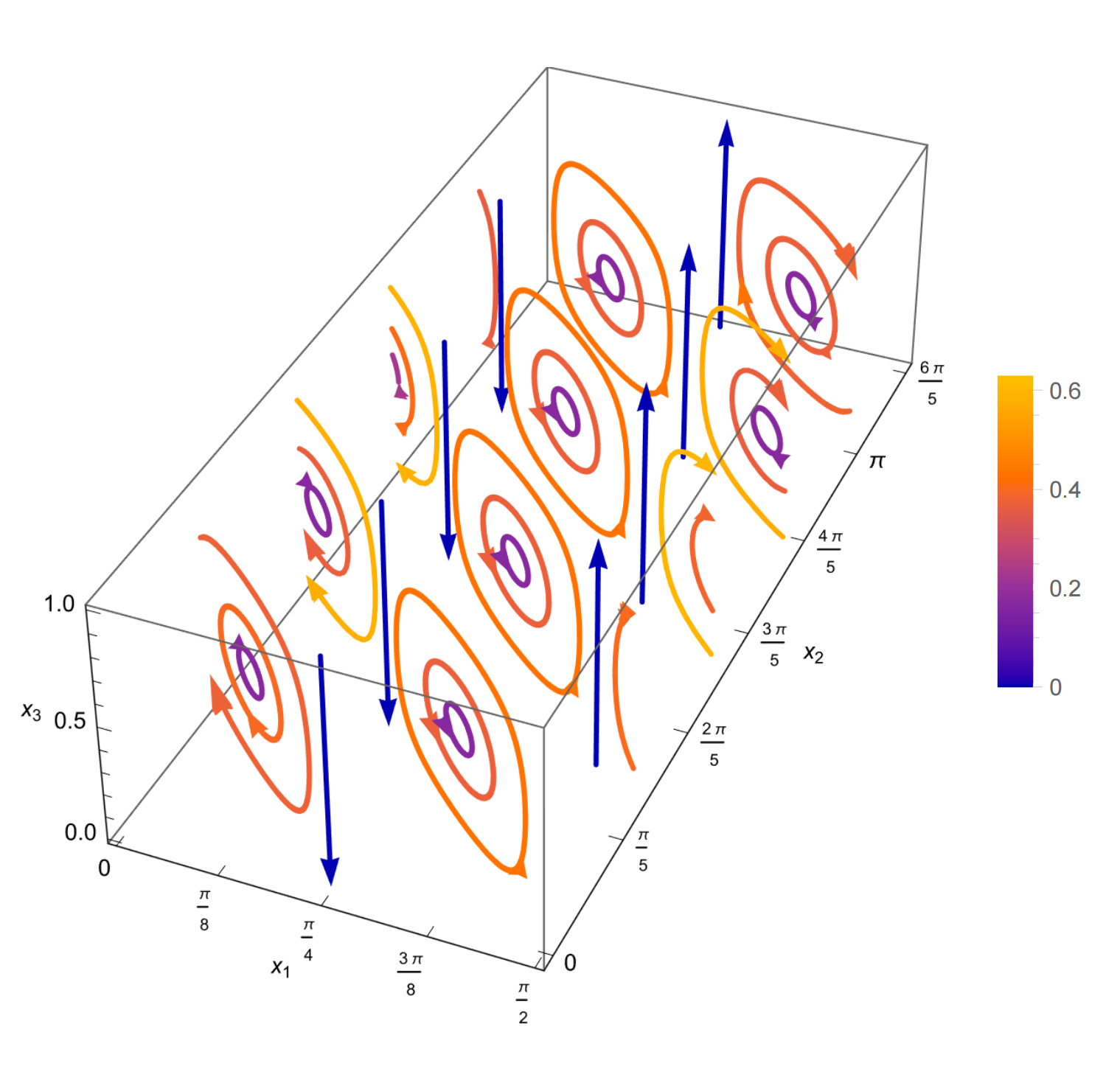}}
    \end{minipage}
    \hfill
    \begin{minipage}[t]{0.45\linewidth}
        \centering
        {\includegraphics[width=1.5in]{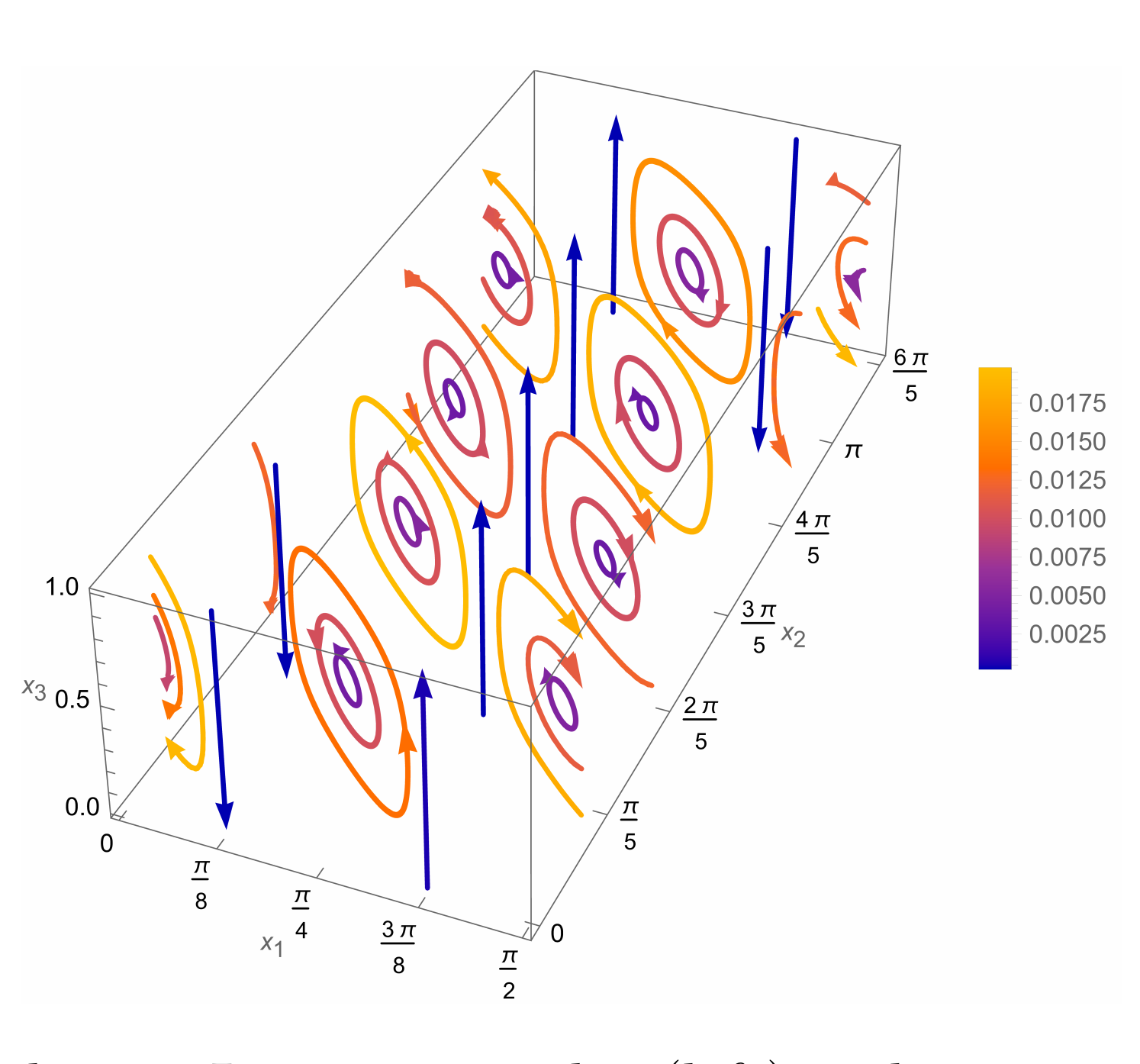}}
    \end{minipage}
    \caption{The approximate bifurcated stable solution \( \Psi_{3} \) – stream (left) and magnetic induction (right). In both figures, half of the period in \( x_2 \) is drawn, where \( \Ta = 2700\), \(\Q = 500\), \(\Pr = 0.9 \), \(L_1 =1\), \(L_2 = 1.2 \) and \( \Ra = 2746.12>\Ra_{c_1} \).}
    \label{Ta2000Q100Pr06L11L215velocityandmag}
\end{figure}

\begin{figure}[H]
    \centering
    \includegraphics[height=1.5in]{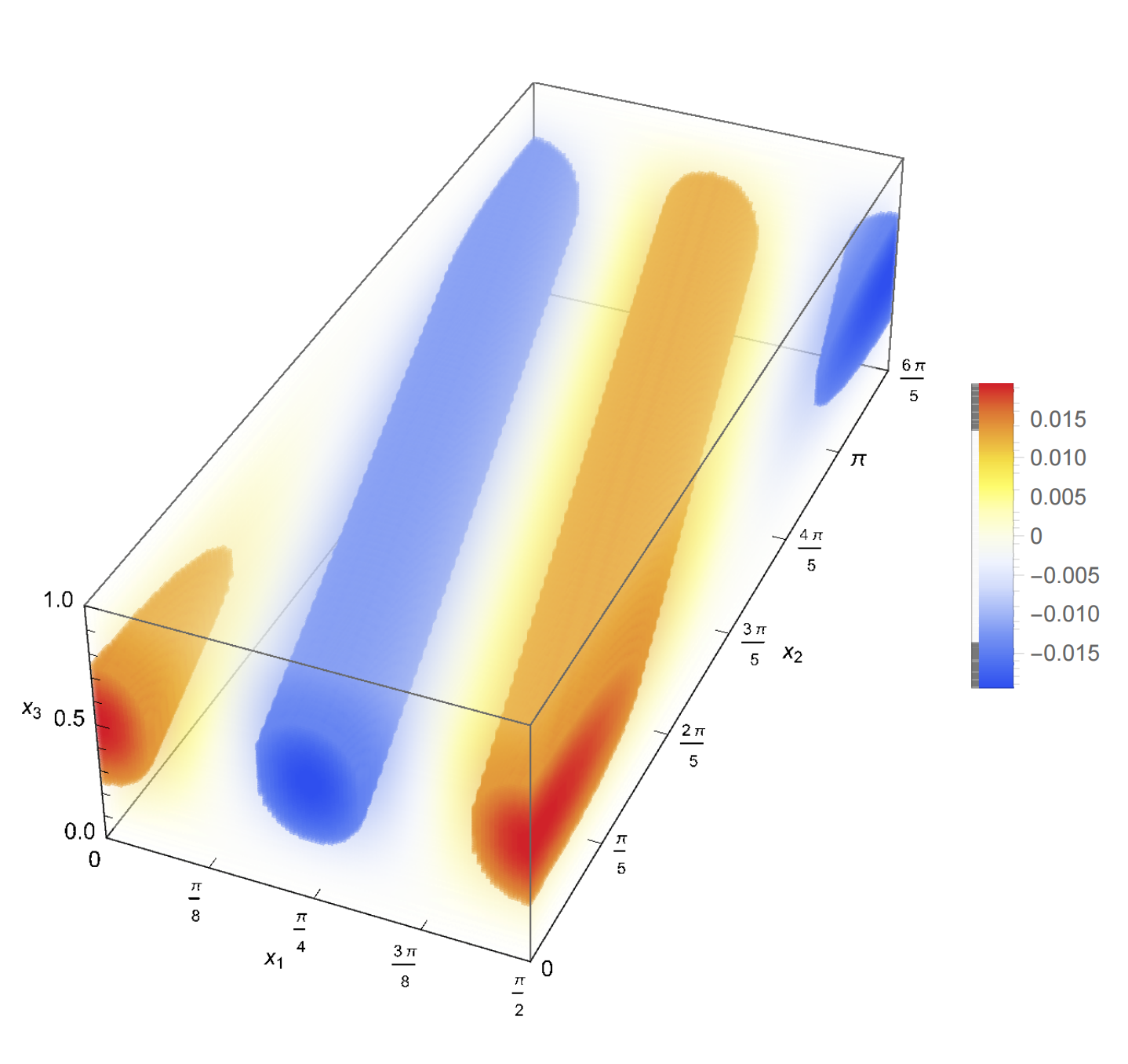}
    \caption{Temperature of the approximate bifurcated solution \( \Psi_{3} \) (half period in \( x_2 \)), where \( \Ta = 2700\), \(\Q = 500\), \(\Pr = 0.9 \), \(L_1 =1\), \(L_2 = 1.2 \) and \( \Ra =2746.12>\Ra_{c_1} \).}
    \label{Ta2000Q100Pr06L11L215thermal}
\end{figure}

From the first two
 scenarios of
\autoref{double-real-theorem}, we have known that after a continuous 
transition where a real eigenvalue with multiplicity two becomes critical,
 the system \eqref{abstract1} bifurcates on $\Ra>\Ra_{c_1}$ to a local attractor
  $\mathcal{A}$ which contains several non-degenerate equilibrium points.
  Among these equilibrium points, stable points are potential states
  after the continuous transition. For the purpose of illustrating 
  these potential states, we plot \( \Psi_1 \) and \( \Psi_3 \) 
 with $(\Ta,\Q,\Pr,L_1,L_2)$ given in the first 
  row of \autoref{tab1} and $\Ra=2785.75$, see \autoref{Ta2800Q500Pr075L11L212Psi1velocityandmag} - \autoref{Ta2800Q500Pr075L11L212Psi3}. 
We also plot \( \Psi_3 \)
with $(\Ta,\Q,\Pr,L_1,L_2)$ given in the second 
  row of \autoref{tab1} and $\Ra=2746.12$, see
\autoref{Ta2000Q100Pr06L11L215velocityandmag} and \autoref{Ta2000Q100Pr06L11L215thermal}.



\section{Conclusions}\label{section7}
In this article, we study a rotating magnetic convective model from the perspective of phase transition dynamics \cite{Ma1,Ma2}. 
First, we verify the PES condition for the RMC model \eqref{M2} ,
showing that it will undergo a dynamic transition at some critical control parameter.
To determine the type of the
transition and search for the potential states after which, we
obtain the explicit expression of its critical control parameter and
the multiplicity of first eigenvalue. Then, we establish several nonlinear transition theorems with one or several explicit transition numbers by using the center manifold reduction,
 which provide the detailed description on the dynamic types of 
 transition and bifurcation.
 
 To show our theoretical results on the transition types for the 
 RMC model \eqref{M2} at certain parameter region,  we perform some careful numerical evaluations of these transition numbers. Our numerical results show the model not only undergoes a continuous type transition,  but also a jump type transition occur at certain parameter regime, which is a new phenomenon not present in the classical RB convection \cite{HHW2019}. Hence, rotation and magnetic field have a substantial impact on the dynamic transition of RB convection.
 
We find that the number of nonzero stable states after a continuous transition
dependent on the type of first eigenvalue, and values of transition numbers.
 If first eigenvalue is real and simple, there are two nontrivial steady-states 
 which are potential stable states.
  If first eigenvalue is complex, there is a periodic solution
 which is the unique potential stable state.
  If first eigenvalue is real and has multiplicity two,  there are four or eight nontrivial steady-states which are potential stable states, depending on $\Gamma_1-\Gamma_3$.

Here, we roughly discuss the physical meanings of continuous and jump types. According to the analysis of spectrum, when $\Ra$ is less than critical value, the trivial steady-state of model is always stable, which is the state that we observe.
By increasing the control parameter $\Ra$, a continuous transition means when $\Ra$ crosses the critical value, the trivial state of the system is smoothly and continuously replaced by a nonzero stable state. The nonzero stable state can be computed approximately in the vicinity of the critical control parameter, and as $\Ra$ approaches critical value, the nonzero stable state will go to the zero equilibrium. 
In contrast, a jump type means the system jump to a nonzero state abruptly,
 which is not in the vicinity of the zero equilibrium. 

 In the present article,  we only consider the transitions from first real simple eigenvalue, first complex eigenvalues and first real eigenvalue with multiplicity two. According to our numerical results, see \autoref{Ta2700Q100Pr06L1052L2052}, there are many other cases such as transition from two pair of complex conjugate eigenvalues or a pair of complex eigenvalues plus simple real eigenvalue, which are non-generic. We will consider them in the future.


\end{document}